\documentclass{article} 
\usepackage{etoolbox}
\newtoggle{springer}
\togglefalse{springer}
\usepackage{graphicx}
\usepackage{algorithm}
\usepackage{algorithmic}
\usepackage{hyperref}
\usepackage{color}
\usepackage{amsmath}
\usepackage{rotating}

\iftoggle{springer}
{
}
{
\usepackage{natbib}
\usepackage{amsthm}
}

\newcommand{\figSize}{0.32\textwidth}
\newcommand{\BEAS}{\begin{eqnarray*}}
\newcommand{\EEAS}{\end{eqnarray*}}
\newcommand{\BEA}{\begin{eqnarray}}
\newcommand{\EEA}{\end{eqnarray}}
\newcommand{\BEQ}{\begin{equation}}
\newcommand{\EEQ}{\end{equation}}
\newcommand{\BIT}{\begin{itemize}}
\newcommand{\EIT}{\end{itemize}}
\newcommand{\BNUM}{\begin{enumerate}}
\newcommand{\ENUM}{\end{enumerate}}
\newcommand{\BA}{\begin{array}}
\newcommand{\EA}{\end{array}}
\newcommand{\diag}{\mathop{\rm diag}}
\newcommand{\Diag}{\mathop{\rm Diag}}

\newcommand{\idm}{I}
\newcommand{\rb}{\mathbb{R}}
\newcommand{\R}{\mathbb{R}}

\input{macros}

\def \E{{\mathbb E}}
\def \P{{\mathbb P}}
\def \R{{\mathbb R}}

\def \F{{\mathcal F}}

\definecolor{blu}{rgb}{0,0,1} 
\def\blu#1{{\color{blu}#1}}

\iftoggle{springer}
{
\smartqed  
\newcommand{\citep}{\cite} 
\newcommand{\citet}{\cite} 
\journalname{Math.~Program.~Ser.~A}
\bibliographystyle{spmpsci}      
}
{
\oddsidemargin .25in  
\evensidemargin .25in \marginparwidth 0.07 true in
\topmargin -0.5in \addtolength{\headsep}{0.25in}
\textheight 8.5 true in    
\textwidth 6.0 true in    
\widowpenalty=10000 \clubpenalty=10000
\parindent 0pt
\topsep 4pt plus 1pt minus 2pt
\partopsep 1pt plus 0.5pt minus 0.5pt
\itemsep 2pt plus 1pt minus 0.5pt
\parsep 2pt plus 1pt minus 0.5pt
\parskip .5pc
\usepackage{natbib}
\newtheorem{theorem}{Theorem}
\newtheorem{lemma}{Lemma}
\bibliographystyle{plainnat}
}
\newcommand{\auth}[1]{\iftoggle{springer}{#1~}{}}
\newcommand{\citesee}[2]{\iftoggle{springer}{\cite[see #1]}{\citep[see][#1]}{#2}}
\newcommand{\citethm}[2]{\iftoggle{springer}{\cite[#1]}{\citep[][#1]}{#2}}
\newcommand{\smaller}{\iftoggle{springer}{\scriptsize}{}}

\begin{document}

\title{Minimizing Finite Sums with the Stochastic Average Gradient
}


\iftoggle{springer}
{
\author{Mark Schmidt \and Nicolas Le Roux \and Francis Bach}
\authorrunning{M.~Schmidt, N.~Le Roux, F.~Bach} 
\institute{Mark Schmidt \at
              Department of Computer Science, University of British Columbia\\
              201 2366 Main Mall, Vancouver BC V6T 1Z4\\
              \email{schmidtm@cs.ubc.ca}           
           \and
           Nicolas Le Roux\at
              Criteo\\
	      32 rue Blanche, 75009 Paris\\
              \email{nicolas@le-roux.name}           
	  \and
	  Francis Bach\at
	  INRIA - SIERRA Project-Team\\
	  Laboratoire d'Informatique de l'Ecole Normale Superieure\\
	  23, avenue d'Italie, CS 81321, 75214 Paris Cedex 13\\
              \email{francis.bach@ens.fr}           
}
\date{Received: date / Accepted: date}
}
{
\author{
Mark Schmidt\\
\texttt{schmidtm@cs.ubc.ca}
\and
Nicolas Le Roux\\
\texttt{nicolas@le-roux.name}
\and
Francis Bach\\
\texttt{francis.bach@ens.fr}
\and
\\
INRIA - SIERRA Project - Team\\
D\'epartement d'Informatique de l'\'Ecole Normale Sup\'erieure\\
Paris, France
}
\date{January 19, 2015}
}

\maketitle

\begin{abstract}
We analyze the stochastic average gradient (SAG) method for optimizing the sum of a finite number of smooth convex functions. Like stochastic gradient (SG) methods, the SAG method's iteration cost is independent of the number of terms in the sum. However, by incorporating a memory of previous gradient values the SAG method achieves a faster convergence rate than black-box SG methods. The convergence rate is improved from $O(1/\sqrt{k})$ to $O(1/k)$ in general, and when the sum is strongly-convex the convergence rate is improved from the sub-linear $O(1/k)$ to a linear convergence rate of the form $O(\rho^k)$ for $\rho < 1$. Further, in many cases the convergence rate of the new method is also faster than black-box deterministic gradient methods, in terms of the number of gradient evaluations. This extends our earlier work~\citep{roux2012stochastic}, which only lead to a faster rate for well-conditioned strongly-convex problems.
Numerical experiments indicate that the new algorithm often dramatically outperforms existing SG and deterministic gradient methods, and that the performance may be further improved through the use of non-uniform sampling strategies.
\iftoggle{springer}
{
}{}
\end{abstract}

\section{Introduction}
\label{intro}

A plethora of the optimization problems arising in practice involve computing a minimizer of a finite sum of functions measuring misfit over a large number of data points. 
A classical example is least-squares regression,
\[
\minimize{x\in\Real^p}\frac{1}{n}\sum_{i=1}^n (a_i^Tx - b_i)^2,
\]
where the $a_i \in \Real^p$ and $b_i \in \Real$ are the data samples
associated with a regression problem. Another important example is logistic regression,
\[
\minimize{x\in\Real^p}\quad \frac{1}{n}\sum_{i=1}^n \log(1+\exp(-b_ia_i^\top x)),
\]
where the $a_i \in \Real^p$ and $b_i \in \{-1,1\}$ are the data samples
associated with a binary classification problem. A key challenge arising in modern
applications is that the number of data points $n$ (also known as \emph{training examples})
can be extremely large, while there is often a large amount of redundancy between examples.
 The most wildly successful class of algorithms for taking advantage 
of the \emph{sum} structure for problems where $n$ is very large are \emph{stochastic gradient} (SG) methods~\citep{robbins1951stochastic,bottou-lecun-2004}.
Although the theory behind SG methods allows them to be applied more generally, 
SG methods are often used to solve the problem of optimizing a finite sample average,
\begin{equation}
\label{eq:1}
\minimize{x\in\Real^p}\quad g(x) \defd \frac{1}{n}\sum_{i=1}^n f_i(x).
\end{equation}
In this work, we focus on such \emph{finite data} problems where
each $f_i$ is \emph{smooth} and \emph{convex}. 

In addition to this basic setting, we will also be interested in cases where the sum $g$ has the additional property that it is \emph{strongly-convex}. This often arises due to the
use of a strongly-convex regularizer such as the squared $\ell_2$-norm, resulting
in problems of the form
\begin{equation}
\label{eq:L2}
\minimize{x\in\Real^p}\quad \frac{\lambda}{2}\|x\|^2 + \frac{1}{n}\sum_{i=1}^n l_i(x),
\end{equation}
where each $l_i$ is a data-misfit function (as in least-squares and logistic regression) and  the positive scalar $\lambda$ controls the strength of the regularization. 
These 
problems can be put in the framework of~\eqref{eq:1} by using the choice
\[
f_i(x) \defd \frac{\lambda}{2}\|x\|^2 + l_i(x).
\]
The resulting function $g$ will be strongly-convex provided that the individual loss functions $l_i$ are convex.
An extensive list of convex loss functions used in a statistical data-fitting context is given
by~\auth{Teo et al.}~\citet{teo2007scalable}, and non-smooth loss functions
(or regularizers) can also be put in this framework by using smooth approximations (for example, see~\citet{nesterov2005smooth}).

For optimizing problem~\eqref{eq:1}, the standard \emph{deterministic} or \emph{full gradient} (FG) method, which dates back to \auth{Cauchy}\citet{cauchy1847methode}, uses iterations
of the form
\begin{equation}
\label{eq:FG}
x^{k+1} = x^k - \alpha_kg'(x^k) = x^k - \frac{\alpha_k}{n}\sum_{i=1}^nf_i'(x^k),
\end{equation}
where $\alpha_k$ is the step size on iteration $k$.
Assuming that a minimizer $x^\ast$ exists, then under standard assumptions the sub-optimality achieved on iteration $k$ of the FG method with a constant step size is given by
\[
g(x^k) - g(x^\ast) = O(1/k),
\]
when $g$ is convex~\citesee{Corollary~2.1.2}{nesterov2004introductory}. This results in a \emph{sublinear} convergence rate. When $g$ is strongly-convex, the error also satisfies
\[
g(x^k) - g(x^\ast) = O(\rho^k),
\]
for some $\rho < 1$ which depends on the condition number of $g$~\citesee{Theorem~2.1.5}{nesterov2004introductory}.
This results in a \emph{linear} convergence rate, which is also known as a \emph{geometric} or \emph{exponential}
rate because the error is cut by a fixed fraction on each iteration.
Unfortunately, the FG method can be unappealing
when $n$ is large because its iteration cost scales linearly in $n$.

 The main appeal of SG methods is that they have an iteration cost which is \emph{independent} of $n$, making them suited for modern problems where $n$ may be very large. The basic SG method for optimizing~\eqref{eq:1}
uses iterations of the form
\begin{equation}
\label{eq:SG}
x^{k+1} = x^k - \alpha_kf_{i_k}'(x^k),
\end{equation}
where at each iteration an index $i_k$ is sampled
uniformly from the set $\{1,\dots,n\}$.
The randomly chosen gradient $f_{i_k}'(x^k)$ yields an unbiased estimate of
the true gradient $g'(x^k)$ and one can show under standard assumptions (see~\citep{nemirovski2009robust}) that,
for a suitably chosen decreasing step-size sequence $\{\alpha_k\}$, the SG
iterations have an expected sub-optimality for convex objectives of
\[
\mathbb{E}[g(x^k)] - g(x^\ast) = O(1/\sqrt{k}),
\]
and an expected sub-optimality for strongly-convex objectives of
\[
\mathbb{E}[g(x^k)] - g(x^\ast) = O(1/k).
\]
In these rates, the expectations are taken with respect 
to the selection of the $i_k$ variables. These sublinear rates are slower than the corresponding
rates for the FG method, and under certain assumptions these convergence rates are \emph{optimal} 
in a model of computation where the algorithm only accesses
the function through unbiased measurements of its objective and gradient (see~\citet{nemirovsky1983problem,nemirovski2009robust,agarwal2010information}).
Thus, we should not expect to be able to obtain the convergence rates of the FG method 
if the algorithm only relies on unbiased gradient measurements. Nevertheless, by using the stronger
assumption that the functions are sampled from a finite dataset,
in this paper we show that we
can achieve the convergence rates of FG methods while preserving the iteration complexity
of SG methods.

The primary contribution of this work is the analysis of a new algorithm that we call the \emph{stochastic average
gradient} (SAG) method, a randomized variant of the incremental aggregated
gradient (IAG) method of \auth{Blatt et al.}\citet{blatt2008convergent}. The SAG method has the low
iteration cost of SG methods, but achieves the convergence rates stated above for the  FG
method. The SAG iterations take the form
\begin{equation}
\label{eq:SAG}
x^{k+1} = x^k - \frac{\alpha_k}{n}\sum_{i=1}^ny^k_i,
\end{equation}
where at each iteration a random index $i_k$ is selected and we
set
\begin{equation}
\label{eq:yi}
y^k_i = \begin{cases}
f_i'(x^k) & \textrm{if $i = i_k$,}\\
y^{k-1}_i & \textrm{otherwise.}
\end{cases}
\end{equation}
That is, like the FG method, the step incorporates a gradient with respect
to
each function. But, like the SG method, each iteration only computes
the gradient with respect to a single example and the cost of the
iterations is independent of $n$. Despite the low cost of the SAG
iterations, we show in this paper that with a constant step-size \emph{the SAG iterations have an $O(1/k)$ convergence rate for convex objectives and a linear convergence rate for strongly-convex objectives}, like the FG method.
That is, by having access to
$i_k$ and by keeping a \emph{memory} of the most recent gradient value
computed for each index $i$, this iteration achieves a faster
convergence rate than is possible for standard SG methods.
Further, in terms of effective passes through the data, we will also see that
for many problems the convergence rate of the SAG method is also faster than is possible
for standard FG methods.

One of the main contexts where minimizing the sum of smooth convex functions arises is machine learning. In this context, $g$ is often an \emph{empirical} risk (or a regularized empirical risk), which is a sample average approximation to the \emph{true} risk that we are interested in. It is known that with $n$ training examples the empirical risk minimizer (ERM) has an error for the true risk of $O(1/\sqrt{n})$ in the convex case and $O(1/n)$ in the strongly-convex case. Since these rates are achieved by doing one pass through the data with an SG method, in the worst case the SAG algorithm applied to the empirical risk cannot improve the convergence rate in terms of the true risk over this simple method. Nevertheless, \auth{Srebro \& Sridharan}\citet{srebro2011erm} note that ``overwhelming empirical evidence shows that for almost all actual data, the ERM \emph{is} better. However, we have no understanding of why this happens". Although our analysis does not give insight into the better performance of ERM, our analysis shows that the SAG algorithm will be preferable to SG methods for finding the ERM and hence for many machine learning applications.

The next section reviews several closely-related algorithms from the
literature, including previous attempts to combine the appealing aspects of
FG and SG methods. However, despite $60$ years of extensive research on SG methods,
with a significant portion of the applications focusing on finite datasets,
we believe that this is the first general method that achieves the convergence rates of FG methods while preserving the
iteration cost of standard SG methods. Section~\ref{convergence} states the
(standard) assumptions underlying our analysis and gives our convergence rate
results.
Section~\ref{sec:implementation} discusses practical
implementation issues including how we adaptively set the step size and how we can reduce the storage cost needed by the algorithm. For example, we can reduce the memory requirements from $O(np)$ to $O(n)$ in the common scenario where each $f_i$ only depends on a linear function of $x$, as in least-squares and logistic regression.
Section~\ref{experiments} presents a numerical
comparison of an implementation based on SAG to competitive SG and FG methods, indicating that the method may be very useful for problems where we can only afford to do a few passes through
a data set.

A preliminary conference version of this work appears in~\citet{roux2012stochastic}, and we extend this work in various ways. Most notably, the analysis in the prior work focuses only on showing linear convergence rates in the strongly-convex case while the present work also gives an $O(1/k)$ convergence rate for the general convex case. In the prior work we show (Proposition~1) that a small step-size 
gives a slow linear convergence rate (comparable to the rate of FG methods in terms of effective passes through the data), while we also show (Proposition~2)
that a much larger step-size yields a much faster convergence rate, but this requires that $n$ is sufficiently large compared to the condition number of the problem.
In the present work (Section~\ref{convergence}) our analysis yields a very fast convergence rate using a large step-size (Theorem~\ref{thm}), even when this condition required by the prior work is not satisfied.
Surprisingly, for ill-conditioned problems our new analysis shows that using SAG iterations can be nearly $n$ times as fast as the standard gradient method.
To prove this stronger result, Theorem~1 employs a Lyapunov function that generalizes the Lyapunov functions used in Propositions~1 and~2 of the previous work. This new Lyapunov function leads to a unified proof for both the convex and the strongly-convex cases, and for both well-conditioned and ill-conditioned problems. However, this more general Lyapunov function leads to a more complicated analysis. To significantly simplify the formal proof, we use a computed-aided strategy to verify the non-negativity of certain polynomials that arise in the proof.
Beyond this significantly strengthened result, in this work we also argue that yet-faster convergence rates may be achieved by \emph{non-uniform} sampling (Section~\ref{sec:lipschitz}) and present numerical results showing that this can lead to drastically improved performance (Section~\ref{exp:lipschitz}).

\iftoggle{springer}
{
Due to space restrictions, some details are omitted in this article. This notably includes the proof of the main theorem, some additional experimental results, and a thorough discussion of the many interesting works that have followed after~\citet{roux2012stochastic}. These extra materials are made available in the extrended arXiv version of the paper located at:\\ \url{http://arxiv.org/abs/1309.2388}.
}

\section{Related Work}
\label{related}

There are a large variety of approaches available to accelerate the
convergence of SG methods, and a full review of this immense literature
would be outside the scope of this work. Below, we comment on the
relationships between the new method and several of the most closely-related
ideas.

\textbf{Momentum}: SG methods that incorporate a momentum term use
iterations of the form
\[
x^{k+1} = x^k - \alpha_k f_{i_k}'(x^k) + \beta_k(x^k - x^{k-1}),
\]
see~\citet{tseng1998incremental}.
It is common to set all $\beta_k = \beta$ for some constant $\beta$, and in
this case we can rewrite the SG with momentum method as
\[
\textstyle
x^{k+1} = x^k - \sum_{j=1}^k\alpha_j\beta^{k-j}f_{i_j}'(x^j).
\]
We can re-write the SAG updates~\eqref{eq:SAG} in a similar form as
\begin{equation}
\label{eq:SAG2}
\textstyle
x^{k+1} = x^k - \sum_{j=1}^k\alpha_kS(j,i_{1:k})f_{i_j}'(x^j),
\end{equation}
where the selection function $S(j,i_{1:k})$ is equal to $1/n$ if $j$ is the maximum iteration number where example $i_j$ was selected and is set to $0$ otherwise.
Thus, momentum uses a \emph{geometric weighting} of previous gradients while
the SAG iterations \emph{select and average} the most recent evaluation of each previous
gradient. While momentum can lead to improved practical performance, it
still requires the use of a decreasing sequence of step sizes and is not
known to lead to a faster convergence rate.

\textbf{Gradient Averaging}:
Closely related to momentum is using the sample average of all previous
gradients,
\[
\textstyle
x^{k+1} = x^k - \frac{\alpha_k}{k}\sum_{j=1}^kf_{i_j}'(x_j),
\]
which is similar to the SAG iteration in the form~\eqref{eq:SAG} but where \emph{all} previous
gradients are used.
This approach is used in the dual averaging method of \auth{Nesterov}\citet{nesterov2009primal} and, while this averaging procedure and its variants
lead to convergence for a constant step size and
can improve the constants in the convergence rate~\citep{xiao2010dual}, it does not improve on the sublinear convergence rates for SG methods.

\textbf{Iterate Averaging}: Rather than averaging the gradients, some authors propose to perform the basic SG iteration but use an average over certain $x^k$ values as the final estimator.
With a suitable choice of step-sizes, this gives the same asymptotic efficiency as Newton-like second-order SG methods and also leads to increased robustness of the convergence rate to the exact sequence of step sizes~\citep{polyak1992acceleration,bach2011non}. Bather's method~\citep[\S1.3.4]{kushner1997stochastic} combines gradient averaging with online iterate averaging and also displays appealing asymptotic properties. 
Several authors have recently shown that suitable iterate averaging schemes obtain an $O(1/k)$ rate for strongly-convex optimization even for non-smooth objectives~\citep{hazan2010beyond,rakhlin2012making}.  However, none of these methods improve on the $O(1/\sqrt{k})$ and $O(1/k)$ rates for SG methods. 

\textbf{Stochastic versions of FG methods}: Various options are available to accelerate the convergence of the FG method for smooth functions, such as the accelerated full gradient (AFG) method of \auth{Nesterov}\citet{nesterov1983method}, as well as classical techniques based on quadratic approximations such as diagonally-scaled FG methods, non-linear conjugate gradient, quasi-Newton, and Hessian-free Newton methods (see~\citep{nocedal2006numerical}).
There has been a substantial amount of work on developing stochastic variants of these algorithms, with several of the notable recent examples including~\citet{bordes2009sgd,kwok2009asg,sunehag2009variable,ghadimi2010optimal,martens2010deep,xiao2010dual}. \auth{Duchi et al.}\citet{duchi2010adaptive} have recently shown an improved regret bound using a diagonal scaling that takes into account previous gradient magnitudes. Alternately, if we split the convergence rate into a deterministic and stochastic part, these methods can improve the dependency of the convergence rate of the deterministic part~\citep{kwok2009asg,ghadimi2010optimal,xiao2010dual}. However, we are not aware of any existing method of this flavor that improves on the $O(1/\sqrt{k})$ and $O(1/k)$ dependencies on the stochastic part. Further, many of these methods typically require carefully setting parameters (beyond the step size) and often aren't able to take
advantage of sparsity in the gradients $f_i'$.

\textbf{Constant step size}: If the SG iterations are used for strongly-convex optimization with a \emph{constant} step size (rather than a decreasing sequence), then \auth{Nedic and Bertsekas}\citet[Proposition 3.4]{nedic2001convergence} showed that the convergence rate of the method can be split into two parts. The first part depends on $k$ and converges linearly to~$0$. The second part is independent of $k$ and does not converge to~$0$. Thus, with a constant step size, the SG iterations have a linear convergence rate up to some tolerance, and in general after this point the iterations do not make further progress. Indeed, up until the recent work of \auth{Bach and Moulines}\citet{bach2013nonStrongly}, convergence of the basic SG method with a constant step size had only been shown for the strongly-convex quadratic case (with averaging of the iterates)~\citep{polyak1992acceleration}, or under extremely strong assumptions about the relationship between the functions $f_i$~\citep{solodov1998incremental}. This contrasts with the method we present in this work which converges to the optimal solution using a constant step size \emph{and} does so with a linear rate (without additional assumptions). 

\textbf{Accelerated methods}: Accelerated SG methods, which despite their name are not related to the aforementioned AFG method, take advantage of the fast convergence rate of SG methods with a constant step size. In particular, accelerated SG methods use a constant step size by default, and only decrease the step size on iterations where the inner-product between successive gradient estimates is negative~\citep{kesten1958accelerated,delyon1993accelerated}. This leads to convergence of the method and allows it to potentially achieve periods of faster convergence where the step size stays constant. However, the overall convergence rate of the method is not improved.

\textbf{Hybrid Methods}: Some authors have proposed variants of the SG
method for problems of the form~\eqref{eq:1} that seek to gradually
transform the iterates into the FG method in order to achieve a faster
convergence rate. \auth{Bertsekas}\citet{bertsekas1997new} proposes to go through the data
cyclically with a specialized weighting that allows the method to achieve a
linear convergence rate for strongly-convex quadratic functions. However,
the weighting is numerically unstable and the linear convergence rate
presented treats full passes through the data as iterations. A related strategy is to
 group the functions $f_i$ into `batches' of increasing size and perform SG iterations on the batches. \auth{Friedlander and Schmidt}\citet{friedlander2011hybrid} give conditions under which this strategy achieves the $O(1/k)$ and $O(\rho^k)$
 convergence rates of FG methods. However, in both
cases the iterations that achieve the faster rates have a cost that is not
independent of $n$, as opposed to SAG iterations.

\textbf{Incremental Aggregated Gradient}:
\auth{Blatt et al.}\citet{blatt2008convergent} present the most closely-related
algorithm to the SAG algorithm, the IAG method. This method is identical to the SAG
iteration~\eqref{eq:SAG}, but uses a cyclic choice of $i_k$ rather than
sampling the $i_k$ values. This distinction has several important
consequences. In particular, Blatt et al.~are only able to show that the
convergence rate is linear for strongly-convex quadratic functions (without
deriving an explicit rate), and their analysis treats full passes through
the data as iterations.
Using a non-trivial extension of their analysis and a novel proof technique
involving bounding the gradient and iterates simultaneously by a Lyapunov potential function, in this work \emph{we give an  $O(1/k)$ rate for general convex functions and an explicit
linear convergence rate for general strongly-convex
functions using the SAG iterations that only examine a single function}.
Further, as our analysis and experiments show, the SAG iterations allow a much
larger step size than is required for convergence of the IAG method.
This leads to more robustness to the selection of the step size and also, if
suitably chosen, leads to a faster convergence rate and substantially improved practical
performance. 
This shows that the simple change (random selection vs.~cycling) can dramatically improve optimization performance.

\textbf{Special Problem Classes}: For certain highly-restricted classes of problems, it is possible to show faster convergence rates for methods that only operate on a single function $f_i$.
 For example, \auth{Strohmer and Vershynin}~\citet{strohmer2009randomized} show that the randomized Kaczmarz method with a particular sampling scheme achieves a linear convergence rate for the problem of solving consistent linear systems.  It is also known that the SG method with a constant step-size has the $O(1/k)$ and $O(\rho^k)$ convergence rates of FG methods if $\norm{f_i'(x)}$ is bounded by a linear function of $\norm{g'(x)}$ for all $i$ and $x$~\citet{schmidt2012fast}. This is the strong condition required by \auth{Solodov}\citet{solodov1998incremental} to show convergence of the SG method with a constant step size. Unlike these previous works, our analysis in the next section applies to general $f_i$ that satisfy standard assumptions, and only requires gradient evaluations of the functions $f_i$ rather than dual block-coordinate steps.  

{\bf Subsequent work}: Since the first version of this work was released, there has been an explosion of research into stochastic gradient methods with faster convergence rates. It has been shown that similar rates can be achieved for certain constrained and non-smooth problems, that similar rates can be achieved without the memory requirements, that Newton-like variants of the method may be possible, and that similar rates can be achieved with other algorithms.
\iftoggle{springer}
{
 In Section~\ref{discussion} of the extended version of this paper, we survey these recent developments.
}
{
 In Section~\ref{discussion}, we survey these recent developments.
}

\section{Convergence Analysis}
\label{convergence}

In our analysis we assume that each
function $f_i$ in~\eqref{eq:1} is convex and
differentiable, and that each gradient $f'_i$ is Lipschitz-continuous with
constant $L$, meaning that for all $x$ and $y$ in $\Real^p$ and each $i$ we have
\begin{equation}
\label{eq:L}
\norm{f'_i(x)-f'_i(y)}\leqslant L\norm{x-y}.
\end{equation}
This is a fairly weak assumption on the $f_i$ functions, and in cases where
the $f_i$ are twice-differentiable it is equivalent to saying that the
eigenvalues of the Hessians of each $f_i$ are bounded above by~$L$. We will
also assume the existence of at least one minimizer $x^\ast$ that achieves
the optimal function value. We denote the average iterate by $\bar{x}^k = \frac{1}{k}\sum_{i=0}^{k-1}x^i$, and the variance of the gradient norms at the 
optimum~$x^\ast$ by $\sigma^2 = \frac{1}{n} \sum_i \|f_i'(x^\ast)\|^2$. Our convergence results consider two different initializations for the $y_i^0$ variables: setting $y_i^0 = 0$ for all $i$, or setting them to the centered gradient at the initial point $x^0$ given by $y_i^0 = f_i'(x^0) - g'(x^0)$. We note that all our convergence results are expressed
in terms of expectations with respect to the internal randomization of the algorithm (the selection of the random variables $i_k$), and not with respect to the data which is assumed to be deterministic and fixed.

In addition to this basic convex case discussed above, we will also consider the case where the average function $g =
\frac{1}{n}\sum_{i=1}^n f_i$ is strongly-convex with constant $\mu > 0$, meaning that
the function $x \mapsto g(x) - \frac{\mu}{2}\norm{x}^2$ is convex. For twice-differentiable $g$, this is equivalent to requiring that the eigenvalues of the Hessian of $g$ are bounded below by $\mu$.
This is a stronger assumption that is often not satisfied in practical applications. Nevertheless, in many applications we are free to
choose a regularizer of the parameters, and thus we can add an $\ell_2$-regularization term as in~\eqref{eq:L2}
to transform any convex problem into a strongly-convex problem (in
this case we have $\mu \geq \lambda$).
Note that strong-convexity implies the existence of a unique $x^\ast$ that achieves
the optimal function value. 

Under these standard assumptions, we now state our convergence result.
\begin{theorem}
\label{thm}
With a constant step size of $\alpha_k = \frac{1}{16L}$, the SAG iterations satisfy for $k \geq 1$:
\begin{align*}
\mathbb{E}[g(\bar{x}^k)] - g(x^*) &\leqslant \frac{32n}{k}C_0,
\end{align*}
where if we initialize with $y_i^0 = 0$ we have
\[
C_0 = g(x^0) - g(x^*) + \frac{4L}{n}\norm{x^0 - x^\ast}^2 + \frac{\sigma^2}{16L},
\]
and if we initialize with $y_i^0 = f_i'(x^0) - g'(x^0)$ we have
\[
C_0 = \frac{3}{2}\left[g(x^0) - g(x^*)\right] + \frac{4L}{n}\norm{x^0 - x^\ast}^2.
\]
Further, if $g$ is $\mu$-strongly convex we have
\begin{align*}
\mathbb{E}[g(x^k)] - g(x^*) &\leqslant \Big(1 - \min\Big\{\frac{\mu}{16L},\frac{1}{8n}\Big\}\Big)^kC_0.
\end{align*}
\end{theorem}
The proof is given in Appendix~B of the extended version of this paper, and involves finding a  Lyapunov function for a non-linear stochastic dynamical system defined on the $y_i^k$ and $x^k$ variables that converges to zero at the above rates, and showing that this function dominates the expected sub-optimality $[\mathbb{E}[g(x^k)] - g(x^*)]$. This is the same approach used to show Proposition~1 and~2 in the conference version of the paper~\citep{roux2012stochastic}, but in this work we use a more general Lyapunov function that gives a much faster rate for ill-conditioned problems and also allows us to analyze problems that are not strongly-convex. To simplify the analysis of this more complicated Lyapunov function, our new proof verifies positivity of certain polynomials that arise in the bound using a computer-aided approach.

Note that while the first part of Theorem~\eqref{thm} is stated for the average $\bar{x}^k$, with a trivial change to the proof technique it can be shown to also hold for any iterate $x^k$ where $g(x^k)$ is lower than the average function value up to iteration $k$, $\frac{1}{k}\sum_{i=0}^{k-1}g(x^i)$. Thus, in addition to $\bar{x}^k$ the result also holds for the best iterate. We also note that our bounds are valid for any $L$ greater than or equal to the minimum $L$ satisfying~\eqref{eq:L}, implying an $O(1/k)$ and linear convergence rate for any $\alpha_k \leqslant 1/16L$ (but the bound becomes worse as $L$ grows). Although initializing each $y_i^0$ with the centered gradient may have an additional cost and slightly worsens the dependency on the initial sub-optimality $(g(x^0) - g(x^\ast))$, it removes the dependency on the variance $\sigma^2$ of the gradients at the optimum.
While we have stated Theorem~\ref{thm} in terms of the function values, in the strongly-convex case we also obtain a convergence rate on the iterates because we have
\[
\frac{\mu}{2}\norm{x^k-x^\ast}^2 \leqslant g(x^k) - g(x^\ast).
\]

Theorem~\ref{thm} shows that the SAG iterations are
advantageous over SG methods in later iterations because they obtain a faster convergence rate. However, the SAG iterations have a worse constant factor because of the dependence on $n$. We can improve the dependence on $n$ using an appropriate choice of $x^0$. In particular, following~\citet{roux2012stochastic} we can set $x^0$ to the result of $n$ iterations of an appropriate SG method. In this setting, the expectation of $g(x^0) - g(x^*)$ is $O(1/\sqrt{n})$ in the convex setting, while both $g(x^0) - g(x^*)$ and $\norm{x^0 - x^\ast}^2$ would be in $O(1/n)$ in the strongly-convex setting. If we use this initialization of $x^0$ and set $y_i^0 = f_i'(x^0) - g'(x^0)$, then in terms of $n$ and $k$
 the SAG convergence rates take the form $O(\sqrt{n}/k)$ and $O(\rho^k/n)$ in the convex and strongly-convex settings, instead of the $O(n/k)$ and $O(\rho^k)$ rates implied by Theorem~\ref{thm}. 
However, in our experiments we do not use an SG initialization but rather use a minor variant of SAG in the early iterations (discussed in the next section), which appears more difficult to analyze but which gives better empirical performance. 

An interesting consequence of using a step-size of $\alpha_k = 1/16L$ is that it makes the method \emph{adaptive} to the strong-convexity constant $\mu$. That is, for problems with a higher degree of \emph{local} strong-convexity around the solution $x^\ast$, the algorithm will automatically take 
advantage of this and yield a faster local rate. This can even lead to a local linear convergence rate if the problem is strongly-convex near the optimum but not globally strongly-convex. This adaptivity to the problem difficulty is in contrast to SG methods whose sequence of step sizes typically depend on global constants and thus do not adapt to local strong-convexity.

We have observed in practice that the IAG method with a step size of $\alpha_k = \frac{1}{16L}$ may diverge. While the step-size needed for convergence of the IAG iterations is not precisely characterized, we have observed that it requires a step-size of approximately $1/nL$ in order to converge. Thus, the SAG iterations can tolerate a step size that is roughly $n$ times larger, which leads to increased robustness to the selection of the step size. Further, as our analysis and experiments indicate, the ability to use a
large step size leads to improved performance of the SAG iterations. Note that using randomized selection with a larger step-size leading to vastly improved performance
is not an unprecedented phenomenon; the analysis of \auth{Nedic and Bertsekas}~\citet{nedic2001convergence} shows
that the iterations of the basic stochastic gradient method with a constant step-size can achieve the same error bound
as full cycles through the data of the cyclic variant of the method by using
steps that are $n$ times larger (see the discussion after Proposition 3.4). Related results also appear in~\citet{collins2008exponentiated,blockFrankWolfe} showing the advantage of stochastic optimization strategies over deterministic optimization strategies in the context of certain dual optimization problems.

The convergence rate of the SAG iterations in the strongly-convex case takes a somewhat surprising form. 
For ill-conditioned problems where $n \leqslant \frac{2L}{\mu}$, $n$ does not appear in the convergence rate and  \emph{the SAG algorithm has nearly the same convergence rate as the FG method} with a step size of $1/16L$, even though it uses iterations which are $n$ times cheaper. This indicates that the basic gradient method (under a slightly sub-optimal step-size) is not slowed down by using out-of-date gradient measurements for ill-conditioned problems. Although $n$ appears in the convergence rate in the well-conditioned setting where $n > \frac{2L}{\mu}$, if we perform $n$ iterations of SAG (i.e., one effective pass through the data), the error is multiplied by $(1 - 1/8n)^n \leqslant \exp(-1/8)$, which is independent of $n$. Thus, in this setting each pass through the data reduces the excess objective by a constant multiplicative factor that is independent of the problem. 

It is interesting to compare the convergence rate of SAG in the strongly-convex case with the known convergence rates for first-order methods~\citesee{\S2}{nesterov2004introductory}.  In Table~\ref{tab:rates}, we use two examples to compare the convergence rate of SAG to the convergence rates of the standard FG method, the faster AFG method, and the lower-bound for any first-order strategy (under certain dimensionality assumptions) for optimizing a function $g$ satisfying our assumptions.  In this table, we compare the rate obtained for these FG methods to the rate obtained by running $n$ iterations of SAG, since this requires the same number of evaluations of $f_i'$. 
Case~$1$ in this table focuses on a well-conditioned case where the rate of SAG is $(1-1/8n)$, while Case~$2$ focuses on an ill-conditioned example where the rate is $(1-\mu/16L)$. Note that in the latter case the $O(1/k)$ rate for the method may be faster.

In Table~\ref{tab:rates} we see that performing $n$ iterations of SAG can actually lead to a rate that is faster than the lower bound for FG methods. Thus, for certain problems \emph{SAG can be substantially faster than any FG method that does not use the structure of the problem}. However, we note that the comparison is somewhat problematic because $L$ in the SAG rates is the Lipschitz constant of the $f_i'$ functions, while in the FG method we only require that $L$ is an upper bound on the Lipschitz continuity of $g'$ so it may be much smaller.  To give a concrete example that takes this into account and also considers the rates of dual methods and coordinate-wise methods, in Appendix A of the extended version of this paper we attempt to more carefully compare the rates obtained for SAG with the rates of primal and dual FG and coordinate-wise methods for the special case of $\ell_2$-regularized least-squares regression.

\begin{table}
\centering
\begin{tabular}{|l|c|c|c|c|}
\hline
Algorithm & Step Size & Theoretical Rate & Rate in Case 1 & Rate Case 2\\
\hline
FG & $\frac{1}{L}$ & $\left(1-\frac{\mu}{L}\right)^2$ & $0.9998$ & $1.000$\\
FG & $\frac{2}{\mu+L}$ & $\left(1-\frac{2\mu}{L+\mu}\right)^2$ & $0.9996$ & $1.000$\\
AFG & $\frac{1}{L}$ & $\left(1-\sqrt{\frac{\mu}{L}}\right)$ & $0.9900$ & $0.9990$ \\
Lower-Bound & --- & $\left(1-\frac{2\sqrt{\mu}}{\sqrt{L}+\sqrt{\mu}}\right)^2$ & $0.9608$ & $0.9960$\\
\hline
SAG ($n$ iters) & $\frac{1}{16L}$ & $\left(1-\min\{\frac{\mu}{16L},\frac{1}{8n}\}\right)^n$ & $0.8825$ & $0.9938$\\
\hline
\end{tabular}
\caption{Comparison of convergence rates of first-order methods to the convergence rates of $n$ iterations of SAG.  
In the examples we take $n=100000$, $L=100$, $\mu=0.01$ (Case 1), and $\mu=0.0001$ (Case 2).}
\label{tab:rates}
\end{table}

\section{Implementation Details}
\label{sec:implementation}


In Algorithm~\ref{alg:SAG} we give pseudo-code for an implementation of the basic method, where we use a variable $d$ to track the quantity $(\sum_{i=1}^ny_i)$. 
This section focuses on further implementation details that are useful in practice.
 In particular, we discuss modifications that lead to better practical performance than the basic Algorithm~\ref{alg:SAG}, including ways to reduce the storage cost, how to handle regularization, how to set the step size, using mini-batches, and using non-uniform sampling. Note that an implementation of the algorithm that incorporates many of these aspects is available from the first author's webpage.

\begin{algorithm}
\caption{Basic SAG method for minimizing $\frac{1}{n}\sum_{i=1}^nf_i(x)$ with step size $\alpha$.}
\label{alg:SAG}
\begin{algorithmic}
\STATE $d = 0$, $y_i = 0$ for $i = 1,2,\dots,n$
\FOR{$k=0,1,\dots$}
\STATE Sample $i$ from $\{1,2,\dots,n\}$
\STATE $d = d - y_i + f_i'(x)$
\STATE $y_i = f_i'(x)$
\STATE $x = x - \frac{\alpha}{n} d$
\ENDFOR
\end{algorithmic}
\end{algorithm}

\subsection{Structured gradients and just-in-time parameter updates}
\label{sec:sparse}

For many problems the storage cost of $O(np)$ for the $y_i^k$ vectors is prohibitive, but we can often use the structure of the gradients $f_i'$ to reduce this cost. For example, a commonly-used specialization of~\eqref{eq:1} is \emph{linearly-parameterized} models which take form
\begin{equation}
\label{eq:atx}
\minimize{x\in\Real^p}\quad g(x) \defd \frac{1}{n}\sum_{i=1}^n f_i(a_i^\top x).
\end{equation}
Since each $a_i$ is constant, for these problems we only need to store the scalar $f_{i_k}'(u_i^k)$ for $u_i^k = a_{i_k}^\top x^k$ rather than the full gradient $a_i f_i'(u_i^k)$. This reduces the storage cost from $O(np)$ down to $O(n)$.



For problems 
where the vectors $a_i$ are sparse, an individual gradient $f_i'$ will inherit the sparsity pattern of the corresponding $a_i$. However,
the update of $x$ in Algorithm~\ref{alg:SAG} appears unappealing since in general $d$ will be dense, resulting in an iteration cost of $O(p)$. Nevertheless, we can take advantage of the simple form of the SAG updates to implement a `just-in-time' variant of the SAG algorithm where the iteration cost is proportional to the number of non-zeroes in $a_{i_k}$. In particular, we do this by not explicitly storing the full vector $x^k$ after each iteration.  Instead, on each iteration we only compute the elements $x_j^k$ corresponding to non-zero elements of $a_{i_k}$, by applying the \emph{sequence of updates} to each variable $x_j^k$ since the last iteration where it was non-zero in $a_{i_k}$. This sequence of updates can be applied efficiently since it simply involves changing the step size. For example, if variable $j$ has been zero in $a_{i_k}$ for $5$ iterations, then we can compute the needed value $x_j^k$ using
\[
x_j^k = x_j^{k-5} - \frac{5\alpha}{n}\sum_{i=1}^n(y_i^{k})_j.
\]
This update 
allows SAG to be efficiently applied to sparse data sets where $n$ and $p$ are both in the millions or higher but the number of non-zeros is much less than $np$.

\subsection{Re-weighting on early iterations}
\label{sec:reweight}

In the update of $x$ in Algorithm~\ref{alg:SAG}, we normalize the direction $d$ by the total number of data points~$n$.  When initializing with $y_i^0 = 0$ we believe this leads to steps that are too small on early iterations of the algorithm where we have only seen a fraction of the data points, because many $y_i$ variables contributing to $d$ are set to the uninformative zero-vector.  
Following \auth{Blatt et al.}\citet{blatt2008convergent}, the more logical normalization is to divide $d$ by $m$, the number of data points that we have seen at least once (which converges to $n$ once we have seen the entire data set), leading to the update $x = x - \frac{\alpha}{m}d$.  Although this modified SAG method appears more difficult to analyze, in our experiments we found that running the basic SAG algorithm from the very beginning with this modification outperformed the basic SAG algorithm as well as the SG/SAG hybrid algorithm mentioned in the Section~\ref{convergence}.  In addition to using the gradient information collected during the first $k$ iterations, this modified SAG algorithm is also advantageous over hybrid SG/SAG algorithms because it only requires estimating a single constant step size.

\subsection{Exact and efficient regularization}
\label{sec:regul}

In the case of regularized objectives like~\eqref{eq:L2}, the cost of computing the gradient of the regularizer is independent of $n$. 
Thus, we can use the exact gradient of the regularizer in the update of $x$, and only use $d$ to approximate the sum of the $l_i'$ functions. 
By incorporating the gradient of the regularizer explicitly, the update for $y_i$ in Algorithm~\ref{alg:SAG} becomes $y_i = l_i'(x)$, and in the case of $\ell_2$-regularization the update for $x$ becomes
\[
x = x - \alpha\left(\frac{1}{m}d + \lambda x\right) = \left(1-\alpha\lambda\right)x - \frac{\alpha}{m}d.
\]
If the loss function gradients $l_i'$ are sparse as in Section~\ref{sec:sparse}, then these modifications lead to a reduced storage requirement even though the gradient of the regularizer is dense.
 Further,
although the update of $x$ again appears to require dense vector operations, we can implement the algorithm efficiently if the $a_i$ are sparse. In particular, to allow efficient multiplication of $x$ by the scalar~$(1-\alpha\lambda)$, it is useful to represent $x$ in the form $x = \kappa z$, where $\kappa$ is a scalar and $z$ is a vector (as done by \auth{Shalev-Shwartz et al.}\citet{shalev2007pegasos}). Under this representation, we can multiply $x$ by a scalar in $O(1)$ by simply updating $\kappa$ (though to prevent $\kappa$ becoming too large or too small we may need to occasionally re-normalize by setting $z = \kappa z$ and $\kappa=1$). To efficiently implement the vector subtraction operation, we can use a variant of the just-in-time updates from Section~\ref{sec:sparse}.  In Algorithm~\ref{alg:SAG-L2}, we give pseudo-code for a variant of SAG that includes 
all of these modifications, and thus uses no full-vector operations.  This code uses a vector $y$ to keep track of the scalars $l_i'(u_i^k)$, a vector $C$ to determine whether a data point has previously been visited, a vector $V$ to track the last time a variable was updated, and a vector $S$ to keep track of the cumulative sums needed to implement the just-in-time updates.

\begin{algorithm}
\caption{SAG variant for minimizing $\frac{\lambda}{2}\norm{x}^2 + \frac{1}{n}\sum_{i=1}^nl_i(a_i^\top x)$, with step size $\alpha$ and $a_i$ sparse.}
\label{alg:SAG-L2}
\begin{algorithmic}
\STATE \COMMENT{Initialization, note that $x = \kappa z$.}
\STATE $d = 0, y_i = 0$ for $i = 1,2,\dots,n$
\STATE $z=x$, $\kappa=1$
\STATE $m=0$,  $C_i = 0$ for $i=1,2,\dots,n$
\STATE $S_{-1} = 0$, $V_j = 0$ for $j=1,2,\dots,p$
\FOR{$k = 0,1,\dots$}
\STATE Sample $i$ from $\{1,2,\dots,n\}$
\IF{$C_i = 0$}
\STATE \COMMENT{This is the first time we have sampled this data point.}
\STATE $m=m+1$
\STATE $C_i = 1$
\ENDIF
\FOR{$j$ non-zero in $a_i$}
\STATE \COMMENT{Just-in-time calculation of needed values of $z$.}
\STATE $z_j = z_j - (S_{k-1} - S_{V_j-1})d_j$
\STATE $V_j = k$
\ENDFOR
\STATE \COMMENT{Update the memory $y$ and the direction $d$.}
\STATE Let $J$ be the support of $a_i$
\STATE $d_J = d_J - a_{iJ}(y_i - l_i'(\kappa a_{iJ}^Tz_J))$
\STATE $y_i = l_i'(\kappa a_{iJ}^Tz_J)$
%
\STATE \COMMENT{Update $\kappa$ and the sum needed for $z$ updates.}
\STATE $\kappa = \kappa(1-\alpha\lambda)$
\STATE $S_k = S_{k-1} + \alpha/(\kappa m)$
\ENDFOR
\STATE \COMMENT{Final $x$ is $\kappa$ times the just-in-time update of all $z$.}
\FOR{$j = 1,2,\dots,p$}
\STATE $x_j = \kappa(z_j - (S_{k-1} - S_{V_j-1})d_j)$
\ENDFOR
\end{algorithmic}
\end{algorithm}

\subsection{Warm starting}

In many scenarios we may need to solve a set of closely-related optimization problems.  For example, we may want to apply Algorithm~\ref{alg:SAG-L2} to a regularized objective of the form~\eqref{eq:L2} for several values of the regularization parameter $\lambda$.
Rather than solving these problems independently, we might expect to obtain better performance by warm-starting the algorithm.
Although initializing $x$ with the solution of a related problem can improve performance, we can expect an even larger performance improvement if we also use the gradient information collected from a run of SAG for a close value of $\lambda$.  For example, in Algorithm~\ref{alg:SAG-L2} we could initialize $x$, $y_i$, $d$, $m$, and $C_i$ based on a previous run of the SAG algorithm. In this scenario, Theorem~\ref{thm} suggests that it may be beneficial in this setting to center the $y_i$ variables around $d$.

\subsection{Larger step sizes}
\label{sec:stepSizes}

In our experiments we have observed that utilizing a step size of $1/L$, as in standard FG methods, always converged and often performed better than the step size of $1/16L$ suggested by our analysis. Thus, in our experiments we used $\alpha_k = 1/L$ even though we do not have a formal analysis of the method under this step size. 
We also found that a step size of $2/(L + n\mu)$, which in the strongly-convex  case corresponds to the best fixed step size for the FG method in the special case of $n=1$~\citesee{Theorem 2.1.15}{nesterov2004introductory}, sometimes yields even better performance (though in other cases it performs poorly).

\subsection{Line-search when $L$ is not known}
\label{sec:linesearch}

In general the Lipschitz constant $L$ will not be known, but we may obtain a reasonable approximation of a valid $L$ by evaluating $f_i$ values while running the algorithm. In our experiments, we used a basic line-search
 where we start with an initial estimate $L^0$, and double this estimate whenever we do not satisfy the  inequality
\[
f_{i_k}(x^k - \frac{1}{L^k}f_{i_k}'(x^k)) \leqslant f_{i_k}(x^k) - \frac{1}{2L^k}\norm{f_{i_k}'(x^k)}^2,
\]
which must be true if $L^k$ is valid.
An important property of this test is that it depends on $f_{i_k}$ but not on $g$, and thus the cost of performing this test is independent of $n$.
 To avoid instability caused by comparing very small numbers, we only do this test when $\norm{f_{i_k}'(x^k)}^2 > 10^{-8}$. 
Since $L$ is a global quantity but the algorithm will eventually remain within a neighbourhood of the optimal solution, it is possible that a smaller estimate of $L$ (and thus a larger step size) can be used as we approach~$x^\ast$. To potentially take advantage of this, we initialize with the slightly smaller $L^k = (L^{k-1}2^{-1/n})$ at each iteration, so that the estimate of $L$ is halved if we do $n$ iterations (an effective pass through the data) and never violate the inequality. Note that in the case of $\ell_2$-regularized objectives, we can perform the line-search to find an estimate of the Lipschitz constant of $l_i'$ rather than $f_i'$, and then simply add $\lambda$ to this value to take into account the effect of the regularizer. 

Note that the cost of this line-search is \emph{independent} of $n$, making it suitable for large problems. Further, for linearly-parameterized models of the form $f_i(a_i^Tx)$, it is also possible to implement the line-search so that its cost is also independent of the number of variables $p$. To see why, if we use $\delta^k = a_{i_k}^Tx^k$ and the structure in the gradient then the left side is given by
\[
f_{i_k}\left(a_{i_k}^T\left(x^k - \frac{1}{L^k}f_{i_k}'(x^k)\right)\right) = f_{i_k}\left(\delta^k - \frac{f_{i_k}'(\delta^k)}{L^k}\norm{a_{i_k}}^2\right).
\]
Thus, if we pre-compute the squared norms $\norm{a_i}^2$ and note that $\delta^k$ and $f_{i_k}'(\delta^k)$ are already needed by the SAG update, then each iteration only involves operations on scalar values and the single-variable function $f_{i_k}$.

\subsection{Mini-batches for vectorized computation and reduced storage}
\label{sec:miniBatch}

Because of the use of vectorization and parallelism in modern architectures, practical SG implementations often group functions into `mini-batches' and perform SG iterations on the mini-batches. We can also use mini-batches within the SAG iterations to take advantage of the same vectorization and parallelism. Additionally, for problems with dense gradients mini-batches can dramatically decrease the storage requirements of the algorithm, since we only need to store a vector $y_i$ for each mini-batch rather than for each example. Thus, for example, using a mini-batch of size $100$ leads to a $100$-fold reduction in the storage cost.

A subtle issue that arises when using mini-batches is that the value of $L$ in the Lipschitz condition~\eqref{eq:L} is based on the mini-batches instead of the original functions $f_i$. For example, consider the case where we have a batch $\mathcal{B}$ and the minimum value of $L$ in~\eqref{eq:L} for each $i$ is given by $L_i$. In this case, a valid value of $L$ for the function $x \mapsto \frac{1}{|\mathcal{B}|}\sum_{i\in\mathcal{B}}f_i(x)$ would be $\max_{i\in\mathcal{B}}\{L_i\}$. We refer to this as $L_\textrm{max}$. But we could also consider using $L_\textrm{mean} = \frac{1}{|\mathcal{B}|}\sum_{i\in\mathcal{B}}L_i$. The value $L_\textrm{mean}$ is still valid and will be smaller than $L_\textrm{max}$ unless all $L_i$ are equal. We could even consider the minimum possible value of $L$, which we refer to as $L_\textrm{Hessian}$ because (if each $f_i$ is twice-differentiable) it is equal to the maximum eigenvalue of $\frac{1}{|\mathcal{B}|}\sum_{i\in\mathcal{B}}f_i''(x)$ across all $x$. Note that $L_\textrm{Hessian} \leq L_\textrm{mean} \leq L_\textrm{max}$, although $L_\textrm{Hessian}$ will typically be more difficult to compute than $L_\textrm{mean}$ or $L_\textrm{max}$ (although a line-search as discussed in the previous section can reduce this cost). Due to the potential of using a smaller $L$, \emph{we may obtain a faster convergence rate by using larger mini-batches}. However, in terms of passes through the data this faster convergence may be offset by the higher iteration cost associated with using mini-batches.

\subsection{Non-uniform example selection}
\label{sec:lipschitz}

In standard SG methods, it is crucial to sample the functions $f_i$ uniformly, at least asymptotically, in order to yield an unbiased gradient estimate and subsequently achieve convergence to the optimal value (alternately, the bias induced by non-uniform sampling would need to be asymptotically corrected). In SAG iterations, however, the weight of each gradient is constant and equal to $1/n$, regardless of the frequency at which the corresponding function is sampled.  We might thus consider sampling the functions $f_i$ non-uniformly, without needing to correct for this bias.
Though we do not yet have any theoretical proof as to why a non-uniform sampling might be beneficial, intuitively we would expect that we do not
need to sample functions $f_i$ whose gradient changes slowly as often as functions $f_i$ whose gradient changes more quickly. Indeed, we provide here an argument to justify a non-uniform sampling strategy based on the Lipschitz constants of the individual gradients $f_i'$ and we note that in subsequent works this intuition has proved correct for related algorithms~\citep{xiao2014proximal,schmidt2015sag4crf}.

Let $L_i$ again be the Lipschitz constant of $f_i'$, and assume that the functions are placed in increasing order of Lipschitz constants, so that $L_1 \leqslant L_2 \leqslant \ldots \leqslant L_n$.  In the ill-conditioned setting where the convergence rate depends on $\frac{\mu}{L}$, a simple  way to improve the rate by decreasing $L$ is to replace $f_n$ by two functions $f_{n1}$ and $f_{n2}$ such that
\begin{align*}
f_{n1}(x) = f_{n2}(x) &= \frac{ f_n(x) }{2}\\
g(x) &= \frac{1}{n} \left(\sum_{i=1}^{n-1} f_i(x) + f_{n1}(x) + f_{n2}(x)\right)\\
&= \frac{1}{n+1} \left(\sum_{i=1}^{n-1} \frac{n+1}{n}f_i(x) + \frac{n+1}{n}f_{n1}(x) + \frac{n+1}{n}f_{n2}(x)\right) \; .
\end{align*}
We have thus replaced the original problem by a new, equivalent problem where:
\begin{list}{\labelitemi}{\leftmargin=1.7em}
\item $n$ has been replaced by $(n+1)$,
\item $L_i$ for $i \leqslant (n-1)$ is $\frac{L_i(n+1)}{n}$,
\item $L_n$ and $L_{n+1}$ are equal to $\frac{L_n(n+1)}{2n}$.
\end{list}
Hence, if $L_{n-1} < \frac{nL_n}{n+1}$, this problem has the same $\mu$ but a smaller $L$ than the original one, improving the bound on the convergence rate. 
By duplicating $f_n$, we increase its probability of being sampled from $\frac{1}{n}$ to $\frac{2}{n+1}$,  but we also replace $y_n^k$ by a noisier version, i.e.~$y_{n1}^k + y_{n2}^k$.  Using a noisier version of the gradient appears detrimental, so we assume that the improvement comes from increasing the frequency at which $f_n$ is sampled, and that logically we might obtain a better rate by simply sampling $f_n$ more often in the original problem and not explicitly duplicating the data.

We now consider the extreme case of duplicating each function $f_i$ a number of times equal to the Lipschitz constant of their gradient, assuming that these constants are integers. The new problem becomes
\begin{align*}
g(x) &= \frac{1}{n}\sum_{i=1}^{n} f_i(x)\\
&= \frac{1}{n}\sum_{i=1}^{n} \sum_{j = 1}^{L_i} \frac{f_i(x)}{L_i}\\
&= \frac{1}{\sum_k L_k} \sum_{i=1}^{n} \sum_{j = 1}^{L_i} \left(\frac{\sum_k L_k}{n}\frac{f_i(x)}{L_i}\right) \; .
\end{align*}
The function $g$ is now written as the sum of $\sum_k L_k$ functions, each with a gradient with Lipschitz constant $\frac{\sum_k L_k}{n}$. The new problem has the same $\mu$ as before, but now has an $L$ equal to the average of the Lipschitz constants across the $f_i'$, rather than their maximum, thus improving the bound on the convergence rate.
Sampling these functions uniformly is now equivalent to sampling the original $f_i$'s according to their Lipschitz constant.
Thus, we might expect to obtain better performance by, instead of creating a larger problem by duplicating the functions in proportion to their Lipschitz constant, simply sampling the functions from the original problem in proportion to their Lipschitz constants.

Sampling in proportion to the Lipschitz constants of the gradients was explored by \auth{Nesterov}\citet{nesterov2010efficiency} in the context of coordinate descent methods, and is also somewhat related to the sampling scheme used by \auth{Storhmer and Vershynin}\citet{strohmer2009randomized} in the context of their randomized Kaczmarz algorithm. Since the first version of this work was released,~\auth{Needell et al.}\citet{SGkaczmarz} have analyzed sampling according to the Lipschitz constant in the context of SG iterations.
Such a sampling scheme makes the iteration cost depend on $n$, due to the need to generate samples from a general discrete distribution over $n$ variables. However, after an initial preprocessing cost of $O(n)$ we can sample from such distributions in $O(\log n)$ using a simple binary search~\citesee{Example~2.10}{robert2004monte}. 

Unfortunately, sampling the functions according to the Lipschitz constants and using a step size of $\alpha_k = \frac{n}{\sum_{i}L_i}$ does not seem to converge in general. However, by changing the number of times we duplicate each $f_i$, we can interpolate between the Lipschitz sampling and the uniform sampling. In particular, if each function $f_i$ is duplicated $L_i + c$ times, where $L_i$ is the Lipschitz constant of $f_i'$ and $c$ a positive number, then the new problem becomes
\begin{align*}
g(x) &= \frac{1}{n}\sum_{i=1}^{n} f_i(x)\\
&= \frac{1}{n}\sum_{i=1}^{n} \sum_{j = 1}^{L_i + c} \frac{f_i(x)}{L_i + c}\\
&= \frac{1}{\sum_k (L_k + c)} \sum_{i=1}^{n} \sum_{j = 1}^{L_i + c} \left(\frac{\sum_k (L_k + c)}{n}\frac{f_i(x)}{L_i + c}\right) \; .
\end{align*}
Unlike in the previous case, these $\sum_k (L_k + c)$ functions have gradients with different Lipschitz constants. Denoting $L = \max_i L_i$, the maximum Lipschitz constant is equal to $\frac{\sum_k (L_k + c)}{n}\frac{L}{L + c}$ and we must thus use the step size $\alpha = \frac{L+ c}{L\left(\frac{\sum_k L_k}{n} + c\right)}$.

\section{Experimental Results}
\label{experiments}
\label{sec:experiments}

In this section we perform empirical evaluations of the SAG iterations. We first compare the convergence of an implementation of the SAG iterations to a variety of competing methods available. 
\iftoggle{springer}
{
We then seek to evaluate the effect of different of non-uniform sampling. In Section 5 of the extended version of this paper, we present further experiments evaluating the effect of the step-size and the effect of using mini-batches.
}
{
We then seek to evaluate the effect of different algorithmic choices such as the step size, mini-batches, and non-uniform sampling.
}


\subsection{Comparison to FG and SG Methods}
\label{exp:compare}

The theoretical convergence rates suggest the following strategies for
deciding on whether to use an FG or an SG method:
\begin{list}{\labelitemi}{\leftmargin=1.7em}
\item If we can only afford one pass through the data, then an SG method should
be used.
\item If we can afford to do many passes through the data (say, several
hundred), then an FG method should be used.
\end{list}
We expect that the SAG iterations will be most useful between these two extremes, where we can afford to do more than one pass through the data but cannot afford to do enough passes to warrant using FG algorithms like the AFG or L-BFGS methods. To test whether this is indeed the case in practice,
we perform a variety of experiments evaluating the performance of the SAG algorithm in this scenario.

Although the SAG algorithm can be applied more generally, in our experiments we focus on the important and widely-used $\ell_2$-regularized logistic regression problem
\begin{equation}
\label{eq:logreg}
\minimize{x\in\Real^p}\quad \frac{\lambda}{2}\|x\|^2 + \frac{1}{n}\sum_{i=1}^n \log(1+\exp(-b_ia_i^\top x)),
\end{equation} 
as a canonical problem satisfying our assumptions.  In our experiments we set the regularization parameter $\lambda$ to $1/n$, which is in the range of the smallest values that would typically be used in practice, and thus which results in the most ill-conditioned problems of this form that would be encountered.  Our experiments focus on the freely-available benchmark binary classification data sets listed in Table~\ref{table:data}.
The
\emph{quantum} and \emph{protein} data set was obtained from the KDD Cup 2004 website;\footnote{\small \url{http://osmot.cs.cornell.edu/kddcup}} the \emph{covertype} (based on the datset of Blackard, Jock, and
Dean), \emph{rcv1}, \emph{news}, and \emph{rcv1Full} data sets were obtained from the LIBSVM Data website; \footnote{\small\url{http://www.csie.ntu.edu.tw/~cjlin/libsvmtools/datasets}}; 
the \emph{sido} data set was obtained from the Causality Workbench website,\footnote{ \small\url{http://www.causality.inf.ethz.ch/home.php}}
the \emph{spam} data set was prepared by \auth{Carbonetto}\citesee{\S2.6.5}{Carbonetto09} using the TREC 2005 corpus\footnote{\small\url{http://plg.uwaterloo.ca/~gvcormac/treccorpus}}; and the \emph{alpha} data set was obtained from the Pascal Large Scale Learning Challenge website\footnote{\small\url{http://largescale.ml.tu-berlin.de}}. 
We added a (regularized) bias term to all data sets, and for dense features we standardized so that they would have a mean of zero and a variance of one. To obtain results that are independent of the practical implementation of the algorithm, we measure the objective as a function of the number of effective passes through the data, measured as the number of times we evaluate $l_i'$ divided by the number of examples $n$.  If they are implemented to take advantage of the sparsity present in the data sets, the runtimes of all algorithms examined in this section differ by at most a constant times this measure.

\begin{table*}
\centering
\begin{tabular}{|l|r|r|c|}
\hline
Data set & Data Points & Variables & Reference \\
\hline
\emph{quantum} & 50 000 & 78 & \citep{caruana2004kdd}\\
\hline
\emph{protein} & 145 751 & 74 &  \citep{caruana2004kdd}\\
\hline
\emph{covertype} & 581 012 & 54  & \citep{blake1998uci}\\
\hline
\emph{rcv1} & 20 242 & 47 236  & \citep{lewis2004rcv1}\\
\hline
\emph{news} & 19 996 & 1 355 191 & \citep{keerthi2005modified}\\
\hline
\emph{spam} & 92 189 & 823 470 & \citep{cormack2005spam,Carbonetto09}\\
\hline
\emph{rcv1Full} & 697 641 & 47 236  & \citep{lewis2004rcv1}\\
\hline
\emph{sido} & 12 678 & 4 932  & \citep{SIDO}\\
\hline
\emph{alpha} & 500 000 & 500  & Synthetic\\
\hline
\end{tabular}
\caption{Binary data sets used in the experiments.}
\label{table:data}
\end{table*}

In our first experiment we compared the following variety of competitive FG and SG methods:
\begin{list}{\labelitemi}{\leftmargin=1.7em}
\item \textbf{AFG}: A variant of the accelerated full gradient method of \auth{Nesterov}\citet{nesterov1983method}, where iterations of~\eqref{eq:FG} with a step size of $1/L^k$ are interleaved with an extrapolation step. We used an adaptive line-search to estimate a local $L$ based on the variant proposed for $\ell_2$-regularized logistic regression by \auth{Liu et al.}\citet{liu2009large}.
\item \textbf{L-BFGS}: A publicly-available limited-memory quasi-Newton method that has been tuned for log-linear models such as logistic regression~\citep{minFunc}. This method is the most complicated method we considered.
\item \textbf{SG}: The stochastic gradient method described by
iteration~\eqref{eq:SG}. Since setting the step-size in this method is a tenuous issue, we chose the constant step size that gave the best performance (in hindsight) among all powers of $10$ (we found that this constant step-size strategies gave better performance than the variety of decreasing step-size strategies that we experimented with).
\item \textbf{ASG}: The average of the iterations generated by the SG method above, where again we choose the best step size among all powers of $10$.\footnote{Note that we also compared to a variety of other SG methods including the popular Pegasos SG method of \auth{Shalev-Shwartz et al.}\citet{shalev2007pegasos}, SG with momentum, SG with gradient averaging, the regularized dual averaging method of \auth{Xiao}\citet{xiao2010dual} (a stochastic variant of the primal-dual subgradient method of \auth{Nesterov}\citet{nesterov2009primal} for regularized objectives), the accelerated SG method of \auth{Delyon and Juditsky}\citet{delyon1993accelerated}, SG methods that only average the later iterations as in the `optimal' SG method for non-smooth optimization of \auth{Rakhlin et al.}\citet{rakhlin2012making}, the epoch SG method of \auth{Hazan and Kale}\citet{hazan2010beyond}, the `nearly-optimal' SG method of \auth{Ghadimi and Lan}\citet{ghadimi2010optimal}, a diagonally-scaled SG method using the inverse of the coordinate-wise Lipshitz constants as the diagonal terms, and the adaptive diagonally-scaled AdaGrad method of \auth{Duchi et al.}\citet{duchi2010adaptive}. However, we omit results obtained using these algorithms since they never performed substantially better than the minimum between the \emph{SG} and \emph{ASG} methods when their step-size was chosen in hindsight.}
\item \textbf{IAG}: The incremental aggregated gradient method
of \auth{Blatt et al.}\citet{blatt2008convergent} described by iteration~\eqref{eq:SAG} with a
cyclic choice of $i_k$. We used the re-weighting described in Section~\ref{sec:reweight}, we used the exact regularization as described in Section~\ref{sec:regul}, and we chose the step-size that gave the best performance among all powers of $10$.
\item \textbf{SAG-LS}:  The proposed stochastic average gradient method
described by iteration~\eqref{eq:SAG}. We used the re-weighting described in Section~\ref{sec:reweight}, the exact regularization as described in Section~\ref{sec:regul}, and we used a step size of $\alpha_k = 1/L^k$ where $L^k$ is an approximation of the Lipschitz constant for the negative log-likelihoods $l_i(x) = \log(1 + \exp(-b_ia_i^\top x))$. Although this Lipschitz constant is given by $0.25\max_i\{\norm{a_i}^2\}$, we used the line-search described in Section~\ref{sec:linesearch} to estimate it, to test the ability to use SAG as a black-box algorithm (in addition to
avoiding this calculation and potentially obtaining a faster convergence rate by obtaining an estimate that could be smaller than the global value). To initialize the line-search we set $L^0=1$.
\end{list}

\begin{figure}
\centering
\mbox{ \includegraphics[width=\figSize]{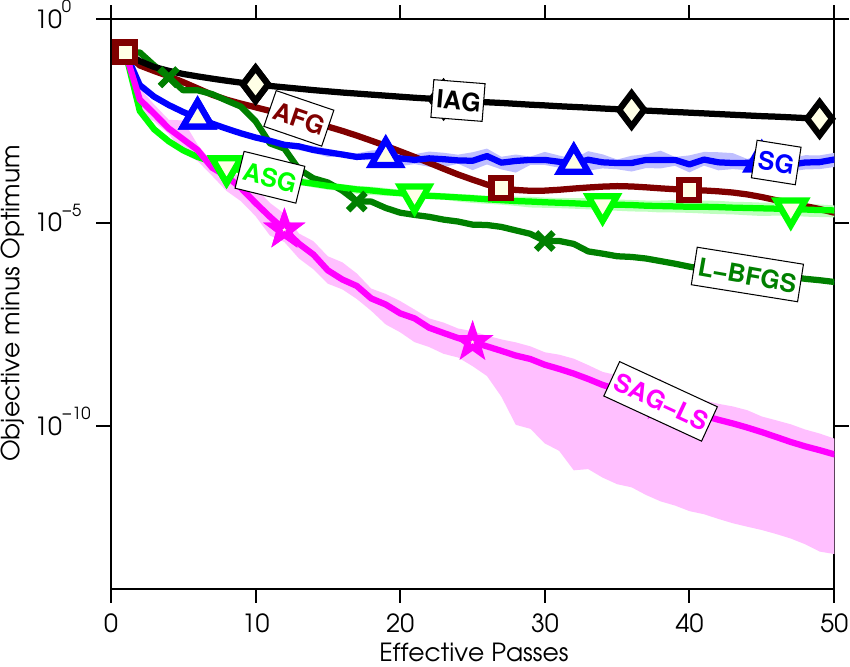} \hspace*{-.1cm}
\includegraphics[width=\figSize]{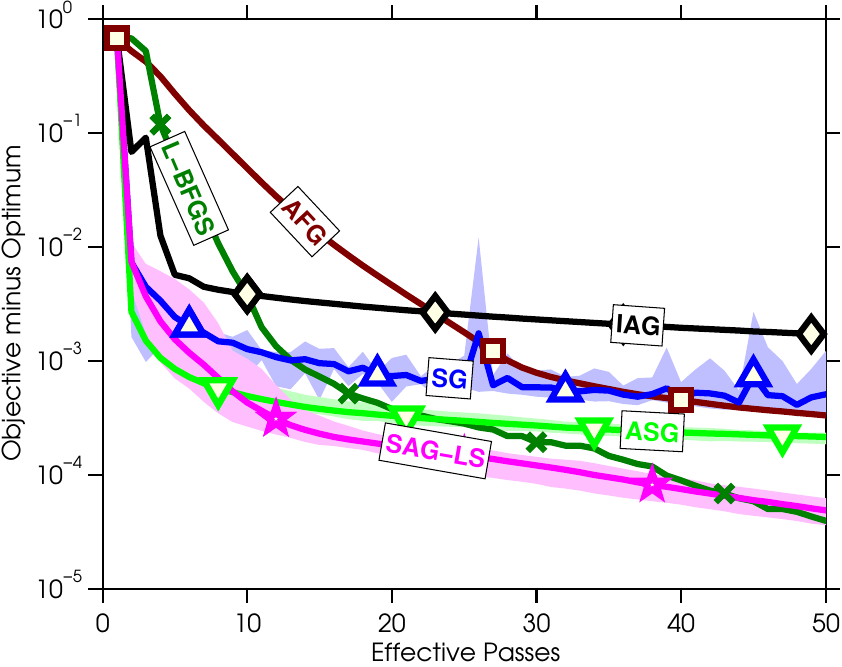} \hspace*{-.1cm}
\includegraphics[width=\figSize]{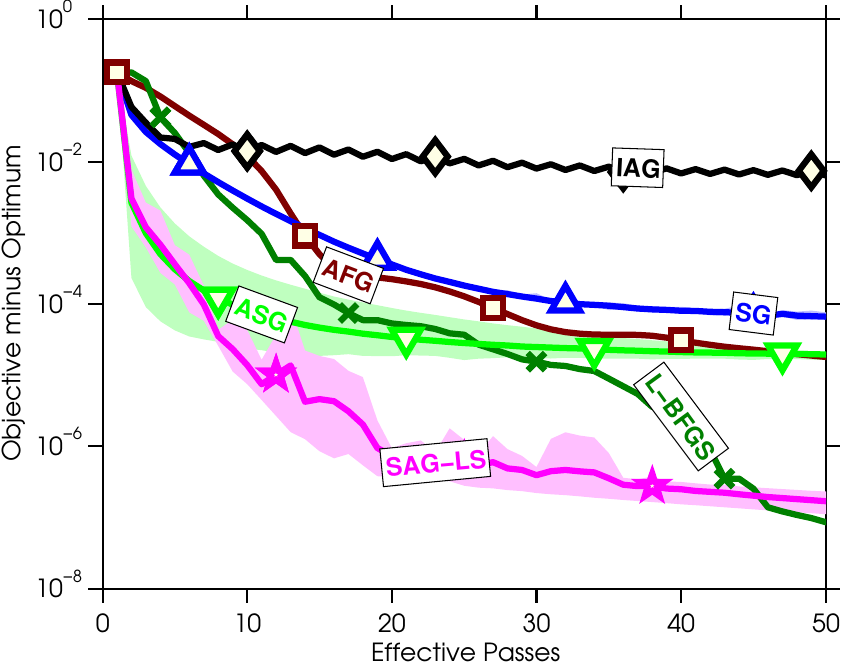} \hspace*{-.1cm}}

\mbox{\includegraphics[width=\figSize]{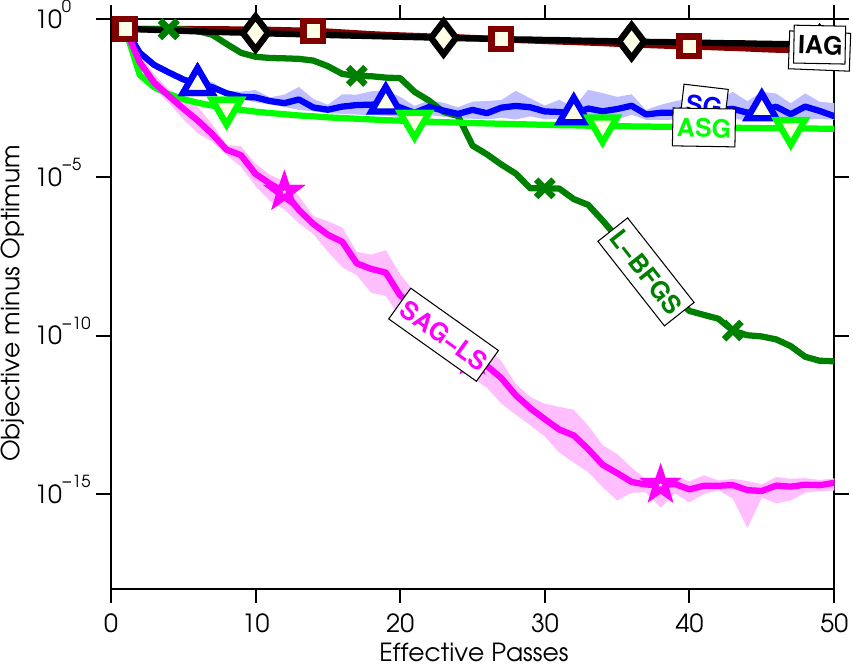} \hspace*{-.1cm}
\includegraphics[width=\figSize]{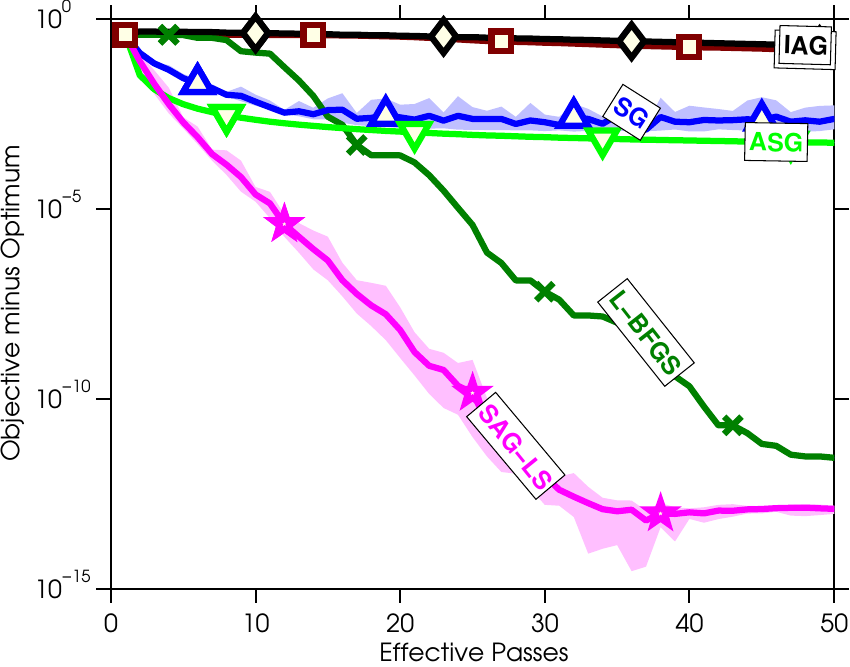} \hspace*{-.1cm}
\includegraphics[width=\figSize]{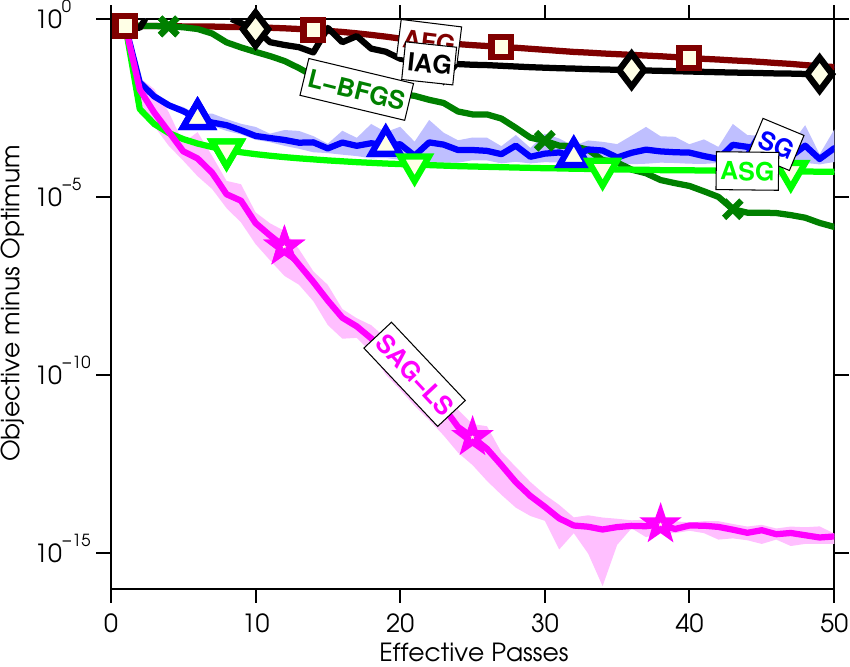} \hspace*{-.1cm}
}

\mbox{\includegraphics[width=\figSize]{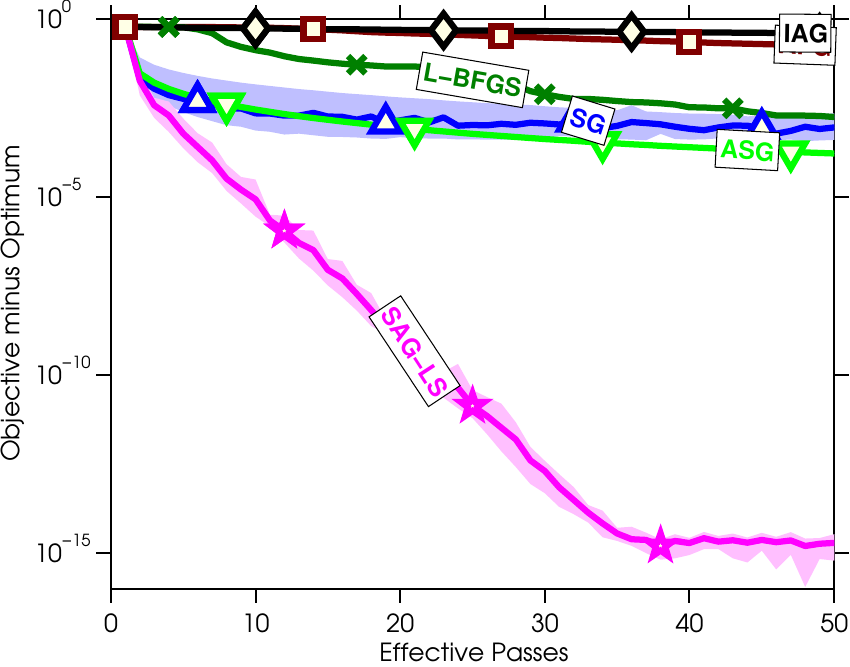} \hspace*{-.1cm}
\includegraphics[width=\figSize]{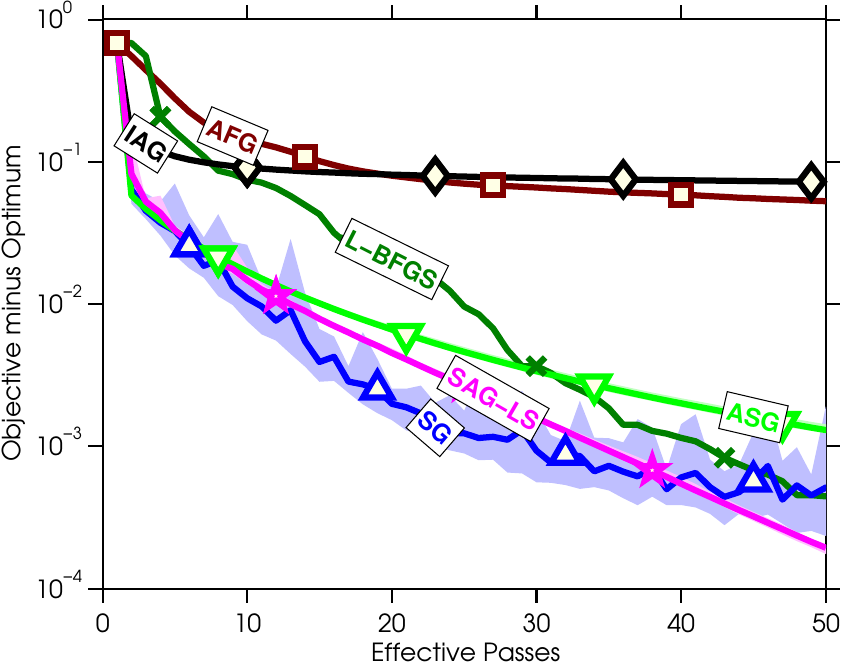} \hspace*{-.1cm}
\includegraphics[width=\figSize]{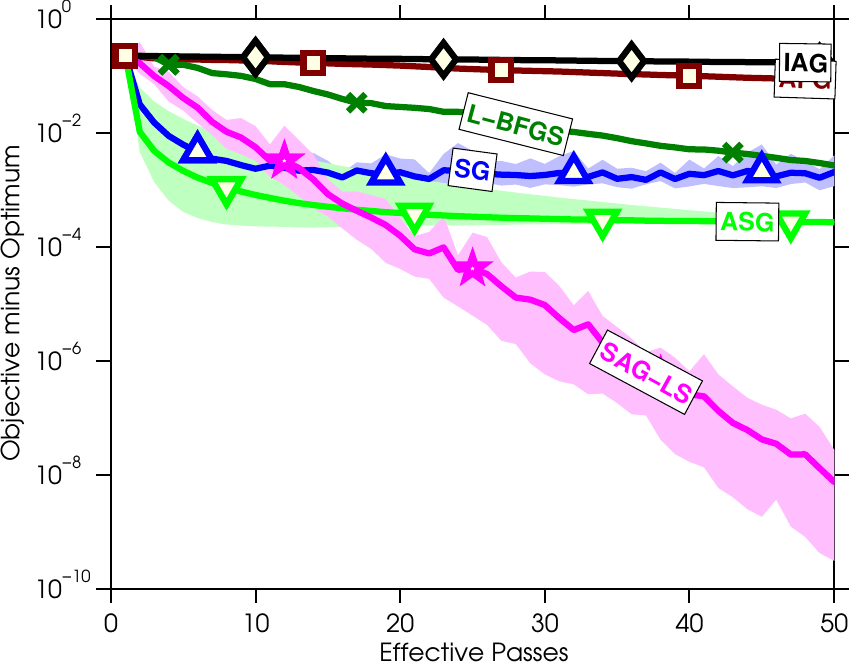} \hspace*{-.1cm}
}
\caption{Comparison of different FG and SG optimization strategies. The top row gives results on the \emph{quantum} (left), \emph{protein} (center) and \emph{covertype} (right) datasets. The middle row gives results on the \emph{rcv1} (left), \emph{news} (center) and \emph{spam} (right) datasets.  The bottom row gives results on the \emph{rcv1Full} (left), \emph{sido} (center), and \emph{alpha} (right) datasets. This figure is best viewed in colour.}
\label{fig:methods}
\end{figure}

We plot the results of the different methods for the first $50$ effective passes through the data in Figure~\ref{fig:methods}. For the stochastic methods, we plot the mean performance as well as the minimum and maximum function values across $10$ choices for the initial random seed.
We can observe several trends across the experiments:
\begin{list}{\labelitemi}{\leftmargin=1.7em}
  \item {\bf FG vs.~SG}: Although the performance of SG methods is known to be catastrophic if the step size is not chosen carefully, after giving the SG methods (\emph{SG} and \emph{ASG}) an unfair advantage (by allowing them to choose the best step-size in hindsight), the SG methods always do substantially better than the FG methods (\emph{AFG} and \emph{L-BFGS}) on the first few passes through the data. However, the SG methods typically make little progress after the first few passes. In contrast, the FG methods make steady progress and eventually the faster FG method (\emph{L-BFGS}) typically passes the SG methods.
\item {\bf (FG and SG) vs.~SAG}: The SAG iterations seem to achieve the best of both worlds. They start out substantially better than FG methods, often obtaining similar performance to an SG method with the best step-size chosen in hindsight. But the SAG iterations continue to make steady progress even after the first few passes through the data. This leads to better performance than SG methods on later iterations, and on most data sets the sophisticated FG methods do not catch up to the SAG method even after $50$ passes through the data.
\item {\bf IAG vs.~SAG}: Even though these two algorithms differ in only the seemingly-minor detail of selecting data points at random (SAG) compared to cycling through the data (IAG), the performance improvement of SAG over its deterministic counterpart IAG is striking (even though the IAG method was allowed to choose the best step-size in hindsight). We believe this is due to the larger step sizes allowed by the SAG iterations, which would cause the IAG iterations to diverge.
\end{list}

\subsection{Comparison to Coordinate Optimization Methods}

For the $\ell_2$-regularized logistic regression problem, an alternative means to take advantage of the structure of the problem and achieve a linear convergence rate with a cheaper iteration cost than FG methods is to use randomized coordinate optimization methods. In particular, we can achieve a linear convergence rate by applying coordinate descent to the primal~\citep{nesterov2010efficiency} or coordinate-ascent to the dual~\citep{schwartz12}. In our second experiment, we included the following additional methods in this comparison:
\begin{list}{\labelitemi}{\leftmargin=1.7em}
\item {\bf PCD}: The randomized primal coordinate-descent method of \auth{Nesterov}\citet{nesterov2010efficiency}, using a step-size of $1/L_j$, where $L_j$ is the Lipschitz-constant with respect to coordinate $j$ of $g'$. Here, we sampled the coordinates uniformly.
  \item {\bf PCD-L}: The same as above, but sampling coordinates according to their Lipschitz constant, which can lead to an improved convergence rate~\citep{nesterov2010efficiency}.
  \item {\bf DCA}: Applying randomized coordinate ascent to the dual, with uniform example selection and an exact line-search~\citep{schwartz12}.
\end{list}
As with comparing FG and SG methods, it is difficult to compare coordinate-wise methods to FG and SG methods in an implementation-independent way. To do this in a way that we believe is fair (when discussing convergence rates), we measure the number of effective passes of the \emph{DCA} method as the number of iterations of the method divided by $n$ (since each iteration accesses a single example as in SG and SAG iterations). We measure the number of effective passes for the \emph{PCD} and \emph{PCD-L} methods as the number of iterations multiplied by $n/p$ so that $1$ effective pass for this method has a cost of $O(np)$ as in FG and SG methods. We ignore the additional cost associated with the Lipschitz sampling for the \emph{PCD-L} method (as well as the expense incurred because the \emph{PCD-L} method tended to favour updating the bias variable for sparse data sets) and we also ignore the cost of numerically computing the optimal step-size for the \emph{DCA} method.

\begin{figure}
\centering
\mbox{ \includegraphics[width=\figSize]{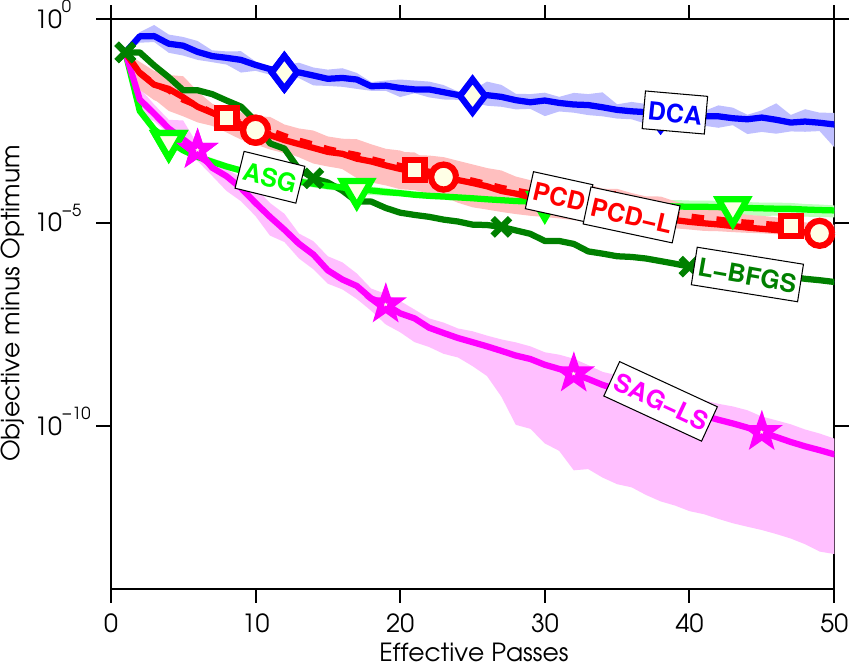} \hspace*{-.1cm}
\includegraphics[width=\figSize]{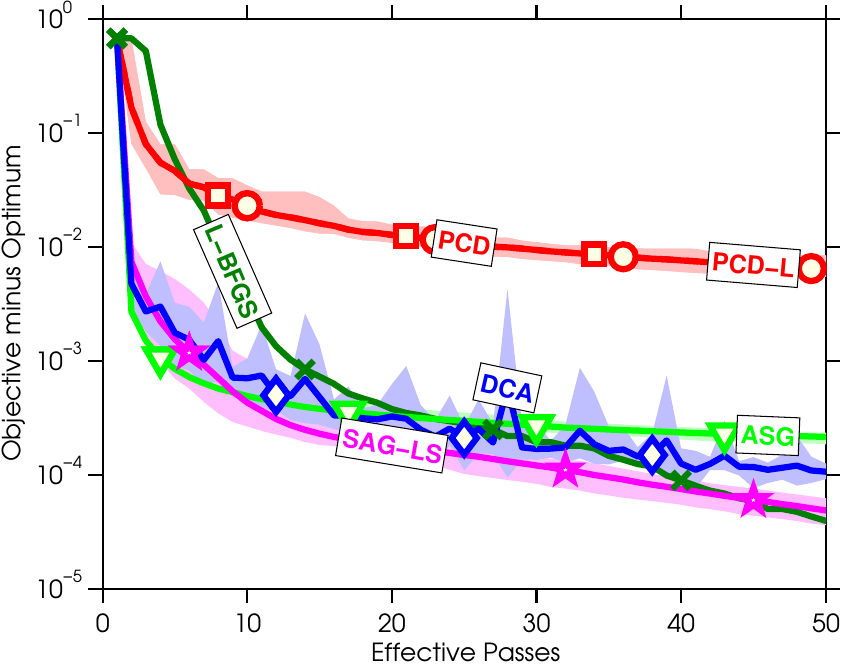} \hspace*{-.1cm}
\includegraphics[width=\figSize]{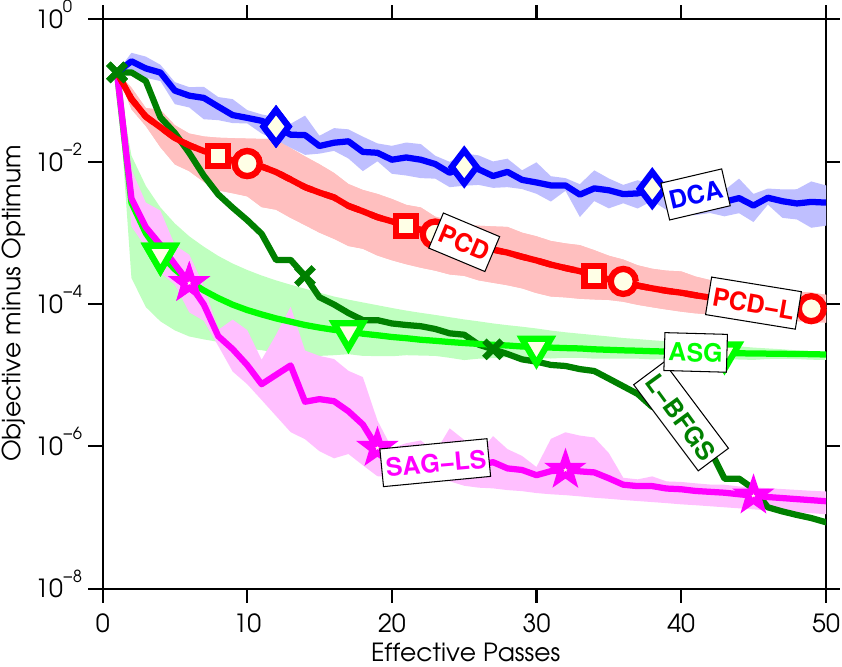} \hspace*{-.1cm}}

\mbox{\includegraphics[width=\figSize]{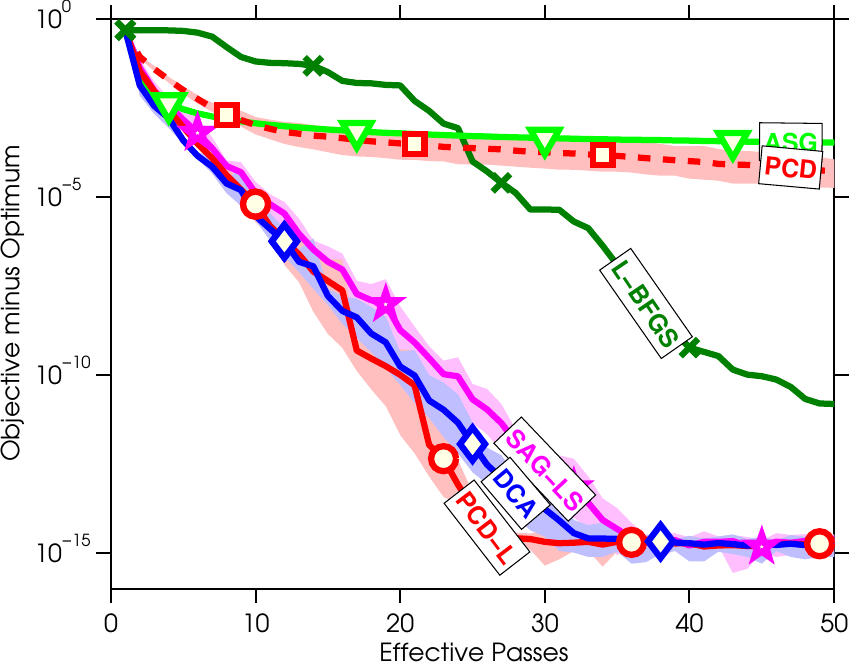} \hspace*{-.1cm}
\includegraphics[width=\figSize]{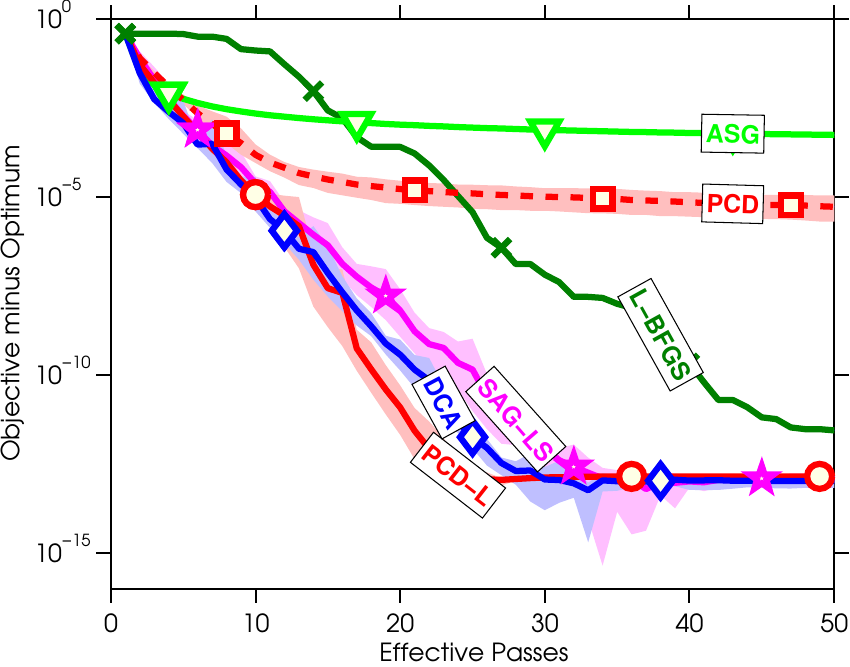} \hspace*{-.1cm}
\includegraphics[width=\figSize]{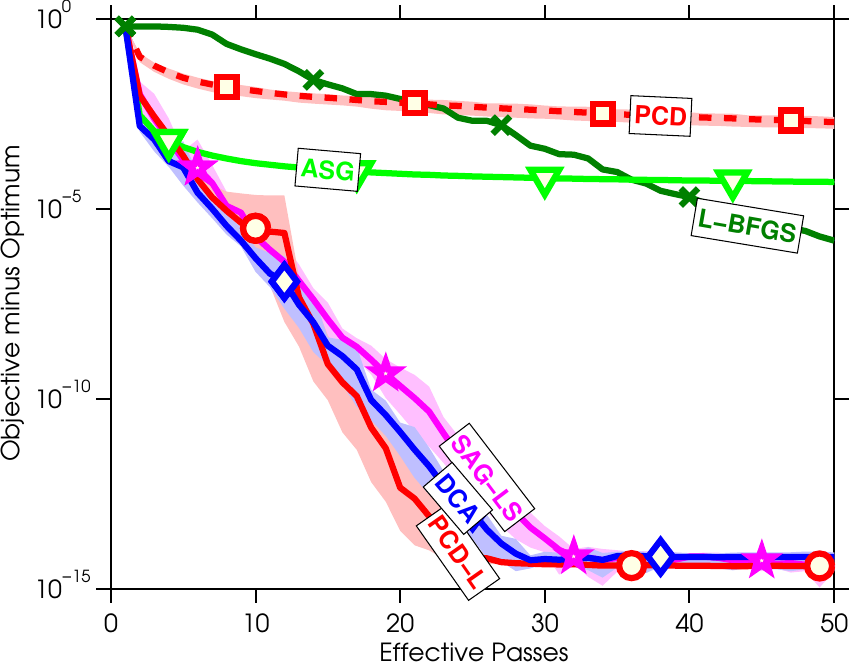} \hspace*{-.1cm}
}

\mbox{\includegraphics[width=\figSize]{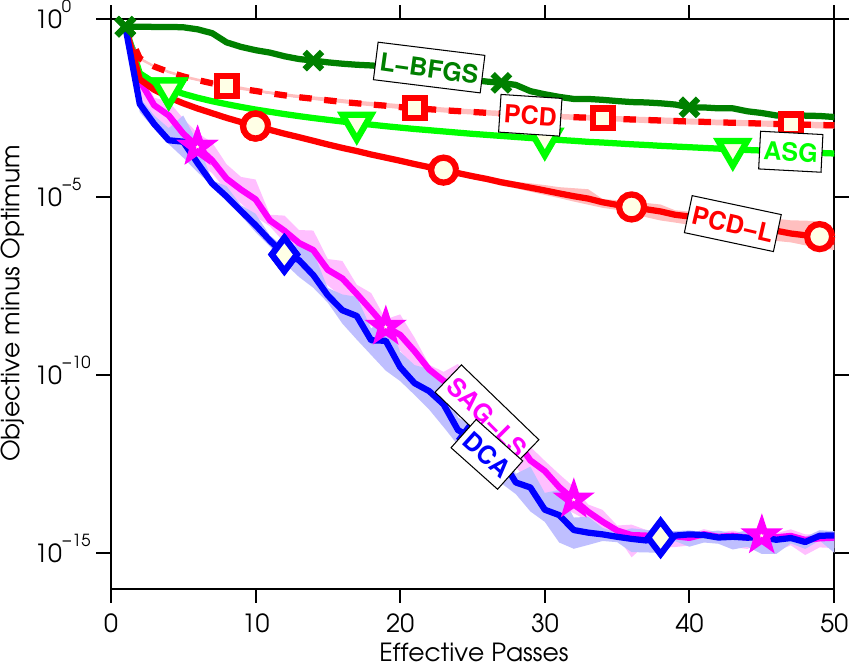} \hspace*{-.1cm}
\includegraphics[width=\figSize]{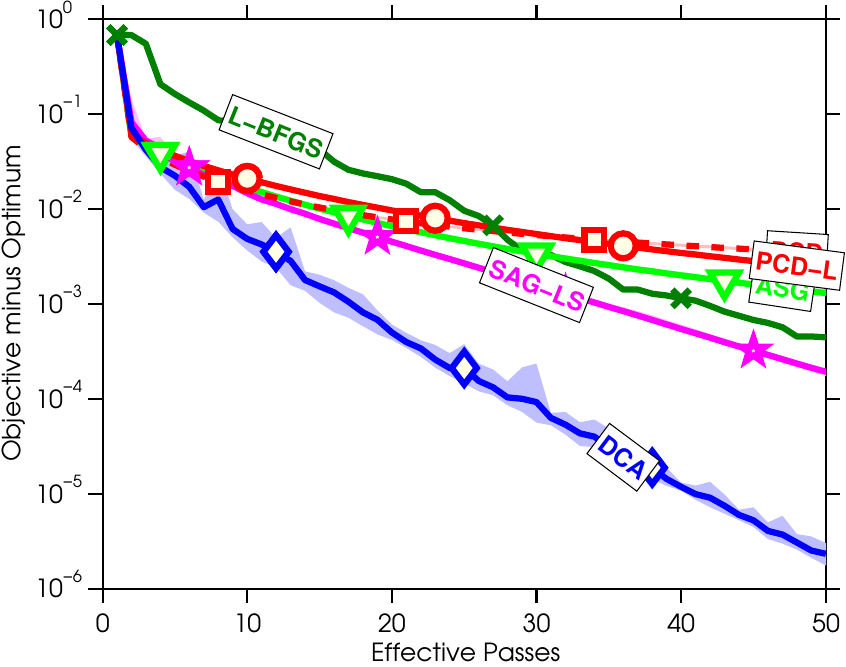} \hspace*{-.1cm}
\includegraphics[width=\figSize]{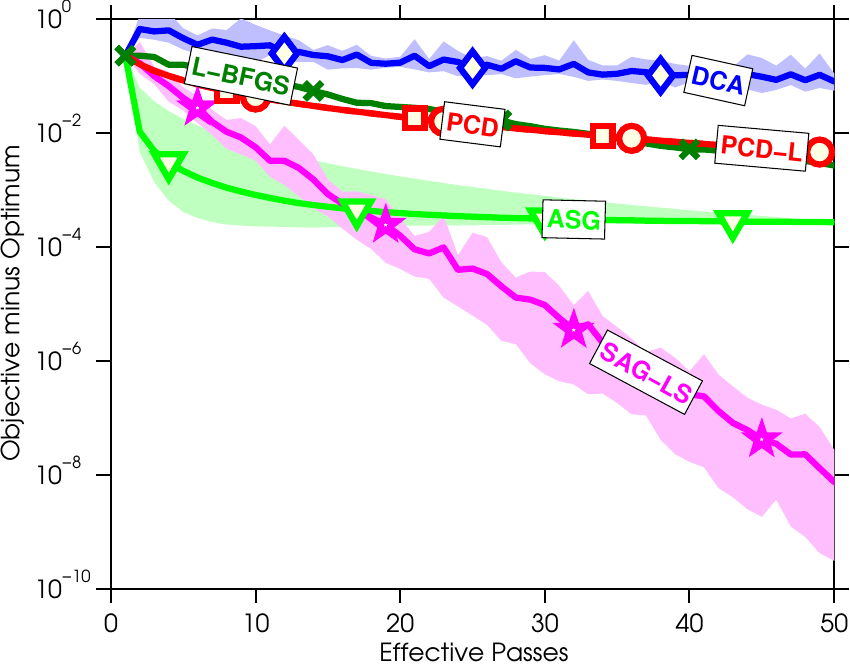} \hspace*{-.1cm}
}
\caption{Comparison of optimization different FG and SG methods to coordinate optimization methods.The top row gives results on the \emph{quantum} (left), \emph{protein} (center) and \emph{covertype} (right) datasets. The middle row gives results on the \emph{rcv1} (left), \emph{news} (center) and \emph{spam} (right) datasets.  The bottom row gives results on the \emph{rcv1Full} (left), \emph{sido} (center), and \emph{alpha} (right) datasets. This figure is best viewed in colour.}
\label{fig:CD}
\end{figure}

We compare the performance of the randomized coordinate optimization methods to several of the best methods from the previous experiment in Figure~\ref{fig:CD}. In these experiments we observe the following trends:
\begin{list}{\labelitemi}{\leftmargin=1.7em}
  \item {\bf PCD vs.~PCD-L}: For the problems with $n > p$ (top and bottom rows of Figure~\ref{fig:CD}), there is little difference between uniform and Lipschitz sampling of the coordinates. For the problems with $p > n$ (middle row of Figure~\ref{fig:CD}), sampling according to the Lipschitz constant leads to a large performance improvement over uniform sampling.
  \item {\bf PCD vs.~DCA}: For the problems with $p > n$, \emph{DCA} and \emph{PCD-L} have similar performance. For the problems with $n > p$, the performance of the methods typically differed but neither strategy tended to dominate the other.
  \item {\bf (PCD and DCA) vs.~(SAG)}: For some problems, the \emph{PCD} and \emph{DCA} methods have performance that is similar to \emph{SAG-LS} and the \emph{DCA} method even gives better performance than \emph{SAG-LS} on one data set. However, for many data sets either the \emph{PCD-L} or the \emph{DCA} method (or both) perform poorly while the \emph{SAG-LS} iterations are among the best or substantially better than all other methods on every data set. This suggests that, while coordinate optimization methods are clearly extremely effective for some problems, the SAG method tends to be a more robust choice across problems.
  \end{list}

\iftoggle{springer}
{

}
{
\subsection{Comparison of Step-Size Strategies} 
  
In our prior work we analyzed the step-sizes $\alpha_k = 1/2nL$ and $\alpha_k = 1/2n\mu$~\citep{roux2012stochastic}, while Section~\ref{convergence} considers the choice $\alpha_k = 1/16L$ and Section~\ref{sec:stepSizes} discusses the choices $\alpha_k = 1/L$ and $\alpha_k = 2/(L + n\mu)$ as in FG methods. In Figure~\ref{fig:stepSize} we compare the performance of these various strategies to the performance of the SAG algorithm with our proposed line-search as well as the IAG and SAG algorithms when the best step-size is chosen in hindsight.
In this plot we see the following trends:
\begin{list}{\labelitemi}{\leftmargin=1.7em}
  \item {\bf Proposition 1 of \citet{roux2012stochastic}}: Using a step-size of $\alpha_k = 1/2nL$ performs poorly, and makes little progress compared to the other methods. This makes sense because Proposition~1 in~\citet{roux2012stochastic} implies that the convergence rate (in terms of effective passes through the data) under this step size will be similar to the basic gradient method, which is known to perform poorly unless the problem is very well-conditioned.
  \item {\bf Proposition 2 of \citet{roux2012stochastic}}: Using a step-size of $\alpha_k = 1/2n\mu$ performs extremely well on the data sets with $p > n$ (middle row). In contrast, for the data sets with $n > p$ it often performs very poorly, and in some cases appears to diverge. This is consistent with Proposition~2 in~\citet{roux2012stochastic}, which shows a fast convergence rate under this step size only if certain conditions on $\{n,\mu,L\}$ hold.
\item {\bf Theorem~\ref{thm}}: Using a step-size of $\alpha_k = 1/16L$ performs consistently better than the smaller step size $\alpha_k = 1/2nL$, but in some cases it performs worse than $\alpha_k = 1/2n\mu$. However, in contrast to $\alpha_k = 1/2n\mu$, the step size $\alpha_k = 1/16L$ always has reasonable performance.
  \item {\bf Section~\ref{sec:stepSizes}}: The step size of $\alpha_k = 1/L$ performs performs extremely well on the data sets with $p > n$, and performs better than the step sizes discussed above on all but one of the remaining data sets. The step size of $\alpha_k = 2/(L + n\mu)$ seems to perform the same or slightly better than using $\alpha_k = 1/L$ except on one data set where it performs poorly.
  \item {\bf Line-Search}: Using the line-search from Section~\ref{sec:linesearch} tends to perform as well or better than the various constant step size strategies, and tends to have similar performance to choosing the best step size in hindsight.
  \item {\bf IAG vs.~SAG}: When choosing the best step size in hindsight, the SAG iterations tend to choose a much larger step size than the IAG iterations. The step sizes chosen for SAG were $100$ to $10000$ times larger than the step sizes chosen by IAG, and always lead to better performance by several orders of magnitude.
  \end{list}

\begin{figure}
\centering
\mbox{ \includegraphics[width=\figSize]{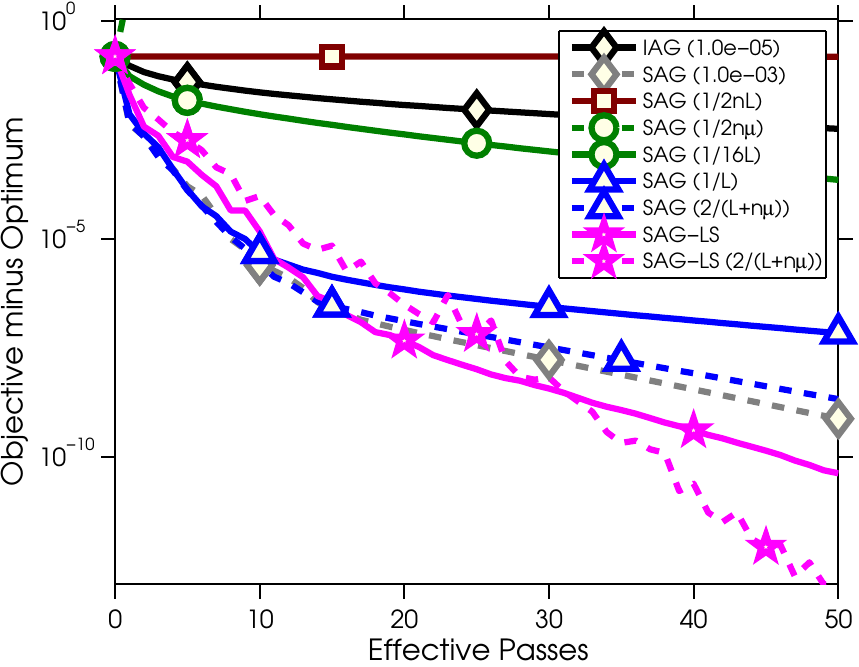} \hspace*{-.1cm}
\includegraphics[width=\figSize]{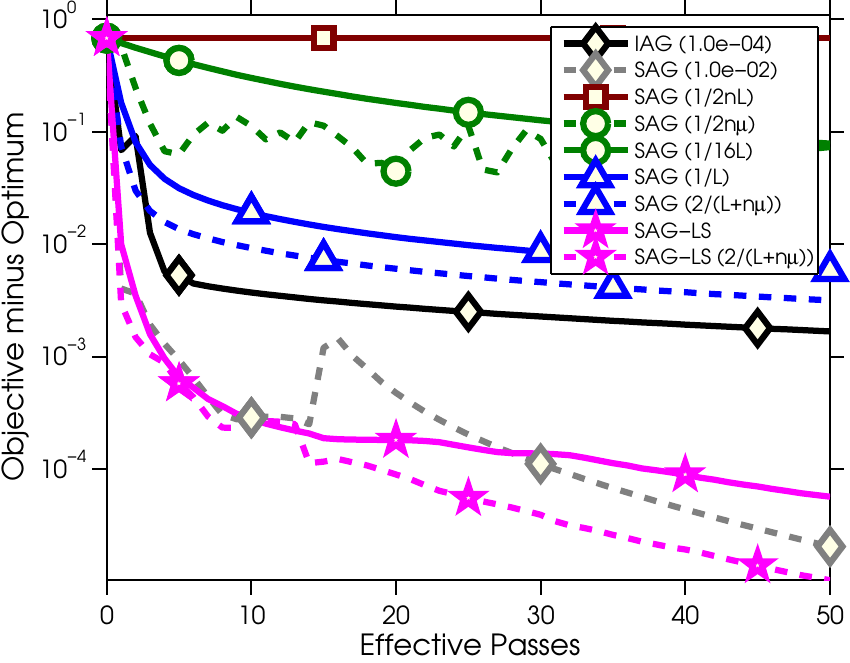} \hspace*{-.1cm}
\includegraphics[width=\figSize]{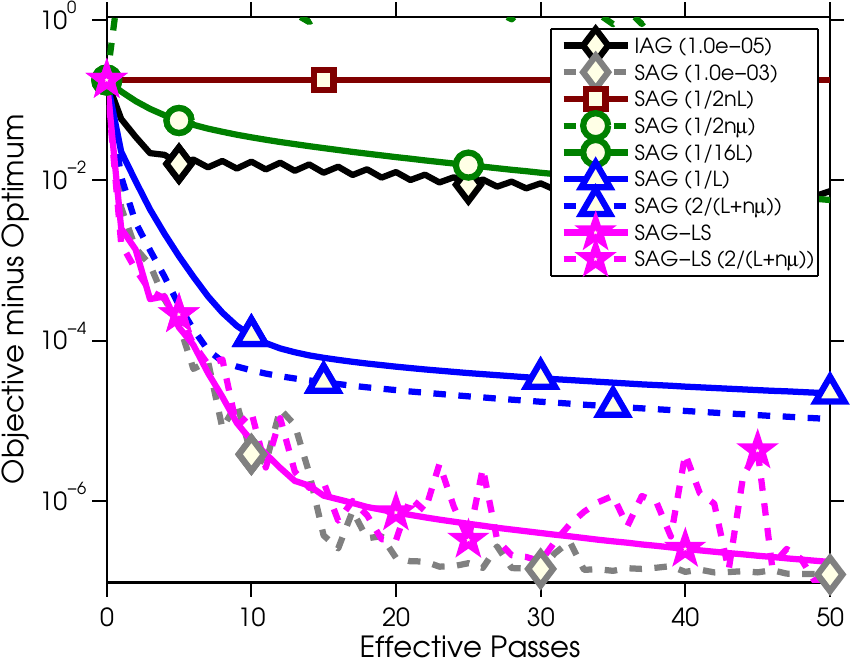} \hspace*{-.1cm}}

\mbox{\includegraphics[width=\figSize]{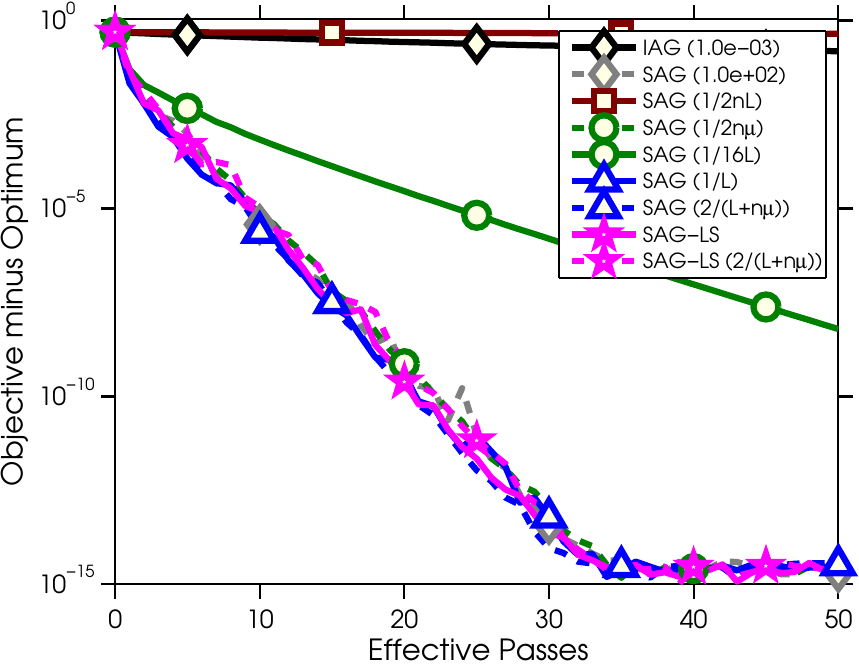} \hspace*{-.1cm}
\includegraphics[width=\figSize]{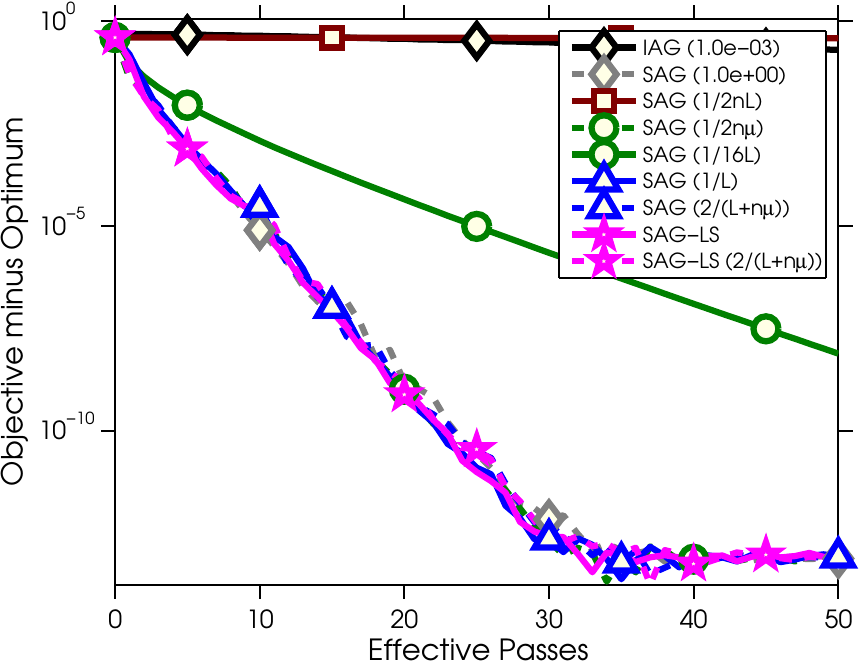} \hspace*{-.1cm}
\includegraphics[width=\figSize]{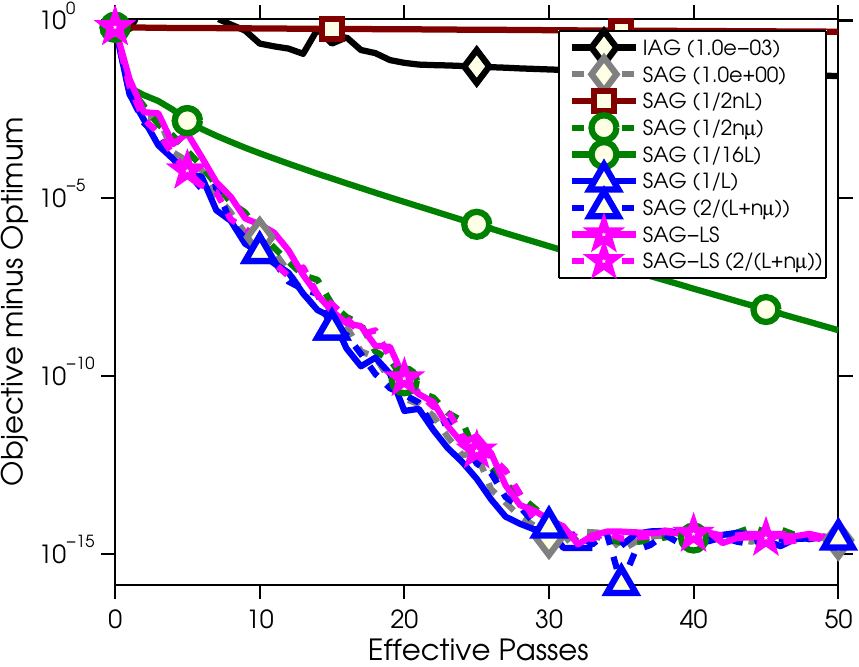} \hspace*{-.1cm}
}

\mbox{\includegraphics[width=\figSize]{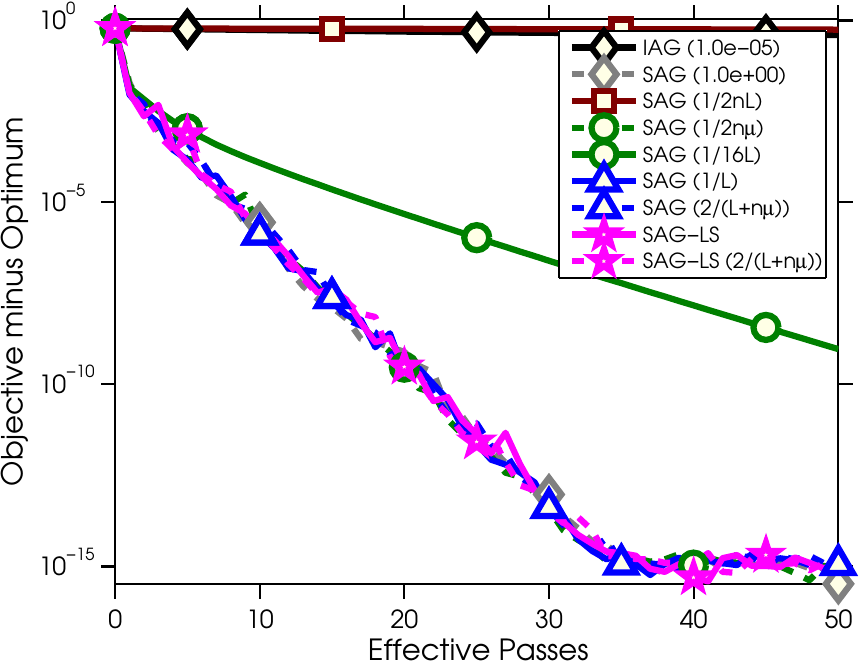} \hspace*{-.1cm}
\includegraphics[width=\figSize]{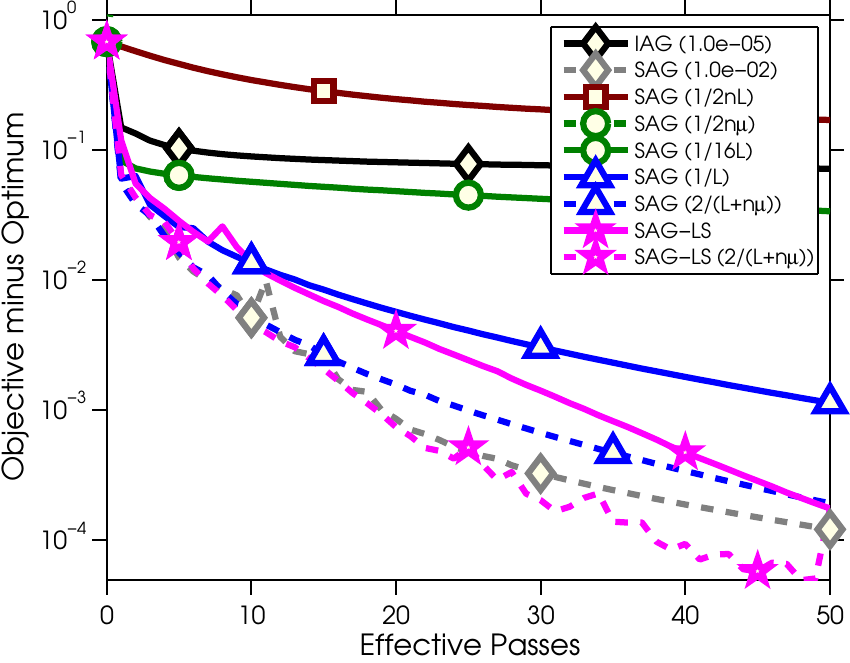} \hspace*{-.1cm}
\includegraphics[width=\figSize]{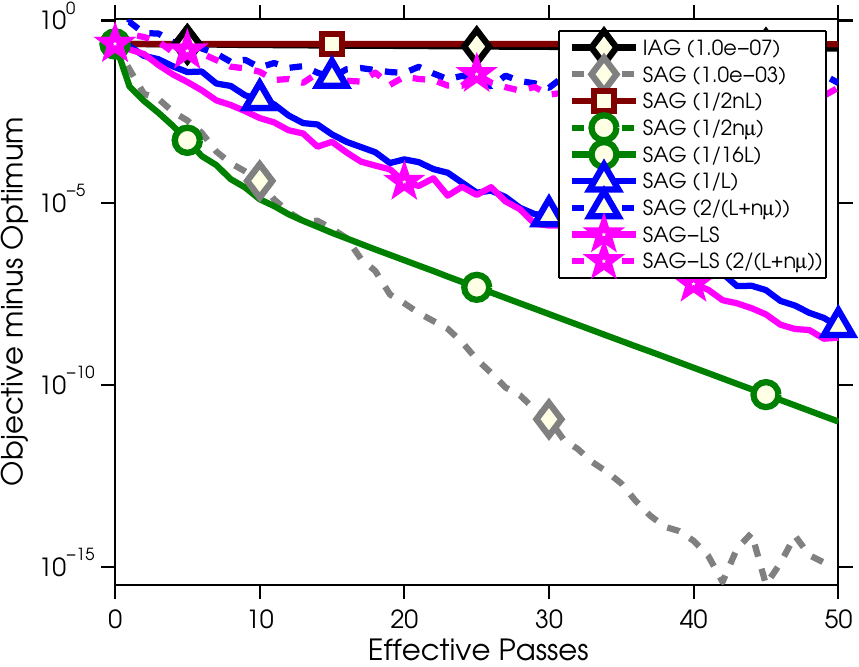} \hspace*{-.1cm}
}
\caption{Comparison of step size strategies for the SAG method. The top row gives results on the \emph{quantum} (left), \emph{protein} (center) and \emph{covertype} (right) datasets. The middle row gives results on the \emph{rcv1} (left), \emph{news} (center) and \emph{spam} (right) datasets.  The bottom row gives results on the \emph{rcv1Full} (left), \emph{sido} (center), and \emph{alpha} (right) datasets. This figure is best viewed in colour.}
\label{fig:stepSize}
\end{figure}

\subsection{Effect of mini-batches}
\label{exp:mini}

As we discuss in Section~\ref{sec:miniBatch}, when using mini-batches within the SAG iterations there is a trade-off between the higher iteration cost of using mini-batches and the faster convergence rate obtained using mini-batches due to the possibility of using a smaller value of $L$. In Figure~\ref{fig:mini_batches}, we compare (on the dense data sets) the excess sub-optimality as a function of the number of examples seen for various mini-batch sizes and the three step-size strategies $1/L_\textrm{max}$, $1/L_\textrm{mean}$, and $1/L_\textrm{Hessian}$ discussed in Section~\ref{sec:miniBatch}.

Several conclusions may be drawn from these experiments:
\begin{list}{\labelitemi}{\leftmargin=1.7em}
\item Even though Theorem~\ref{thm} hints at a maximum mini-batch size of $\frac{n\mu}{2L}$ without loss of convergence speed, this is a very conservative estimate. In our experiments, the original value of $\frac{n\mu}{L}$ was on the order of $10^{-5}$ and mini-batch sizes of up to 500 could be used without a loss in performance. Not only does this yield large memory storage gains, it would increase the computational efficiency of the algorithm when taking into account vectorization.
\item To achieve fast convergence, it is essential to use a larger step-size when larger mini-batches are used. For instance, in the case of the \emph{quantum} dataset with a mini-batch size of $20000$, we have $\frac{1}{L_{\textrm{max}}} = 4.4 \cdot 10^{-4}$, $\frac{1}{L_{\textrm{mean}}} = 5.5 \cdot 10^{-2}$ and $\frac{1}{L_{\textrm{Hessian}}} = 3.7 \cdot 10^{-1}$. In the case of the \emph{covertype} dataset and for a mini-batch size of $20000$, we have $\frac{1}{L_{\textrm{max}}} = 2.1 \cdot 10^{-5}$, $\frac{1}{L_{\textrm{mean}}} = 6.2 \cdot 10^{-2}$ and $\frac{1}{L_{\textrm{Hessian}}} = 4.1 \cdot 10^{-1}$. 
\end{list}

\begin{figure}
\begin{center}
\includegraphics[width=\textwidth]{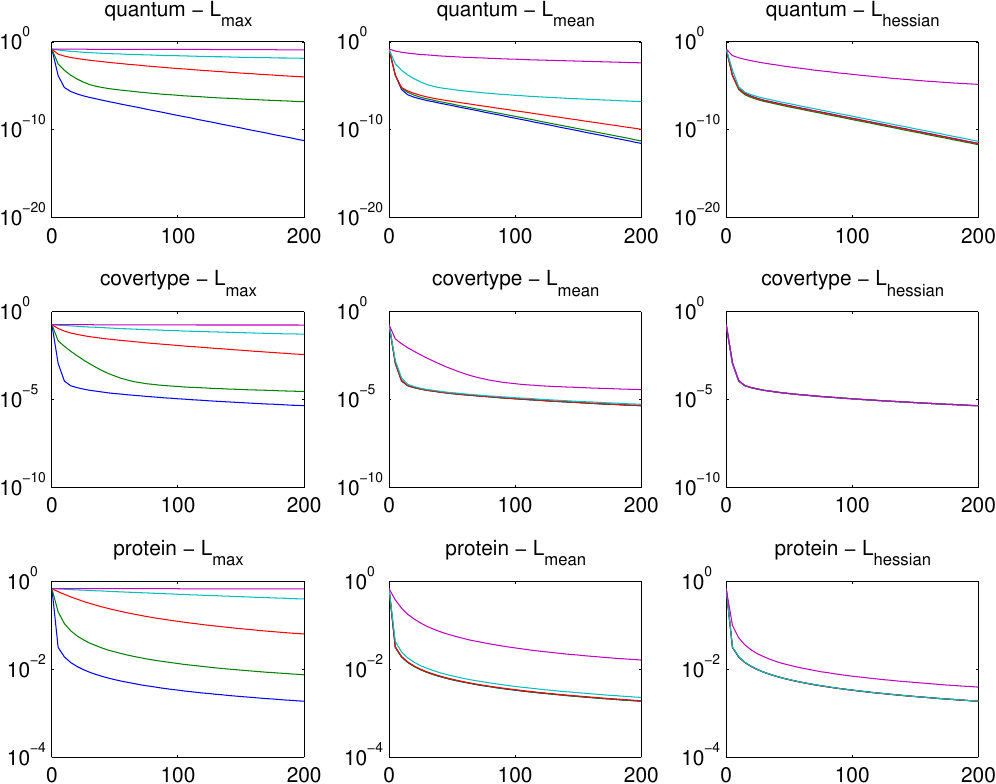}\\
\vspace*{.4cm}
\includegraphics[width=.85\textwidth]{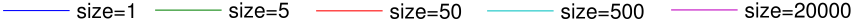}
\end{center}
\caption{Sub-optimality as a function of the number of effective passes through the data for various datasets, step-size selection schemes and mini-batch sizes. The datasets are \emph{quantum} (top), \emph{covertype} (middle) and \emph{protein} (bottom). Left: the step-size is $1/L$ with $L$ the maximum Lipschitz constant of the individual gradients. It is thus the same for all mini-batch sizes. Center: the step-size is $1/L$ where $L$ is obtained by taking the maximum among the averages of the Lipschitz constants within mini-batches. Right: the step-size is $1/L$ where $L$ is obtained by computing the maximum eigenvalue of the Hessian for each mini-batch, then taking the maximum of these quantities across mini-batches.
}
\label{fig:mini_batches}
\end{figure}
}

\subsection{Effect of non-uniform sampling}
\label{exp:lipschitz}

\begin{figure}
\centering
\mbox{ \includegraphics[width=\figSize]{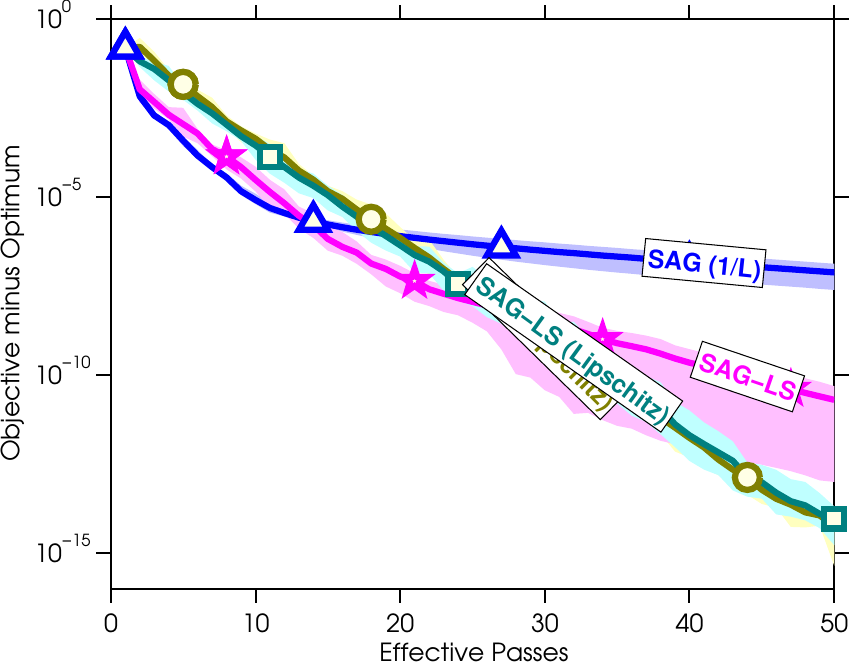} \hspace*{-.1cm}
\includegraphics[width=\figSize]{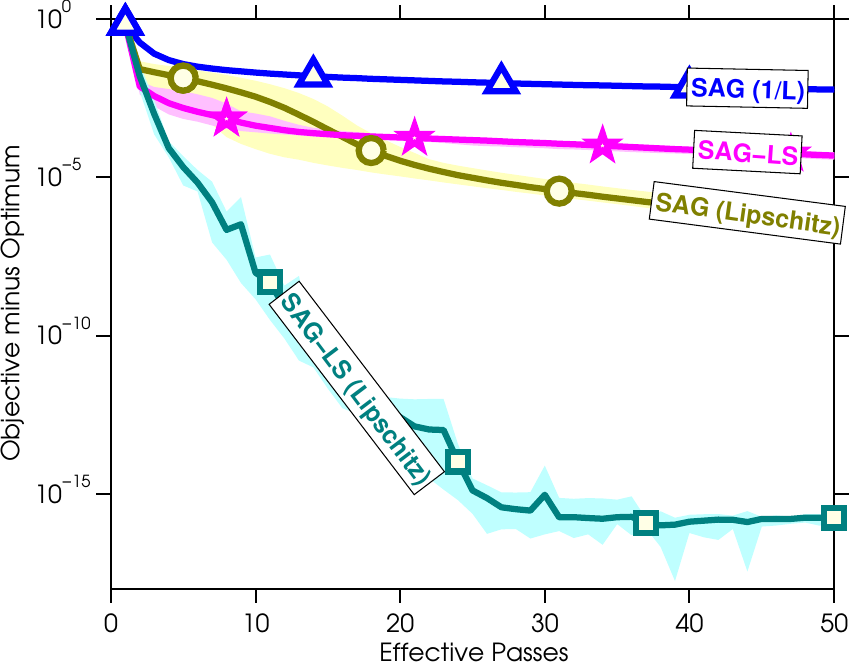} \hspace*{-.1cm}
\includegraphics[width=\figSize]{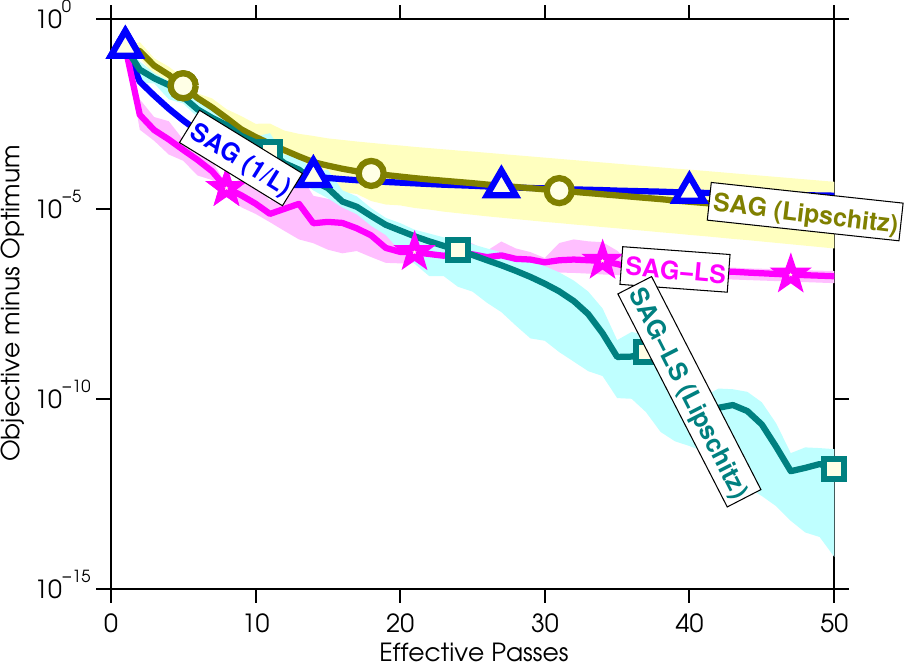} \hspace*{-.1cm}}

\mbox{\includegraphics[width=\figSize]{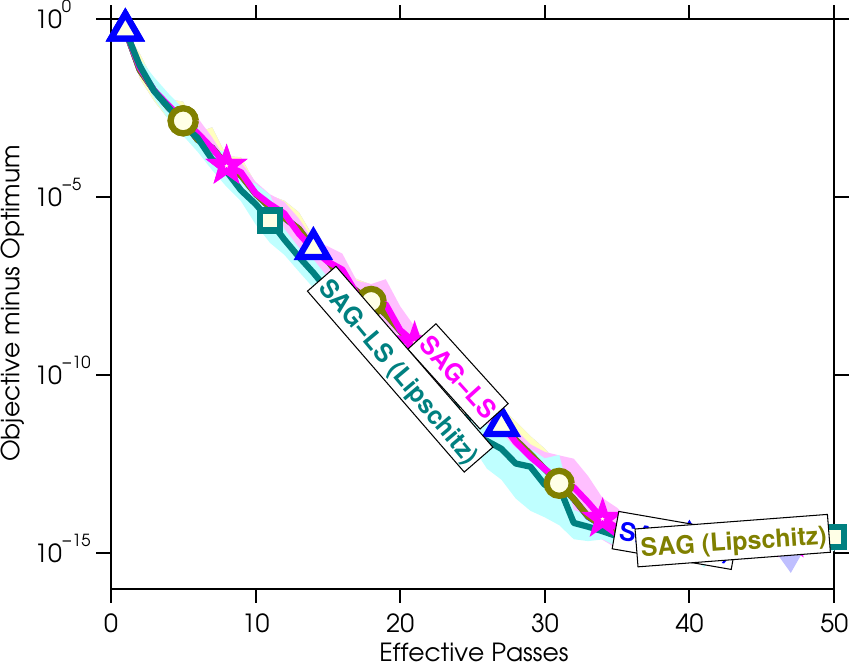} \hspace*{-.1cm}
\includegraphics[width=\figSize]{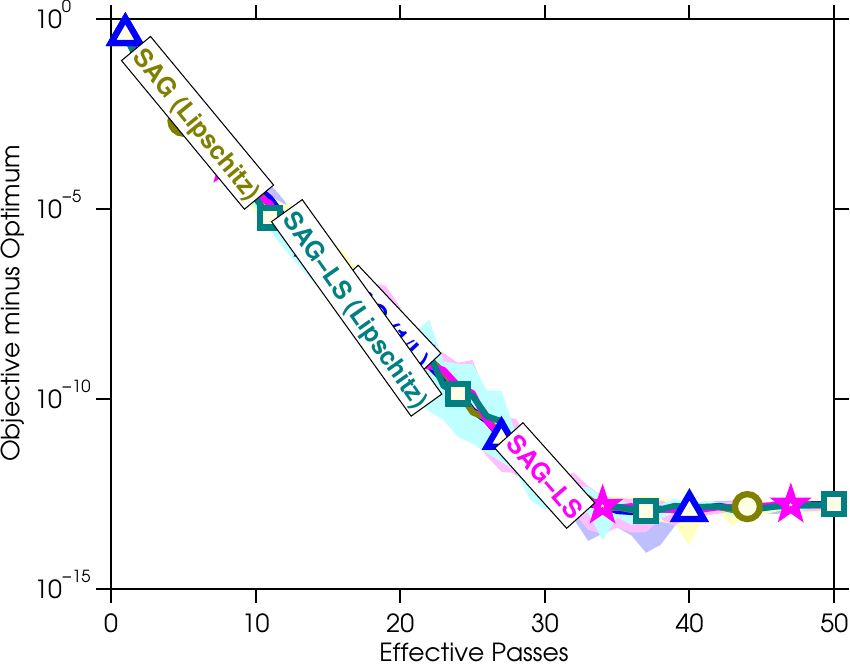} \hspace*{-.1cm}
\includegraphics[width=\figSize]{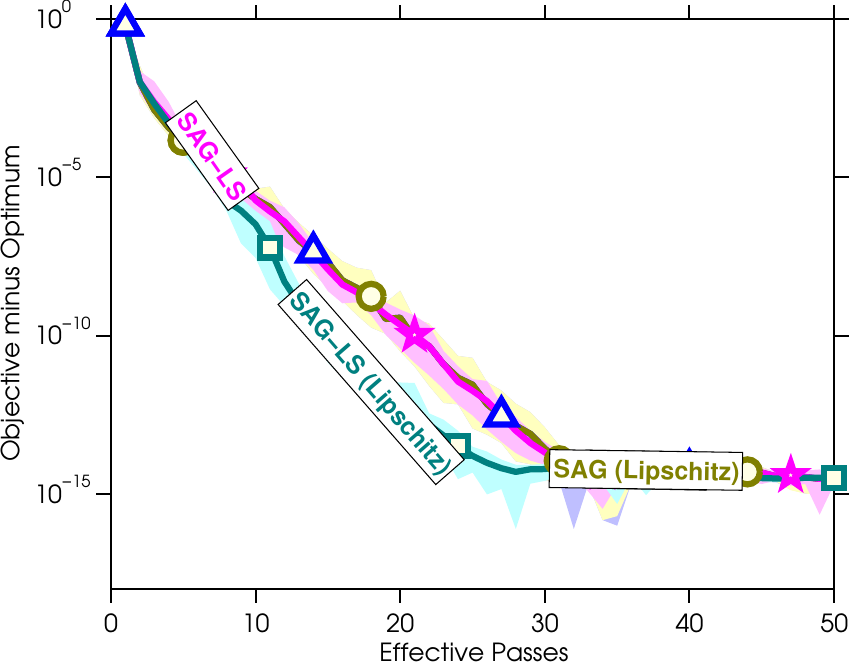} \hspace*{-.1cm}
}

\mbox{\includegraphics[width=\figSize]{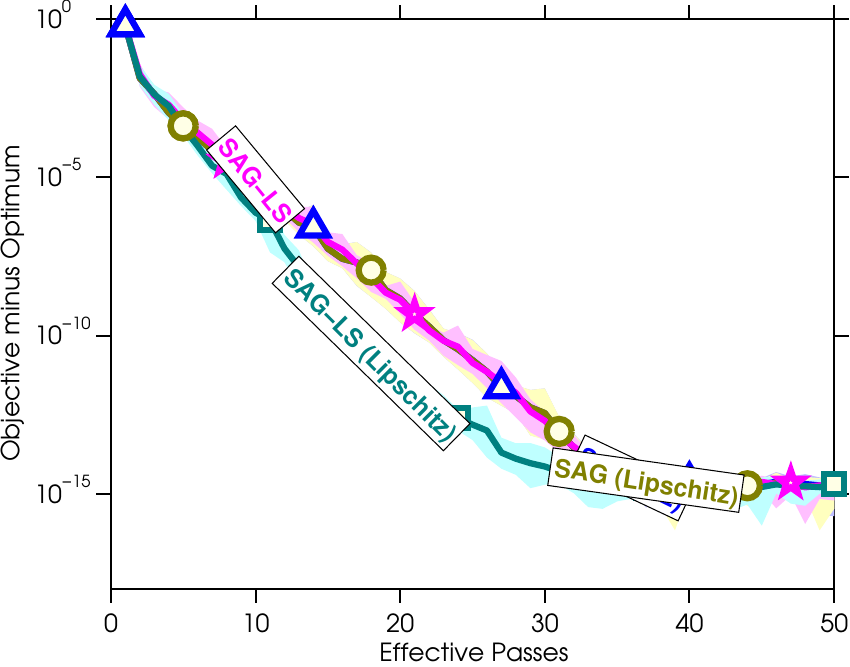} \hspace*{-.1cm}
\includegraphics[width=\figSize]{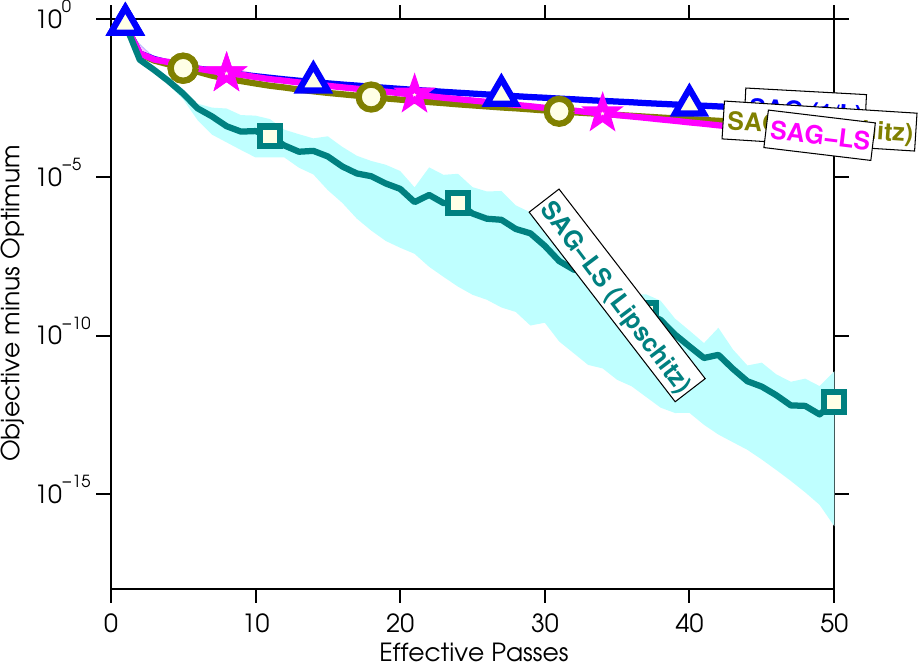} \hspace*{-.1cm}
\includegraphics[width=\figSize]{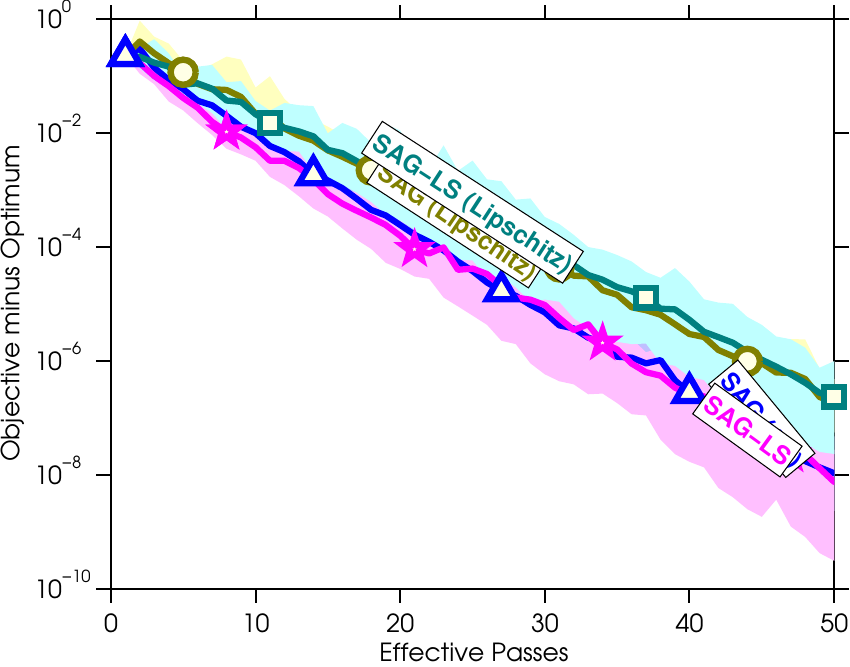} \hspace*{-.1cm}}
\caption{Comparison of uniform and non-uniform sampling strategies for the SAG algorithm. The top row gives results on the \emph{quantum} (left), \emph{protein} (center) and \emph{covertype} (right) datasets. The middle row gives results on the \emph{rcv1} (left), \emph{news} (center) and \emph{spam} (right) datasets.  The bottom row gives results on the \emph{rcv1Full}| (left), \emph{sido} (center), and \emph{alpha} (right) datasets. This figure is best viewed in colour.}
\label{fig:lipschitz}
\end{figure}

In our final experiment, we explored the effect of using the non-uniform sampling strategy discussed in Section~\ref{sec:lipschitz}.  In Figure~\ref{fig:lipschitz}, we compare several of the SAG variants with uniform sampling to the following two methods based on non-uniform sampling:
\begin{list}{\labelitemi}{\leftmargin=1.7em}
  \item {\bf SAG (Lipschitz)}: In this method we sample the functions in proporition to $L_i + c$, where $L_i$ is the global Lipschitz constant of the corresponding $f_i'$ and we set $c$ to the average of these constants, $c=L_\textrm{mean}=(1/n)\sum_i L_i$. Plugging in these values into the formula at the end of Section~\ref{sec:lipschitz} and using $L_\textrm{max}$ to denote the maximum $L_i$ value, we set the step-size to $\alpha_k = 1/2L_\textrm{max} + 1/2L_\textrm{mean}$.
  \item {\bf SAG-LS (Lipschitz)}: In this method we formed estimates of the quantities $L_i$, $L_\textrm{max}$, and $L_\textrm{mean}$. The estimator $L^k_\textrm{max}$ is computed in the same way as the SAG-LS method. To estimate each $L_i$, we keep track of an estimate $L_i^k$ for each $i$ and we set $L^k_\textrm{mean}$ to the average of the $L_i^k$ values among the $f_i$ that we have sampled at least once. We set $L_i^k = L_i^{k-1}$ if example $i$ was not selected and otherwise we initialize to $L_i^k = L_i^{k-1}/2$ and perform the line-search until we have a valid $L_i^k$ (this means that $L^k_\textrm{mean}$ will be approximately halved if we perform a full pass through the data set and never violate the inequality). To ensure that we do not ignore important unseen data points for too long, in this method we sample a previously unseen function with probability $(n-m)/n$, and otherwise we sample from the previously seen $f_i$ in proportion to $L_i^k + L^k_\textrm{mean}$. To prevent relying too much on our initially-poor estimate of $L_\textrm{mean}$, we use a step size of $\alpha_k = \frac{n-m}{n}\alpha_\textrm{max} + \frac{m}{n}\alpha_\textrm{mean}$, where $\alpha_\textrm{max} = 1/L^k_\textrm{max}$ is the step-size  we normally use with uniform sampling and $\alpha_\textrm{mean} = 1/2L^k_\textrm{max} + 1/2L^k_\textrm{mean}$ is the step-size we use with the non-uniform sampling method, so that the method interpolates between these extremes until the entire data set has been sampled.
\end{list}
We make the following observations from these experiments:
\begin{list}{\labelitemi}{\leftmargin=1.7em}
  \item {\bf SAG (1/L) vs.~SAG (Lipschitz)}: With access to global quantities and a constant step size, the difference between uniform and non-uniform sampling was typically quite small. However, in some cases the non-uniform sampling method behaved much better (top row of Figure~\ref{fig:lipschitz}).
  \item {\bf SAG-LS vs.~SAG-LS (Lipschitz)}: When estimating the Lipschitz constants of the individual functions, the non-uniform sampling strategy often gave better performance. Indeed, the adaptive non-uniform sampling strategy gave solutions that are orders of magnitude more accurate than any method we examined for several of the data sets (e.g., the \emph{protein}, \emph{covertype}, and \emph{sido} data sets) In the context of logistic regression, it makes sense that an adaptive sampling scheme could lead to better performance, as many correctly-classified data samples might have a very slowly-changing gradient near the solution, and thus they do not need to be sampled often.
  \end{list}

\section{Discussion}
\label{discussion}

\iftoggle{springer}
{
Since the first version of this work was published~\citep{roux2012stochastic}, there has been an explosion of interest in stochastic methods with improved convergence rates. In Section~6 of the extended version of this paper we review other algorithms that have been discovered to have similar convergence rates, including stochastic dual coordinate ascent (SDCA)~\citep{schwartz12}, incremental surrogate optimization (MISO)~\citep{mairal2013surrogate}, stochastic variance-reduced gradient (SVRG) methods~\citep{mahdavi2013mixedGrad,johnson2013accelerating,konevcny2013semi} which remove the memory requirement, and the SAGA method~\citep{defazio2014saga}. In Section~6 of the extended version we also review many of the possible variants on these basic algorithms that have been explored. This includes accelerated gradient methods that achieve faster convergence rates for ill-conditioned problems~\citep{lin2015universal}, proximal-gradient and alternating direction method of multipliers (ADMM) variants that can solve certain constrained and non-smooth problems~\citep{defazio2014saga,zhong2013fastADMM}, coordinate-wise variants that only update a subset of the variables on each iteration~\citep{konevcny2014semi},  Newton-like variants of the method~\cite{sohl2014fast}, methods where non-uniform sampling provably improves the convergence rate~\citep{schmidt2015sag4crf}, and analyses that give a linear convergence rate without strong convexity~\citep{gong2014linear}.

}
{

Since the first version of this work was published~\citep{roux2012stochastic}, there has been an explosion of interest in stochastic methods with improved convergence rates. In this section we first review other algorithms that have been discovered to have this property, and then we discuss the many possible variants on these basic algorithms that have been explored. As this is a very quickly-evolving area there are likely to be many new developments in the near future, but we note that this literature review is up to date as of January, 2015.

\subsection{Alternative Algorithms }

\textbf{SDCA}: The first algorithm that was shown to have a similar convergence rate and iteration cost to SAG was in fact a much older algorithm: coordinate optimization applied to a dual problem with randomized coordinate selection, referred to as stochastic dual coordinate ascent (SDCA). \auth{Shalev-Shwartz and Zhang}\citet{schwartz12} consider the problem of minimizing an $\ell_2$-regularized finite sum,
\[
\minimize{x\in\Real^p}\frac{\lambda}{2}\norm{x}^2 + \frac{1}{n}\sum_{i=1}^n f_i(a_i^Tx),
\]
where each $f_i$ is convex and each $f_i'$ is Lipschitz-continuous, by optimizing its Fenchel dual,
\[
\maximize{y\in\Real^n}-\frac{1}{2\lambda}\norm{\frac{1}{n}\sum_{i=1}^n y_ia_i}^2 - \frac{1}{n}\sum_{i=1}^nf_i^*(-y_i).
\]
They consider applying exact coordinate optimization to a randomly selected coordinate. Using the primal-dual relationship $x = \frac{1}{n\lambda}\sum_{i=1}^Nx_ia_i$, they show a linear convergence rate in terms of the duality gap. Since there is one dual variable associated with each example $i$, the iteration cost is independent of $n$ and thus the strategy has similar convergence properties to SAG (the memory requirements are identical in this context, see Section~\ref{sec:sparse}). 

This result is related to the earlier result of \auth{Nesterov}\citep{nesterov2010efficiency}, who shows a similar convergence rate for randomized coordinate descent.\footnote{Local linear convergence rates of deterministic coordinate descent methods had been established much earlier~\citep{luo1992convergence}.} However, the result of Nesterov cannot be directly applied to the dual problem in general since the dual does not necessarily have a Lipschitz-continuous gradient. A similar result was also reported even earlier by~\auth{Collins et al.}\citet{collins2008exponentiated} without requiring that the gradient of the dual is Lipschitz-continuous, but this result again only applied to the dual objective.
More recently, \auth{Shalev-Shwartz and Zhang}\citet{PSDCA} show linear convergence of SDCA in the more general setting
\[
\minimize{x\in\Real^p}\lambda r(x) + \frac{1}{n}\sum_{i=1}^n f_i(A_i^Tx),
\]
where $A_i$ are matrices and $r$ is $1-$strongly convex. They also relax the requirement of exact coordinate optimization, providing a variety of more practical alternatives. Further, in subsequent work they obtain a convergence rate in the convex case by adding an explicit strongly-convex regularizer to the problem~\citep{ASDCA}.

A disadvantage of SDCA compared to the SAG algorithm is that the SDCA convergence rates depend on $\lambda$ rather than the strong-convexity constant $\mu$. In the worst case we have $\mu=\lambda$, but if $\mu$ is much larger then the convergence rate of SAG is much faster. Further, even in cases where $\mu=\lambda$, the convergence rate of SAG might be much faster if the iterates stay in local region with a higher strong-convexity constant. As an extreme example, due to local strong-convexity SAG might have a linear convergence rate in scenarios where SDCA has a sub-linear convergence. This subtle but practically important issue was a key  focus in the recent work of \auth{Agarwal \& Bottou}\citep{argawal2014bounds}, and indeed the performance of SDCA was very poor on three of the test problems in our experiments. See Figure~\ref{fig:CD} (top left, top right, bottom right).

\textbf{MISO}: \auth{Mairal}\citet{mairal2013surrogate} analyzes a very general surrogate optimization framework, that includes a wide vareity of existing algorithms. He also considers incremental algorithms in this framework, and specialized to the smooth and unconstraeind setting (with a `Lipschitz surrogate') obtains an algorithm (MISO) that is very similar to the SAG algorithm,
\[
x^{k+1} = \frac{1}{n}\sum_{i=1}^{n}x^k_i - \frac{\alpha_k}{n}\sum_{i=1}^ny^k_i,
\]
where $y^k_i$ is defined as in the SAG algorithm~\ref{eq:yi}, and $x_i^k$ is the parameter vector used to compute the corresponding $y_i^k$. Thus, instead of applying the SAG step to $x^k$, MISO applies the step to the average of the previous $x_i^k$ values used to form the current $y_i^k$ variables. \auth{Mairal}\citet{mairal2013surrogate} shows that this algorithm also achieves an $O(1/k)$ rate for convex objectives and a linear convergence rate for strongly-convex objectives. 

However, MISO has the disadvantage that it not only requires storing the $n$ gradient values but also storing $n$ previous iterations (which are less likely to have a nice structure). Further, the linear convergence rate shown for MISO is substantially slower than the convergence rate shown in Theorem~\ref{thm}. In particular, the rate is more similar to the substantially slower Proposition~1 in our prior work~\citep{roux2012stochastic}. Subsequent work on the MISO algorithm has shown an analogous result to Proposition~2 in our prior work; if the $n$ is sufficiently larger than $L\mu$, then using a step-size proportional to $1/\mu$ yields a faster convergence rate~\citep{mairal2014incremental,defaziofinito}. However, using this step-size causes divergence if $\mu$ is not sufficiently large.

\textbf{SVRG}: Another interesting framework that has been considered is known as `mixed optimization', `stochastic variance-reduced gradient' (SVRG), or `semi-stochastic gradient descent' (S2GD)~\citep{mahdavi2013mixedGrad,johnson2013accelerating,zhang2013linear,konevcny2013semi}.\footnote{\citet{wang2013variance} consider a related algorithm that maintains an easily-computable approximation to the current gradient $g'(x^k)$. Although this can improve the constants in the sublinear SG convergence rates, it does not improve the rates.} Unlike FG methods that utilize the full gradient $g'$ and SG methods that consider individual gradients $f_i'$, these mixed optimization methods combine the two. In particular, they use a (possibly regularized) iteration of the form
\[
 x^{k+1} = x^k - \alpha_k(f_i'(x^k) - f_i'(\tilde{x}^k) + \frac{1}{n}\sum_{j=1}^nf_j'(\tilde{x}^k)),
\]
where $\tilde{x}^k$ is the last iterate where the full gradient $g'$ was evaluated. The algorithm alternates between computing the full gradient at $\tilde{x}^k$, and performing some number $m$ of stochastic gradient iterations.
Note this is very similar to the SAG algorithm written in the form
\[
 x^{k+1} = x^k - \alpha_k(\frac{1}{n}f_i'(x^k) - \frac{1}{n}f_i'(x_i^k) + \frac{1}{n}\sum_{j=1}^nf_j'(x_i^k)),
\]
where $x_i^k$ is defined as in the MISO algorithm. These methods are thus similar to SAG in the use of potentially-outdated gradient information, but differ in that the outdated gradients are all computed at the same previous iteration $\tilde{x}^k$ and the weighting of terms is changed. If the step-size and parameter $m$ are set appropriately, this algorithm has a linear convergence rate in the strongly-convex case. It can also achieve an $O(1/k)$ convergence  rate in the general convex case, by applying it to a purturbed problem where a strongly-convex regularizer has been added. 

A key advantage of this strategy is it only requires storing $\tilde{x}^k$, rather than the $n$ gradient values. However, to obtain this improved memory requirement we must evaluate $f_i'$ twice on each iteration. Further, the algorithm innefficiently requires full passes through the data to evaluate $f_j'(\tilde{x}^k)$ for all $j$.

A disadvantage of the SVRG method is that it requires setting two parameters (rather than one) and the convergence rate depends in a non-trivial way on their interaction both with each other and with both the Lipschitz constant and the strong-convexity. As with SDCA, the dependence on the strong-convexity constant is notable. In particular, the algorithm can't be applied without modification to general convex problems, and (although the effect is not as severe as it is with SDCA) the convergence rate may be substantially slower than SAG if there exists hidden strong convexity.

\textbf{Self-Concordant Objectives}: A surprising recent development is that \auth{Bach and Moulines}\citet{bach2013nonStrongly} have shown that the finite sum assumption is not required to obtain the $O(1/k)$ convergence rate in the special case of least squares. They show that an averaged SG method achieves this rate, and give a 2-phase Newton-like SG method that achieves this rate under an assumption similar to self-concordance. However, this assumption does not hold in general for the class of problems considered here.

\textbf{SAGA}: One of the most recent developments in this area is the SAGA algorithm~\cite{defazio2014saga}. This algorithm intelligently combines the updates used in the SAG and SVRG algorithms. This method maintains the appealing properties of SAG, but yields a simpler proof. However, the proof technique used in that work does not yield a simpler analysis of the original SAG algorithm.

\textbf{Lower Bounds}: There has been recent work on determining a lower bound on the convergence rate that can be expected for minimizing finite sums. \auth{Defazio et al.}\citet{defaziofinito} show that the rate must be at least $(1-1/n)$ in the strongly-convex case, while \auth{Agarwal and Bottou}\citet{argawal2014bounds} establish a bound that also depends on the condition number $L/\mu$.

\subsection{Generalizations and Other Issues}

{\bf Accelerated gradient}: AFG methods are variants of the basic FG method that can obtain an $O(1/k^2)$ convergence rate for the general convex case, and faster linear convergence rates that depend on the square root of the condition number $(\sqrt{L/\mu})$ rather than the condition number $(L/\mu)$ in the strongly-convex case~\citesee{\S2.2.1}{nesterov2004introductory}. For strongly-convex objectives, it has been shown that a mini-batch strategy can obtain a better dependence on the condition number in certain regimes for SDCA~\citep{shalev2013accelerated} and more recently for SVRG~\citep{nitanda2014stochastic}. It is possible that similar arguments could hold for SAG algorithms, which could be advantageous over these methods for reasons discussed in the previous section.

\auth{Shalev-Schwartz and Zhang}\citet{ASDCA} have also given an accelerated version of the SDCA method for ill-conditioned problems, that uses SDCA to solve a sequence of regularized problems up to a prescribed optimality. 
\iftoggle{springer}
{
\blu{Since this procedure ultimately relies on the sequence of primal solutions, we can also accelerate SAG using a procedure like this~\citep{lin2015universal}.}
}
{
Since this procedure ultimately relies on the sequence of primal solutions, we can also accelerate SAG using a procedure like this.
}
 However, a difficulty with this procedure (whether SAG or SDCA are used) is the cost of measuring the optimality of the sub-problems.

Rather than using this inner-outer procedure, \auth{Lin et al.}\citet{lin2014accelerated} show that using a deterministic primal iteration (based on the full data set) allows one to construct a primal solution that has the accelerated rate from the result of an accelerated dual coordinate ascent method. \auth{Zhang and Lin}\citet{zhang2014stochastic} give a coordinate-wise variant of the accelerated primal-dual method of \auth{Chambolle and Pock}\citet{chambolle2011first} that also achieves this rate. Based on these results, it is possible that accelerated versions of SAG could be developed that do not rely on an inner-outer procedure.

The convergence rates of these accelerated methods have the same form as the lower bound established by \auth{Agarwal and Bottou}~\citet{argawal2014bounds}. However, Agarwal and Bottou point out that these accelerated convergence rates depend on $\lambda$ rather than $\mu$. Thus, they conclude that the basic SAG algorithm may still be much faster on some problems than accelerated SDCA methods, and indeed show how to construct a problem where SAG is arbitrarily faster than accelerated SDCA. The possibility of developing an accelerated SAG algorithm that is adaptive to $\mu$ remains open.


{\bf Proximal gradient and ADMM}: It is becoming increasingly common to address problems of the form
\[
\minimize{x\in\Real^p}\quad r(x) + g(x) := r(x) + \frac{1}{n}\sum_{i=1}^n f_i(x),
\]
where $f_i$ and $g$ satisfy our assumptions and $r$ is a general convex function that could be non-smooth or enforce that constraints are satisfied.
Proximal-gradient methods for problems
with this structure use iterations of the form
\[
x^{k+1} = \textrm{prox}_{\alpha_k}\left[x^k - \frac{\alpha_k}{n}\sum_{i=1}^nf_i'(x^k)\right],
\]
where the prox$_{\alpha_k}[y]$ operator is the solution to the proximity problem
\[
\minimize{x\in\Real^p}\quad \frac{1}{2}\norm{x-y}^2 + \alpha_k r(x).
\]
Proximal-gradient and accelerated proximal-gradient methods are appealing for solving
non-smooth optimization problems because they achieve the same convergence
rates as FG methods for smooth optimization problems~\citep{nesterov2007gradient,SchmidtLeRouxBach11}.  We
have explored a variant of the proximal-gradient method where the average over the $f_i'(x^k)$ values is replaced by the 
SAG approximation~\eqref{eq:yi}.  Although our analysis does not directly apply to this scenario, we believe
that this proximal-SAG algorithm for composite non-smooth optimization
achieves the same convergence rates as the SAG algorithm for smooth optimization; this is supported by the experiments of \auth{Xiao and Zhang}\citet{xiao2014proximal}. Indeed, there now exist proximal-gradient variants of SDCA~\citep{PSDCA}, MISO~\citep{mairal2013surrogate}, SVRG~\citep{xiao2014proximal}, and SAGA~\citep{defazio2014saga}.

In cases where $r$ is composed with a linear function, we can consider approaches based on the alternating direction method of multipliers (ADMM). This has been explored for SDCA~\citep{suzuki2014stochastic} and MISO~\citep{zhong2013fastADMM}, and a variant based on SAG is also likely to be possible.

{\bf Coordinate-wise}:
The key advantage of SG and SAG methods is that the iteration cost is independent of the number of functions $n$.  
However, in many applications we may also be concerned with the dependence of the iteration cost on the number of variables~$p$.  Randomized coordinate-wise methods offer linear convergence rates with an iteration cost that is linear in $n$
but independent of~$p$~\citep{nesterov2010efficiency}.  We can consider a variant of SAG whose iteration cost
is independent of both $n$ and $p$ by using the update
\[
[y^k_{i}]_j = \begin{cases}
[f_i'(x^k)]_j & \textrm{if $i = i_k$ and $j = j_k$}\\
[y^{k-1}_i]_j & \textrm{otherwise,}
\end{cases}
\]
to each coordinate $j$ of each vector $y_i$ in the SAG algorithm, and where $j_k$ is a sample from the set $\{1,2,\dots,p\}$. \auth{Kon\v{c}n\'{y} et al.}\citet{konevcny2014semi} recently proposed an SVRG algorithm of this flavour.

{\bf Newton-like}:
In cases where $g$ is twice-differentiable, we can also consider Newton-like variants of the SAG algorithm,
\[
x^{k+1} = x^k - \frac{\alpha_k}{n}B^k[\sum_{i=1}^ny^k_i],
\]
where $B^k$ is a positive-definite approximation to the inverse Hessian $g''(x^k)$.  We would expect to obtain a faster
convergence rate by using an appropriate choice of the sequence $\{B^k\}$.  However, in order to not 
increase the iteration cost these matrices should be designed to allow fast multiplication.  For example,
we could choose $B^k$ to be diagonal, which would also preserve any sparsity present in the updates. \auth{Sohl-Dickstein et al.}\citet{sohl2014fast} propose a quasi-Newton method in this class that shows impressive empirical results, although the iteration cost is much higher.


{\bf Relaxing Convexity Assumptions}:
It is likely that the convexity assumptions made in this work could be relaxed. For example, \auth{Gong and Ye}\citet{gong2014linear} show SVRG can obtain a linear convergence under weaker assumptions. Since SAG is adaptive to hidden strong-convexity, the assumptions needed for its linear convergence are likely even weaker. 
Further, the SAG algorithm may also be useful even if $g$ is non-convex (which is different than SDCA). In this case we expect that, similar to the IAG method~\citep{blatt2008convergent}, the algorithm converges to a stationary point under very general conditions.

{\bf Non-Uniform Sampling}:
We have given an argument that non-uniform sampling should benefit the SAG algorithm, and shown empirically that it can lead to a substantial improvement. However, we have not yet given a full analysis of this scheme. Subsequent works have shown that the type of dependency we conjecture here (e.g., dependence on the average Lispschitz constant) can be achieved with non-uniform sampling in the context of SDCA~\citep{qu2014randomized,zhao2014stochastic}, SVRG~\citep{xiao2014proximal}, and SAGA~\citep{schmidt2015sag4crf}

{\bf Step-size selection and termination criteria}:
}
The three major disadvantages of SG methods are: (i) the slow convergence rate, (ii) deciding when to terminate the algorithms, and (iii) choosing the step size while running the algorithm.  This work shows that the SAG iterations achieve a much faster convergence rate,
but the SAG iterations may also be advantageous in terms of termination criteria and choosing step sizes.  In particular, the SAG iterations suggest a natural termination criterion; since the quantity $d$ in Algorithm~\ref{alg:SAG} converges to $f'(x^k)$ as $\norm{x^k-x^{k-1}}$ converges to zero, we can use $\norm{d}$ as an approximation of the optimality of $x^k$ as the iterates converge. Regarding choosing the step-size, a disadvantage of a constant step-size strategy is that a step-size that is too large may cause divergence.  But, we expect that it is possible to design line-search or trust-region strategies that avoid this issue. Such strategies might even lead to faster convergence for functions that are locally well-behaved around their optimum, as indicated in our experiments. Further, while SG methods require specifying a sequence of step sizes and mis-specifying this sequence can have a disastrous effect on the convergence rate~\citesee{\S2.1}{nemirovski2009robust}, our theory shows that the SAG iterations achieve a fast convergence rate for any sufficiently small constant
step size, and our experiments indicate that a simple line-search gives strong performance.

\iftoggle{springer}
{
\begin{acknowledgements}
}
{
\subsection*{Acknowledgements}
}
We would like to thank the anonymous reviewers for their many useful comments.
This work was partially supported by the European Research Council (SIERRA-ERC-239993) and a Google Research Award.
Mark Schmidt is also supported by a postdoctoral fellowship from the Natural Sciences and Engineering Research Council of Canada.
\iftoggle{springer}
{
\end{acknowledgements}
}

\iftoggle{springer}
{

}
{
\appendix

\section*{Appendix A: Comparison of convergence rates}
\addtocounter{section}{1}

We now present a comparison of the convergence rates of primal and dual FG and coordinate-wise methods to the rate of the SAG method in terms of effective passes through the data for $\ell_2$-regularized least-squares regression. In particular, in this appendix we will consider the problem
\[
\minimize{x\in\Real^p} \quad g(x) \defd \frac{\lambda}{2}\norm{x}^2 + \frac{1}{2n}\sum_{i=1}^n(a_i^Tx - b_i)^2,
\]
where to apply the SG or SAG method we can use
\[
f_i(x) := \frac{\lambda}{2}\norm{x}^2 + \frac{1}{2}(a_i^Tx - b_i)^2.
\]
If we use $b$ to denote a vector containing the values $b_i$ and $A$ to denote a matrix withs rows $a_i$, we can re-write this problem as
\[
\minimize{x\in\Real^p} \quad \frac{\lambda}{2}\norm{x}^2 + \frac{1}{2n}\norm{Ax - b}^2.
\]
The Fenchel dual of this problem is
\[
\minimize{y\in\Real^n} \quad d(y) \defd \frac{n}{2}\norm{y}^2 + \frac{1}{2\lambda}y^\top AA^\top y + y^\top b.
\]
We can obtain the primal variables from the dual variables by the formula $x = (-1/\lambda)A^\top y$.
Convergence rates of different primal and dual algorithms are often expressed in terms of the following Lipschitz constants:
\begin{align*}
&L_g = \lambda + M_\sigma/n & & \textrm{(Lipschitz constant of $g'$)}\\
&L_g^i = \lambda + M_i \quad & & \textrm{(Lipschitz constant for all $f_i'$)}\\
&L_g^j = \lambda + M_j/n \quad & & \textrm{(Lipschitz constant of all $g_j'$)}\\
&L_d = n + M_\sigma/\lambda & & \textrm{(Lipschitz constant of $d'$)}\\
&L_d^i = n + M_i/\lambda \quad & & \textrm{(Lipschitz constant of all $d_i'$)}\\
\end{align*}
Here, we use $M_\sigma$ to denote the maximum eigenvalue of $A^\top A$, $M_i$ to denote the maximum squared row-norm $\max_i\{\norm{a_i}^2\}$, and $M_j$ to denote the maximum squared column-norm $\max_j\{\sum_{i=1}^n(a_i)^2_j\}$. We use $g_j'$ to refer to element of $j$ of $g'$, and similarly for $d_i'$. The convergence rates will also depend on the primal and dual strong-convexity constants:
\begin{align*}
&\mu_g = \lambda + m_\sigma/n & & \textrm{(Strong-convexity constant of $g$)}\\
&\mu_d = n + m_\sigma'/\lambda & & \textrm{(Strong-convexity constant of $d$)}
\end{align*}
Here, $m_\sigma$ is the minimum eigenvalue of $A^\top A$, and $m_\sigma'$ is the minimum eigenvalue of $AA^\top$.

\subsection{Full Gradient Methods}

Using a similar argument to~\cite[Theorem 2.1.15]{nesterov2004introductory}, if we use the basic FG method with a step size of $1/L_g$, then $(f(x^k) - f(x^\ast))$ converges to zero with rate
\[
\left(1-\frac{\mu_g}{L_g}\right)^2 
= \left(1 - \frac{\lambda + m_\sigma/n}{\lambda + M_\sigma/n}\right)^2 
= \left(1 - \frac{n\lambda + m_\sigma}{n\lambda + M_\sigma}\right)^2 
\leqslant \exp\left(-2\frac{n\lambda + m_\sigma}{n\lambda + M_\sigma}\right),
\]
while a larger step-size of $2/(L_g+\mu_g)$ gives a faster rate of
\[
\left(1 - \frac{\mu_g + \mu_g}{L_g + \mu_g}\right)^2
= \left(1 - \frac{n\lambda + m_\sigma}{n\lambda + (M_\sigma + m_\sigma)/2}\right)^2
\leqslant \exp\left(-2\frac{n\lambda + m_\sigma}{n\lambda + (M_\sigma + m_\sigma)/2}\right),
\]
and we see that the speed improvement is determined by how much smaller $m_\sigma$ is than $M_\sigma$.

If we use the basic FG method on the dual problem with a step size of $1/L_d$, then $(d(x^k) - d(x^\ast))$ converges to zero with rate
\[
\left(1-\frac{\mu_d}{L_d}\right)^2
= \left(1-\frac{n + m_\sigma'/\lambda}{n + M_\sigma/\lambda}\right)^2 
= \left(1-\frac{n\lambda + m_\sigma'}{n\lambda + M_\sigma}\right)^2
\leqslant \exp\left(-2\frac{n\lambda + m'_\sigma}{n\lambda + M_\sigma}\right),
\]
and with a step-size of $2/(L_d+\mu_d)$ the rate is
\[
\left(1 - \frac{\mu_d + \mu_d}{L_d + \mu_d}\right)^2
= \left(1 - \frac{ n\lambda + m_\sigma'}{n\lambda + (M_\sigma + m_\sigma')/2}\right)^2
\leqslant \exp\left(-2\frac{n\lambda + m'_\sigma}{n\lambda + (M_\sigma + m'_\sigma)/2}\right).
\]
Thus, whether we can solve the primal or dual method faster depends on $m_\sigma$ and $m_\sigma'$. In the over-determined case where $A$ has independent columns, a primal method should be preferred. In the under-determined case where $A$ has independent rows, we can solve the dual more efficiently. However, we note that a convergence rate on the dual objective does not necessarily yield the same rate in the primal objective. If $A$ is invertible (so that $m_\sigma=m_\sigma'$) or it has neither independent columns nor independent rows (so that $m_\sigma = m_\sigma' = 0$), then there is no difference between the primal and dual rates.

The AFG method achieves a faster rate. Applied to the primal with a step-size of $1/L_g$ it has a rate of~\citethm{Theorem 2.2.2}{nesterov2004introductory}
\[
\left(1-\sqrt{\frac{\mu_g}{L_g}}\right)
= \left(1-\sqrt{\frac{\lambda + m_\sigma/n}{\lambda + M_\sigma/n}}\right)
= \left(1-\sqrt{\frac{n\lambda + m_\sigma}{n\lambda + M_\sigma}}\right)
\leqslant \exp\left(-\sqrt{\frac{n\lambda + m_\sigma}{n\lambda + M_\sigma}}\right),
\]
and applied to the dual with a step-size of $1/L_d$ it has a rate of
\[
\left(1-\sqrt{\frac{\mu_d}{L_d}}\right)
= \left(1-\sqrt{\frac{n + m_\sigma'\lambda}{n + M_\sigma/\lambda}}\right)
= \left(1-\sqrt{\frac{n\lambda + m_\sigma'}{n\lambda + M_\sigma}}\right)
\leqslant \exp\left(-\sqrt{\frac{n\lambda + m'_\sigma}{n\lambda + M_\sigma}}\right).
\]
\subsection{Coordinate-Descent Methods}

The cost of applying one iteration of an FG method is $O(np)$. For this same cost we could apply $p$ iterations of a coordinate descent method to the primal, assuming that selecting the coordinate to update has a cost of $O(1)$. If we select coordinates uniformly at random, then the convergence rate for $p$ iterations of coordinate descent with a step-size of $1/L_g^j$ is~\citethm{Theorem 2}{nesterov2010efficiency}
\[
\left(1-\frac{\mu_g}{pL_g^j}\right)^p
= \left(1-\frac{\lambda + m_\sigma/n}{p(\lambda + M_j/n)}\right)^p
= \left(1-\frac{n\lambda + m_\sigma}{p(n\lambda + M_j)}\right)^p
\leqslant \exp\left(-\frac{n\lambda + m_\sigma}{n\lambda + M_j}\right).
\]
Here, we see that applying a coordinate-descent method can be much more efficient than an FG method if $M_j << M_\sigma$. This can happen, for example, when the number of variables $p$ is much larger than the number of examples $n$. Further, it is possible for coordinate descent to be faster than the AFG method if the difference between $M_\sigma$ and $M_j$ is sufficiently large.

For the $O(np)$ cost of one iteration of the FG method, we could alternately perform $n$ iterations of coordinate descent on the dual problem. With a step size of $1/L_d^i$ this would obtain a rate on the dual objective of
\[
\left(1-\frac{\mu_d}{nL_d^i}\right)^n
= \left(1-\frac{n + m_\sigma'/\lambda}{n(n + M_i/\lambda)}\right)^n
= \left(1-\frac{n\lambda + m_\sigma'}{n(n\lambda + M_i)}\right)^n
\leqslant \exp\left(-\frac{n\lambda + m'_\sigma}{n\lambda + M_i}\right),
\]
which will be faster than the dual FG method if $M_i << M_\sigma$. This can happen, for example, when the number of examples $n$ is much larger than the number of variables $p$. The difference between the primal and dual coordinate methods depends on $M_i$ compared to $M_j$ and $m_\sigma$ compared to $m_\sigma'$. 

The analysis of SDCA gives a convergence rate in terms of the duality gap (and hence the primal) rather than simply in terms of the dual~\citet{schwartz12}. Using the SDCA analysis, we obtain a rate in the duality gap of
\[
\left(1-\frac{1}{n + M_i/\lambda}\right)^n
= \left(1-\frac{\lambda}{n\lambda + M_i}\right)^n
\leqslant \exp\left(-\frac{n\lambda}{n\lambda + M_i}\right),
\]
which is the same as the dual rate given above when $m_\sigma'=0$, and is slower otherwise.

\subsection{Stochastic Average Gradient}

For the $O(np)$ cost of one iteration of the FG method, we can perform $n$ iterations of SAG. With a step size of $1/16L$, performing $n$ iterations of the SAG algorithm has a rate of
\iftoggle{springer}
{
\begin{align*}
\left(1 - \min\left\{\frac{\mu_g}{16L_g^i},\frac{1}{8n}\right\}\right)^n
& = \left(1 - \min\left\{\frac{\lambda + m_\sigma/n}{16(\lambda + M_i)},\frac{1}{8n}\right\}\right)^n
\\
& \leqslant \exp\left(-\frac{1}{16}\min\left\{  \frac{n\lambda + m_\sigma}{\lambda + M_i},2 \right\}\right)
\end{align*}
}
{
\[
\left(1 - \min\left\{\frac{\mu_g}{16L_g^i},\frac{1}{8n}\right\}\right)^n
= \left(1 - \min\left\{\frac{\lambda + m_\sigma/n}{16(\lambda + M_i)},\frac{1}{8n}\right\}\right)^n
\leqslant \exp\left(-\frac{1}{16}\min\left\{  \frac{n\lambda + m_\sigma}{\lambda + M_i},2 \right\}\right)
\]
}
In the case where $n \leqslant 2L_g^i/\mu_g$, this is most similar to the rate obtained with the dual coordinate descent method, though there is a constant factor of $16$ and the rate depends on the primal strong convexity constant $m_\sigma$ rather than the dual $m_\sigma'$. However, in this case the SAG rate will often be faster because the term $n\lambda$ in the denominator is replaced by $\lambda$. An interesting aspect of the SAG rate is that unlike other methods the convergence rate of SAG reaches a limit: the convergence rate improves as $n$ grows and as the condition number  $L_g^i/\mu_g$ decreases but no further improvement in the convergence rate is obtained beyond the point where $n = 2L_g^i/\mu_g$. Thus, while SAG may not be the method of choice for very well-conditioned problems, it is a robust choice as it will always tend to be among the best methods in most situations.



\section*{Appendix B: Proof of the theorem}
\addtocounter{section}{1}
\addtocounter{subsection}{-3}

In this Appendix, we give the proof of Theorem~\ref{thm}. 

\subsection{Problem set-up and notation}

Recall that we use $g = \frac{1}{n}\sum_{i=1}^n f_i$ to denote a $\mu$-strongly convex function, where the functions $f_i$ for $i=1,\ldots,n$ are convex functions from $\rb^p$ to $\rb$ with $L$-Lipschitz continuous gradients. In this appendix we will use the convention that $\mu \geqslant 0$ so that the regular convex case (where $\mu = 0$) is allowed. We assume that a minimizer of $g$ is attained by some parameter $x^\ast$ (such a value always exists and is unique when $\mu>0$).

Recall that the SAG algorithm performs the recursion (for $k \geqslant 1$):
\BEAS
x^k & = &  x^{k-1} - \frac{\alpha}{n} \sum_{i=1}^n y_i^{k},
\EEAS
where an integer $i_k$ is selected uniformly at random  from $\{1,\dots,n\}$ and we set
\[
y^k_i = \begin{cases}
f_i'(x^{k-1}) & \textrm{if $i = i_k$,}\\
y^{k-1}_i & \textrm{otherwise.}
\end{cases}
\]
We will use the notation
\begin{align*}
y^k = \left(
\begin{array}{c}
y_1^k\\
\vdots\\
y_n^k
\end{array}\right) \in \rb^{n p},
\qquad
\theta^k =
\left(
\begin{array}{c}
y_1^k\\
\vdots\\
y_n^k\\
x^k
\end{array}\right) \in \rb^{(n+1)p},
\qquad
\theta^\ast =
\left(
\begin{array}{c}
f_1'(x^\ast)\\
\vdots\\
f_n'(x^\ast)\\
x^\ast
\end{array}\right) \in \rb^{(n+1)p} \; ,
\end{align*}
and we will also find it convenient to use
\begin{align}
\label{eq:e}
e = \left(
\begin{array}{c}
\idm\\
\vdots\\
\idm
\end{array}\right) \in \rb^{n p \times p},
\qquad
f'(x) = \left(
\begin{array}{c}
f_1'(x)\\
\vdots\\
f_n'(x)
\end{array}\right) \in \rb^{n p}.
\end{align}
With this notation, note that $g'(x) = \frac{1}{n}e^\top f'(x)$ and $x^k = x^{k-1} - \frac{\alpha}{n}e^\top y^k$.

For a square $n \times n$ matrix $M$, we use $\diag(M)$ to denote a vector of size $n$ composed of the diagonal of $M$, while for a vector $m$ of dimension $n$,
$\Diag(m)$ is the $n \times n$ diagonal matrix with $m$ on its diagonal. Thus $\Diag(\diag(M))$ is a diagonal matrix with the diagonal elements of $M$ on its diagonal, and $\diag(\Diag(m))=m$.

In addition, if $M$ is a $tp \times tp$ matrix and $m$ is a $tp \times p$ matrix, then we will use the convention that:
\begin{itemize}
\item $\diag(M)$ is the $tp \times p$ matrix being the concatenation of the $t$ ($p \times p$)-blocks on the diagonal of $M$;
\item $\Diag(m)$ is the $tp \times tp$ block-diagonal matrix whose ($p \times p$)-blocks on the diagonal are equal to the ($p \times p$)-blocks of $m$.
\end{itemize}
Finally,~$\mathcal{F}_{k}$ will denote the $\sigma$-field of information up to (and including) time $k$. In other words, $\mathcal{F}_k$ is the $\sigma$-field generated by $i_1,\dots,i_k$. Given $\mathcal{F}_{k-1}$, we can write the expected values of the $y^k$ and $x^k$ variables in the SAG algorithm as
\BEAS
\E (y ^k | \mathcal{F}_{k-1} ) &  = &\big(1 - \frac{1}{n} \big) y^{k-1} + \frac{1}{n} f'(x^{k-1}), \\
\E ( x^k | \mathcal{F}_{k-1} ) & = & x^{k-1} - \frac{\alpha}{n} \big(1 - \frac{1}{n} \big) e^\top y^{k-1}  - \frac{\alpha}{n^2} e^\top f'(x^{k-1}).
\EEAS

The proof is based on finding a Lyapunov function $\mathcal{L}$ from $\rb^{(n+1)p}$ to $\rb$ such that the sequence $\E \mathcal{L}(\theta^k)$ decreases at an appropriate rate, and $\E \mathcal{L}(\theta^k)$ dominates $[g(x^k) - g(x^\ast)]$. We derive these results in a parameterized way, leaving a variety of coefficients undetermined but tracking the constraints required of the coefficients as we go. Subsequently, we guide the setting of these coefficients by using a second-order cone solver, and verify the validity of the resulting coefficients using a symbolic solver to check positivity of certain polynomials. Finally, the constants in the convergence rate are given by the initial values of the Lyapunov function, based on the choice of $y^0$. 

The Lyapunov function contains a term of the form~\mbox{$(\theta^k - \theta^*)^\top \left(\begin{array}{cc} A & B \\ B^\top & C \end{array}\right) (\theta^k - \theta^*)$} for some values of $A$, $B$ and $C$. Our analysis makes use of the following lemma, derived in~\citethm{Appendix~A.4}{roux2012stochastic}, showing how this quadratic form evolves through the SAG recursion in terms of the elements of $\theta^{k-1}$.
  \begin{lemma}
  \label{lemma}
If $ P = \left(
\begin{array}{cc} A & B \\ B^\top &  C
\end{array}
\right)$,
for $A \in \rb^{np \times np}$, $B \in \rb^{np \times p}$ and $C \in \rb^{p \times p}$, then
\begin{align*}
&\E\left[\left.(\theta^k - \theta^*)^\top \left(\begin{array}{cc} A & B \\ B^\top &  C
\end{array}\right) (\theta^k - \theta^*)\right| \mathcal{F}_{k-1}\right] \nonumber\\
&\hspace*{2.5cm}=(y^{k-1} - f'(x^\ast))^\top \left[\left(1 - \frac{2}{n}\right)S + \frac{1}{n}\Diag(\diag(S))\right](y^{k-1} - f'(x^\ast))\\
&\hspace*{2.5cm}+ \frac{1}{n}(f'(x^{k-1}) - f'(x^\ast))^\top \Diag(\diag(S))(f'(x^{k-1}) - f'(x^\ast))\\
&\hspace*{2.5cm}+ \frac{2}{n}(y^{k-1} - f'(x^\ast))^\top \left[S - \Diag(\diag(S))\right]  (f'(x^{k-1}) - f'(x^\ast))\\
&\hspace*{2.5cm}+ 2\left(1 - \frac{1}{n}\right) (y^{k-1} - f'(x^\ast))^\top \left[B - \frac{\alpha}{n}eC\right](x^{k-1} - x^\ast)\\
&\hspace*{2.5cm}+ \frac{2}{n}(f'(x^{k-1}) - f'(x^\ast))^\top\left[B - \frac{\alpha}{n}eC\right](x^{k-1} - x^\ast)\\
&\hspace*{2.5cm}+ (x^{k-1} - x^\ast)^\top C (x^{k-1} - x^\ast) \; ,
\end{align*}
with
\[
S = A - \frac{\alpha}{n} B e^\top - \frac{\alpha}{n} e B^\top + \frac{\alpha^2}{n^2} e Ce^\top  \in \rb^{np \times np} \; .
\]
\end{lemma}
Our proof also uses the following lemma, giving the inverse of a highly-structured matrix that arises in the analysis.
\begin{lemma}
\label{lemma2}
Let $I$ be the identity in $\rb^{np\times np}$ and $e \in \rb^{np \times p}$ be defined by stacking identity matrices in $\rb^{n \times n}$ as in~\eqref{eq:e}. If $\alpha$ and $\beta$ are non-zero real-valued scalars, then it holds that
\[
\left(\alpha\left(I - \frac{1}{n}ee^T\right) + \beta\left(\frac{1}{n}ee^T\right)\right)^{-1} = \frac{1}{\alpha}\left(I - \frac{1}{n}ee^T\right) + \frac{1}{\beta}\left(\frac{1}{n}ee^T\right).
\]
\end{lemma}
\begin{proof}
It is sufficient to verify that the inverse of the left side times the right side is equal to the identity matrix,\begin{align*}
& \left(\alpha\left(I - \frac{1}{n}ee^T\right) + \beta\left(\frac{1}{n}ee^T\right)\right)\left(\frac{1}{\alpha}\left(I - \frac{1}{n}ee^T\right) + \frac{1}{\beta}\left(\frac{1}{n}ee^T\right)\right)\\
& = \left(I - \frac{1}{n}ee^T\right)\left(I - \frac{1}{n}ee^T\right) + \frac{\alpha}{\beta}\left(I - \frac{1}{n}ee^T\right)\left(\frac{1}{n}ee^T\right) \\&+ \frac{\beta}{\alpha}\left(\frac{1}{n}ee^T\right)\left(I - \frac{1}{n}ee^T\right) + \left(\frac{1}{n}ee^T\right)\left(\frac{1}{n}ee^T\right)\\
& = \left(I - \frac{2}{n}ee^T + \frac{1}{n}ee^T\right) + \frac{\alpha}{\beta}\left(\frac{1}{n}ee^T - \frac{1}{n}ee^T\right)
+\frac{\beta}{\alpha}\left(\frac{1}{n}ee^T - \frac{1}{n}ee^T\right) + \frac{1}{n}ee^T\\
& = I.
\end{align*}
where we use that $e^Te = nI$, giving $\left(\frac{1}{n}ee^T\right)\left(\frac{1}{n}ee^T\right) = \frac{1}{n}ee^T$.
\end{proof}

\subsection{General Lyapunov function}

For some $h\geqslant 0$, we consider a Lyapunov function of the form
 $$
 \mathcal{L}(\theta^k) = 
2 h g( x^k + d e^\top y^k) - 2 h  g(x^\ast)
+ (\theta^k - \theta^*)^\top \left(\begin{array}{cc} A & B \\ B^\top &  C
\end{array}\right) (\theta^k - \theta^*)
 $$
 with 
\[
A = a_1ee^\top + a_2\idm, \quad B = be, \quad C = c\idm.
\]
 To achieve the desired convergence rate, our goal is to show for appropriate values of $\delta \geqslant 0 $ and $\gamma \geqslant 0$ that
 \BEAS
 (a) & & \E \big( \mathcal{L}(\theta^k) | \mathcal{F}_{k-1}  \big) \leqslant ( 1 - \delta ) \mathcal{L}(\theta^{k-1}),
 \\
(b) & &  \mathcal{L}(\theta^k) \geqslant \gamma \big[ g(x^{k}) - g(x^\ast) \big],
 \EEAS
 almost surely. Thus, in addition to the algorithm parameter $\alpha$, there are 2 parameters of the result \{$\gamma$, $\delta$\} and 6 parameters of the Lyapunov function \{$a_1$, $a_2$, $b$, $c$, $d$, $h$\}.

 \subsection{Lyapunov upper bound}

To show (a), we derive an upper bound on the quantity
\begin{align*}
& \E \big( \mathcal{L}(\theta^k) | \mathcal{F}_{k-1}  \big) - ( 1 - \delta )  \mathcal{L}(\theta^{k-1})\\
& = 
2 h \E \big[ g( x^k + d e^\top y^k)\big] - 2 h  g(x^\ast)
- (1-\delta)(
2 h g( x^{k-1} + d e^\top y^{k-1}) - 2 h  g(x^\ast))\\
& \smaller{+ \E \big[(\theta^k - \theta^*)^\top \left(\begin{array}{cc} A & B \\ B^\top &  C
\end{array}\right) (\theta^k - \theta^*)\big]
- (1-\delta) (\theta^{k-1} - \theta^*)^\top \left(\begin{array}{cc} A & B \\ B^\top &  C
\end{array}\right) (\theta^{k-1} - \theta^*).}
\end{align*}
   Lemma~\ref{lemma} gives an expression for the expectation over the quadratic term in the last line, and we will have
\BEAS
 S & = & a_2 \idm + \bigg( a_1 - 2\frac{\alpha}{n}b + \frac{\alpha^2}{n^2} c  \bigg) ee^\top \\
\Diag(\diag(S)) & = &  \bigg( a_2 + a_1 - 2\frac{\alpha}{n}b + \frac{\alpha^2}{n^2} c  \bigg) \idm 
\\
S - \Diag(\diag(S)) & = &  \bigg( a_1 - 2\frac{\alpha}{n}b + \frac{\alpha^2}{n^2} c  \bigg) ( ee^\top - \idm).
 \EEAS
 We also use that strong convexity of $g$ implies
   \BEA
   \label{eq:11}
    2 g( x^{k-1}+ d e^\top y^{k-1}) 
    & \geqslant & 2 g(x^{k-1}) + 2 d g'(x^{k-1})^\top e^\top y^{k-1} + \mu d^2 \| e^\top y^{k-1} \|^2.
   \EEA
 Further, using the Lipschitz-continuity of $g'$ and the identity $ x^k + d e^\top y^k = x^{k-1} + ( d - \frac{\alpha}{n} ) e^\top y^k$, followed by applying Lemma~\ref{lemma} using $ A = ee^\top$ to expand the last term, gives us the bound
  \begin{equation}
\label{eq:fBound}
  \begin{aligned}
  &  \E \big[ 2 g( x^k + d e^\top y^k)  | \mathcal{F}_{k-1} \big] \\
  & \leqslant     2 g( x^{k-1}) + 2 \big(d - \frac{\alpha}{n} \big) g'(x^{k-1})^\top 
  \E \big[ e^\top 
  y^{k}  | \mathcal{F}_{k-1} \big]
  + L \big(d - \frac{\alpha}{n} \big)^2 \E \big[ \| e^\top y^{k} \|^2  | \mathcal{F}_{k-1} \big] \\
 & =     2 g( x^{k-1}) + 2 \big(d - \frac{\alpha}{n} \big) g'(x^{k-1})^\top 
  \big[  (1 - \frac{1}{n}) e^\top y^{k-1} + \frac{1}{n} e^\top f'(x^{k-1}) \big]
   \\
  & +  L \big(d - \frac{\alpha}{n} \big)^2 
  \bigg[
   (y^{k-1} - f'(x^\ast))^\top \left[\left(1 - \frac{2}{n}\right)ee^\top + \frac{1}{n}\idm \right](y^{k-1} - f'(x^\ast))  \bigg]
   \\
   &  +   L \big(d - \frac{\alpha}{n} \big)^2 
  \bigg[
   \frac{1}{n} \| f'(x^{k-1}) - f'(x^\ast)\|^2
   \\
 & \hspace{3cm}+ \frac{2}{n}(y^{k-1} - f'(x^\ast))^\top \left[ee^\top - \idm\right]  (f'(x^{k-1}) - f'(x^\ast)) \bigg].
   \end{aligned}
   \end{equation}
 Combining Lemma~\ref{lemma} with Inequalities~\eqref{eq:11} and~\eqref{eq:fBound} yields an upper bound on $\E \big( \mathcal{L}(\theta^k) | \mathcal{F}_{k-1}  \big) - ( 1 - \delta )  \mathcal{L}(\theta^{k-1})$ with a large number of terms. To help in the process of simplifying these terms, Table~\ref{table:Lyap} lists all the terms up to a scalar factor $s$ and possibly a matrix $M$. The table also gives these scalars and matrices for each term, as well as the source of the term.

\begin{sidewaystable}
\label{table:Lyap}
{
\centering
\vspace{300pt}
\begin{tabular}{lllll}
Term  & Scalar $s$ & Matrix $M$ & Source & Group\\
\hline
  $sg(x^\ast)$ & $-2h$ & & $\E \big( \mathcal{L}(\theta^k) | \mathcal{F}_{k-1}  \big)$ & 0\\
         $sg(x^\ast)$ & $2h(1-\delta)$  & &  $L(\theta^{k-1})$& 0\\
        $s(y^{k-1} - f'(x^\ast))^\top M (y^{k-1} - f'(x^\ast))$ & $-(1-\delta)$ & $a_1ee^\top + a_2I$ & $L(\theta^{k-1})$ & 3 and 4\\
        $s(y^{k-1} - f'(x^\ast))^\top M (x^{k-1} -x^\ast)$ & $-2(1-\delta)$ & $be$ & $L(\theta^{k-1})$ & 5\\
        $s(x^{k-1} - x^\ast)\top M (x^{k-1} - x^\ast)$ & $-(1-\delta)$ & $cI$ & $L(\theta^{k-1})$ & 9\\
\hline
$sg(x^{k-1})$ & $-2h(1-\delta)$ & & Inequality~\eqref{eq:11} & 0\\
$sg'(x^{k-1})^\top M y^{k-1}$ & $-2h(1-\delta)d$ & $e^\top$ & Inequality~\eqref{eq:11} & 7\\
$s(y^{k-1})^\top M y^{k-1}$ & $-h(1-\delta)\mu d^2$ & $ee^\top$ & Inequality~\eqref{eq:11} & 4\\
\hline
  $sg(x^{k-1})$ & $2h$ & & Inequality~\eqref{eq:fBound} & 0\\
                  $sg'(x^{k-1})^\top M y^{k-1}$ & $2h(d-\frac{\alpha}{n})$ & $ (1 - \frac{1}{n})e^\top$ & Inequality~\eqref{eq:fBound} & 7\\
                  $sg'(x^{k-1})^\top M f'(x^{k-1})$ & $2h(d-\frac{\alpha}{n})$ & $\frac{1}{n}e^\top$ & Inequality~\eqref{eq:fBound} & 2\\
                  $s(y^{k-1} - f'(x^\ast))^\top M (y^{k-1} - f'(x^\ast))$ & $Lh(d-\frac{\alpha}{n})^2$ & $\left(1-\frac{2}{n}\right)ee^\top + \frac{1}{n}I$ & Inequality~\eqref{eq:fBound} & 3 and 4\\
                  $s(f'(x^{k-1}) - f'(x^\ast))^\top M (f'(x^{k-1}) - f'(x^\ast))$ & $\frac{Lh}{n}(d-\frac{\alpha}{n})^2$ & $I$ & Inequality~\eqref{eq:fBound} & 8\\
                  $s(y^{k-1} - f'(x^\ast))^\top M (f'(x^{k-1}) - f'(x^\ast))$ & $\frac{2Lh}{n}(d-\frac{\alpha}{n})^2$ & $ee^\top - I$ & Inequality~\eqref{eq:fBound} & 6 and 7\\
                  \hline
$s(y^{k-1} - f'(x^\ast))^\top M (y^{k-1} - f'(x^\ast))$ & $1$ & $\left(1-\frac{2}{n}\right)S + \frac{1}{n}\Diag(\diag(S))$ & Lemma~\ref{lemma} & 3 and 4\\
$s(y^{k-1} - f'(x^\ast))^\top M (x^{k-1} -x^\ast)$ & $2\left(1-\frac{1}{n}\right)$ & $(b - \frac{\alpha}{n}c)e$ & Lemma~\ref{lemma} & 5\\
$s(x^{k-1} - x^\ast)\top M (x^{k-1} - x^\ast)$ & $1$ & $cI$ & Lemma~\ref{lemma} & 9\\
$s(f'(x^{k-1}) - f'(x^\ast))^\top M (f'(x^{k-1}) - f'(x^\ast))$ & $\frac{1}{n}$ & $\Diag(\diag(S))$ & Lemma~\ref{lemma} & 8\\
$s(y^{k-1} - f'(x^\ast))^\top M (f'(x^{k-1}) - f'(x^\ast))$ & $\frac{2}{n}$ & $S - \Diag(\diag(S))$ & Lemma~\ref{lemma} & 6 and 7\\
$s(f'(x^{k-1}) - f'(x^\ast))^\top M (x^{k-1} - x^*)$ & $\frac{2}{n}$ & $(b - \frac{\alpha}{n}c)e$ & Lemma~\ref{lemma} & 1
\end{tabular}
\caption{Expressions in upper bound on  $\E \big( \mathcal{L}(\theta^k) | \mathcal{F}_{k-1}  \big) - ( 1 - \delta )  \mathcal{L}(\theta^{k-1})$.}
}
\iftoggle{springer}
{
}
{
\vspace{300pt}
}
\end{sidewaystable}


We combine the expressions from Table~\ref{table:Lyap} using the stated groups in the last column. For example, we add together all the terms in group $0$ to obtain a scalar coefficient $B_0$ for terms in this group. We do this in the straightforward way (using that $g'(x) = \frac{1}{n}e^\top f'(x)$ and $g'(x^\ast)=0$) for groups $0$, $1$, $2$, $5$, $8$, and $9$. For groups $6$ and $7$ we split into terms that have an identity matrix $M$ ($B_6$) and terms with an $ee^\top$ term ($B_7$). Similarly, $B_3$ comes from terms with an identity matrix in $M$ and $B_4$ correspond to terms with an $ee^\top$ term. Combining expressions in this way (and adding/subtracting $(B_3/n)ee^\top$) gives: 
 \BEAS
  & & \E \big( \mathcal{L}(\theta^k) | \mathcal{F}_{k-1}  \big) -  ( 1 - \delta ) \mathcal{L}(\theta^{k-1})
  \\
  & \leqslant & B_0 \big( g(x^{k-1}) -  g(x^\ast)   \big) + B_9 \| x^{k-1} - x^\ast \|^2  
  + (x^{k-1} - x^\ast)^\top g'(x^{k-1})
  B_1 - B_2 \| g'(x^{k-1}) \|^2
\\
    & &   -  (  y^{k-1} - y^\ast)^\top \bigg[
    B_3 ( \idm - \frac{1}{n} ee^\top ) + B_4\frac{1}{n}  ee^\top 
    \bigg]  ( y^{k-1} - y^\ast)  \\
  & & + ( y^{k-1} - y^\ast)^\top \bigg[
  B_5  e  ( x^{k-1} - x^\ast) + B_6 (  f'(x^{k-1}) - f'(x^\ast) ) + B_7 e g'(x^{k-1} )
  \bigg]  \\
  & & +B_8 \| f'(x^{k-1}) - f'(x^\ast)\|^2,
  \EEAS 
 with 
\BEAS
B_0 & = & 2 \delta h\\
B_1 & = & 2 ( b- \frac{\alpha}{n} c ) \\
B_2 & = & 2 ( \frac{\alpha}{n} - d ) h \\
B_3 & = & - \bigg[
  \left(1 - \frac{2}{n}\right)a_2 + \frac{1}{n}\big[
   a_1 + a_2  - 2\frac{\alpha}{n} b + \frac{\alpha^2}{n^2} c 
   \big] - ( 1- \delta) a_2 +  L h \frac{1}{n} \big(d - \frac{\alpha}{n} \big)^2  
 \bigg] \\
B_4 & = & B_3 \\
& - & n \bigg[
   \left(1 - \frac{2}{n}\right)  ( a_1 - 2\frac{\alpha}{n} b + \frac{\alpha^2}{n^2} c )
  - ( 1- \delta) a_1 
   +  L ( 1  - \frac{2}{n} ) h  \big(d - \frac{\alpha}{n} \big)^2  
   - ( 1- \delta) \mu h d^2
   \bigg]\\
B_5 & = &  2\bigg[
  ( \delta - \frac{1}{n} ) b  - \frac{\alpha}{n} (1 - \frac{1}{n} ) c
  \bigg] \\
B_6 & = &  -  \frac{2}{n}    \bigg(
 h  L \big(d - \frac{\alpha}{n} \big)^2   + 
     a_1 - \frac{2\alpha}{n} b + \frac{\alpha^2}{n^2} c
  \bigg) \\ 
B_7 & = & \bigg(
  \bigg[
  2
  \bigg(
h   L \big(d - \frac{\alpha}{n} \big)^2   + 
     a_1 - \frac{2\alpha}{n} b + \frac{\alpha^2}{n^2} c
  \bigg)
  +2 \bigg(  h ( d - \frac{\alpha}{n} )  (1 - \frac{1}{n} )  -   h (1 - \delta) d \bigg)  
  \bigg] 
  \bigg) \\
B_8 & = &  \bigg[
  \frac{1}{n}  \big(
  a_1 + a_2 - 2\frac{\alpha}{n} b + \frac{\alpha^2}{n^2} c 
  \big)  + \frac{L}{n} h \big(d - \frac{\alpha}{n} \big)^2 
  \bigg] \\
B_9 & = &  c \delta .
\EEAS

In this expression, $y^{k-1}  - y^\ast$ only appears through a quadratic form, with a quadratic term
$ -  (  y^{k-1} - y^\ast)^\top \bigg[
    B_3 ( \idm - \frac{1}{n} ee^\top ) + B_4\frac{1}{n}  ee^\top 
    \bigg]  ( y^{k-1} - y^\ast)  $. If $B_3 \geqslant 0$ and $B_4 \geqslant 0$, then this quadratic term is non-positive and we may maximize in closed form with respect to $y^{k-1}$ to obtain an upper bound. Assuming that $B_3\neq 0$ and $B_4 \neq 0$, we use Lemma~\ref{lemma2},
   $\bigg[
    B_3 ( \idm - \frac{1}{n} ee^\top ) + B_4\frac{1}{n}  ee^\top 
    \bigg]^{-1} = \bigg[
    B_3^{-1} ( \idm - \frac{1}{n} ee^\top ) + B_4^{-1} \frac{1}{n}  ee^\top 
    \bigg] $, to obtain:
  \BEAS
   & & \E \big( \mathcal{L}(\theta^k) | \mathcal{F}_{k-1}  \big) -  ( 1 - \delta ) \mathcal{L}(\theta^{k-1})
  \\
   & \leqslant & B_0 \big( g(x^{k-1}) -  g(x^\ast)   \big) + B_9 \| x^{k-1} - x^\ast \|^2  
  + (x^{k-1} - x^\ast)^\top g'(x^{k-1})
  B_1 - B_2 \| g'(x^{k-1}) \|^2
\\
    & &  
    + \frac{1}{4} \frac{B_6^2}{B_3}  (  f'(x^{k-1}) - f'(x^\ast) ) ^\top ( \idm - \frac{1}{n} ee^\top)  (  f'(x^{k-1}) - f'(x^\ast) ) 
    \\
    & &     + \frac{1}{4} \frac{1}{B_4} 
    \bigg[
  B_5  e  ( x^{k-1} - x^\ast) + B_6 (  f'(x^{k-1}) - f'(x^\ast) ) + B_7 e g'(x^{k-1} )
  \bigg] ^\top ( \frac{1}{n} ee^\top ) 
  \\& & \hspace*{3cm}
  \bigg[
  B_5  e  ( x^{k-1} - x^\ast) + B_6 (  f'(x^{k-1}) - f'(x^\ast) ) + B_7 e g'(x^{k-1} )
  \bigg] \\
   & & +B_8 \| f'(x^{k-1}) - f'(x^\ast)\|^2 
\\
& = & B_0 \big( g(x^{k-1}) -  g(x^\ast)   \big) + B_9 \| x^{k-1} - x^\ast \|^2  
  + (x^{k-1} - x^\ast)^\top g'(x^{k-1})
  B_1 - B_2 \| g'(x^{k-1}) \|^2
\\
    & &  
    + \frac{1}{4} \frac{B_6^2}{B_3}  (  f'(x^{k-1}) - f'(x^\ast) ) ^\top ( \idm - \frac{1}{n} ee^\top)  (  f'(x^{k-1}) - f'(x^\ast) ) 
    \\
     & &    {  + \frac{n}{4} \frac{1}{B_4  } 
    \big\| 
  B_5    ( x^{k-1} - x^\ast) +   \big( B_7 + B_6 \big)  g'(x^{k-1} ) \big\|^2}
\\
   & & +B_8 \| f'(x^{k-1}) - f'(x^\ast)\|^2 .
   \EEAS
 By using Lipschitz continuity of each $f_i'$ to observe that
 $\|  f'(x^{k-1}) - f'(x^\ast)  \|^2 = \sum_{i=1}^n
\|  f_i'(x^{k-1}) - f_i'(x^\ast)  \|^2 \leqslant L \sum_{i=1}^n
( f_i'(x^{k-1}) - f_i'(x^\ast) )^\top ( x^{k-1} - x^\ast)   = nL g'(x^{k-1})^\top ( x^{k-1} - x^\ast)$, we obtain:

  \BEAS
  & & \E \big( \mathcal{L}(\theta^k) | \mathcal{F}_{k-1}  \big) -  ( 1 - \delta ) \mathcal{L}(\theta^{k-1})
  \\
  & \leqslant & 
   B_0  \bigg[ (x^{k-1}-x^\ast)^\top g'(x_{k-1}) - \frac{\mu}{2} \| x^{k-1} - x^\ast\|^2  \bigg] + B_9 \| x^{k-1} - x^\ast \|^2  
  \\
& & + (x^{k-1} - x^\ast)^\top g'(x^{k-1})
  B_1 - B_2 \| g'(x^{k-1}) \|^2
\\
    & &  
    +\bigg( \frac{nL}{4} \frac{B_6^2}{B_3} + nL B_8 \bigg) (x^{k-1} - x^\ast)^\top g'(x^{k-1})
    - \frac{n}{4} \frac{B_6^2}{B_3} \| g'(x^{k-1}) \|^2
    \\
     & &    + \frac{n}{4} \frac{B_5^2}{B_4  } \|x^{k-1} - x^\ast\|^2
    + \frac{n}{4} \frac{(B_6 + B_7)^2}{B_4  } \| g'(x^{k-1})\|^2
    \\
  & & +    \frac{2n}{4} \frac{B_5 (B_6+B_7) }{B_4  }  (x^{k-1} - x^\ast)^\top g'(x^{k-1})
 \\
   & = & \| x^{k-1} - x^\ast\|^2 \bigg[
   - B_0 \frac{\mu}{2} + B_9   + \frac{n}{4} \frac{B_5^2}{B_4  }
   \bigg] \\
   & & + (x^{k-1} - x^\ast)^\top g'(x^{k-1}) 
   \bigg[
   B_0 + B_1   + \frac{nL}{4} \frac{B_6^2}{B_3}   +    \frac{2n}{4} \frac{B_5 (B_6+B_7) }{B_4  } + n L B_8 \bigg]\\
   & & +  \| g'(x^{k-1})\|^2
   \bigg[
   -B_2  - \frac{n}{4} \frac{B_6^2}{B_3}   + \frac{n}{4} \frac{(B_6 + B_7)^2}{B_4  }
   \bigg] \\
   & = &  C_0 \| x^{k-1} - x^\ast\|^2 + C_1 (x^{k-1} - x^\ast)^\top g'(x^{k-1}) 
   +  C_2 \| g'(x^{k-1})\|^2,
  \EEAS
with
\BEAS
C_0 & =  &    - B_0 \frac{\mu}{2} + B_9   + \frac{n}{4} \frac{B_5^2}{B_4  }  \\
C_1 & =  & 
   B_0 + B_1   + \frac{nL}{4} \frac{B_6^2}{B_3}   +    \frac{2n}{4} \frac{B_5 (B_6+B_7) }{B_4  } + n L B_8 \\
C_2 & =  &    -B_2  - \frac{n}{4} \frac{B_6^2}{B_3}   + \frac{n}{4} \frac{(B_6 + B_7)^2}{B_4  }.
\EEAS

In order to show the decrease of the Lyapunov function, we need to show that
$ C_0 \| x^{k-1} - x^\ast\|^2 + C_1 (x^{k-1} - x^\ast)^\top g'(x^{k-1}) 
   +  C_2 \| g'(x^{k-1})\|^2 \leqslant 0$ for all $x^{k-1}$.
   
   If we assume that  $C_1\leqslant 0$ and $C_2 \leqslant 0$, then
   we have by strong-convexity
   {\smaller$$
    C_0 \| x^{k-1} - x^\ast\|^2 + C_1 (x^{k-1} - x^\ast)^\top g'(x^{k-1}) 
   +  C_2 \| g'(x^{k-1})\|^2 \leqslant 
    \| x^{k-1} - x^\ast\|^2 \big(
    C_0 + \mu C_1 + \mu^2 C_2
    \big),
   $$}
   and thus the condition is that
$$
C_0 + \mu C_1 + \mu^2 C_2 \leqslant 0.
$$
  
  And if we want to show that $\E \big( \mathcal{L}(\theta^k) | \mathcal{F}_{k-1}  \big) -  ( 1 - \delta ) \mathcal{L}(\theta^{k-1}) \leqslant - C_3 (x^{k-1} - x^\ast)^\top g'(x^{k-1}) $, then it suffices to have $C_1 + C_3 \leqslant 0$ and $C_0 + \mu ( C_1 + C_3) + \mu^2 C_2 \leqslant 0$.

 \subsection{Domination of $g(x^{k}) - g(x^\ast)$}
  
  By again using the strong convexity of $g$ (as in~\eqref{eq:11}), we have
  \BEAS
 & &  \mathcal{L}(\theta^k) - \gamma \big[
  g(x^k) - g(x^\ast)
  \big] \\
  & \geqslant & ( 2 h  -  \gamma)  \big[
  g(x^k) - g(x^\ast)  \big] + c \| x^k - x^\ast\|^2 \\
  & & + (y^k - y^\ast)^\top \bigg[ a_2 ( \idm - \frac{1}{n} ee^\top ) + \frac{1}{n} ee^\top
  ( na_1 + a_2 + n \mu d^2 ) \bigg] (y^k - y^\ast) \\
  & & +  (y^k - y^\ast)^\top e \bigg[
   2 d g'(x^k)  + 2b (x^k - x^\ast)
  \bigg]
\\
& \geqslant & ( 2h  -  \gamma)  \big[
  g(x^k) - g(x^\ast)  \big] + c \| x^k - x^\ast\|^2 \\
  & & - \frac{1}{4} \frac{n}{na_1 + a_2 + n \mu d^2} \| 
   2 d g'(x^k)  + 2b (x^k - x^\ast)\|^2 \mbox{ by minimizing with respect to } y^{k},
\\
& \geqslant &\| x^k - x^\ast\|^2 \bigg[
\frac{L}{2}  ( 2 h -  \gamma)   + c 
-   \frac{n}{na_1 + a_2 + n \mu d^2}
\big( dL + b
\big)^2
 \bigg]  ,
  \EEAS
  using the smoothness of $g$ and the assumption that $ 2h - \gamma \leqslant 0$.
  
  Thus, in order for $\mathcal{L}(\theta^k)$ to dominate $g(x^k) - g(x^*)$ we require the additional constraints
    $$
  \frac{L}{2}  ( 2h -  \gamma)   + c 
-   \frac{n}{na_1 + a_2 + n \mu d^2}
\big( dL + b
\big)^2 \geqslant 0
  $$
$$  na_1 + a_2 + n \mu d^2 \geqslant 0$$
 $$ 2h - \gamma \leqslant 0.$$

  \subsection{Finding constants}
  
  Given the algorithm and result parameters $(\alpha,\gamma,\delta)$, we have the following constraints:
    \begin{align*}
h & \geqslant  0  & \text{(Constraint 1)}\\
2h - \gamma & \leqslant  0 & \text{(Constraint 2)}  \\
  B_3 & > 0  & \text{(Constraint 3)} \\
  B_4 & >  0  & \text{(Constraint 4)} \\
  na_1 + a_2 + n \mu d^2 & \geqslant 0  & \text{(Constraint 5)} \\
  \frac{L}{2}  ( 2h -  \gamma)   + c 
-   \frac{n}{na_1 + a_2 + n \mu d^2}
\big( dL + b
\big)^2  &  \geqslant  0 & \text{(Constraint 6)}  \\
  C_2 & \leqslant   0  & \text{(Constraint  7)}  \\
  C_1 + C_3 & \leqslant   0 & \text{(Constraint 8)}  \\
  C_0 + \mu ( C_1 + C_3) + \mu^2 C_2 & \leqslant  0 & \text{(Constraint 9)} 
  \end{align*}
We also require  $C_1 \leqslant 0$, but this will follow from Constraint 8 since we will have $C_3 \geqslant 0$.
  All of the constraints above are convex in $a_1$, $a_2$, $b$, $c$, $d$, $h$. Thus, given candidate values for the remaining values, the feasibility of these constraints may be checked using a numerical toolbox (as a second-order cone program). However, these parameters should be  functions of $n$ and should be valid for all $\mu$. Thus, given a sampling of values of $\mu$ and $n$, representing these parameters as polynomials in $1/n$, the candidate functions may be found through a second-order cone programs.
  
  By experimenting with this strategy, we have come up with the following values for the constants:
  \BEAS
  a_1 & = &  \frac{1}{32nL} \big( 1 - \frac{1}{2n} \big) \\
  a_2 & = &  \frac{1}{16nL} \big( 1 - \frac{1}{2n} \big)\\
  b & = &  - \frac{1}{4n} \big( 1 - \frac{1}{n} \big) \\
  c & = & \frac{4L}{n} \\
  h & = &   \frac{1}{2} - \frac{1}{n} \\
d & = & \frac{\alpha}{n} \\
  \alpha &  = & \frac{1}{16L} \\
  \delta & = &   \min \big( \frac{1}{8n}, \frac{\mu}{16L} \big) \\
  \gamma & = &   1 \\
  C_3 &  = & \frac{1}{32 n}.
  \EEAS
We will assume $n>1$, since the result for $n=1$ follows because SAG is equivalent to gradient descent in this case.
Under this assumption, we see that Constraints 1 and 2 are satisfied by the above parameterizaiton. 
 Below we verify that Constraints 3-9 are also satisfied using symbolic computations in Matlab (the code used in this section is available on the first author's webpage). Specifically, we parameterize the constraints as polynomials and then verify that the polynomials are positive over an appropriate interval using the function below (which simply checks that no roots of the polynomial $P$ lie in the interval $(x_1,x_2)$, and that $P$ is positive at the mid-point of the interval).

\includegraphics[width=.8\textwidth]{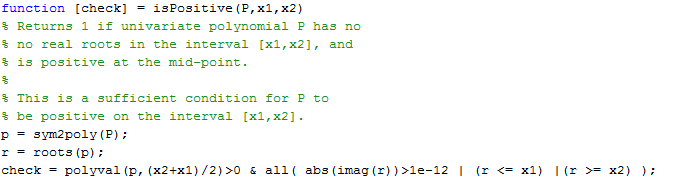}    
           
      \subsection{Verifying the result}
      
      To verify the result under these constants, we consider whether $\delta = 1/8n$ or $\mu/16L$.  In $B_4$, we discard the term of the form $(1-\delta) \mu h d^2$, which does not impact the validity of our result.   
      Moreover, without loss of generality, we may assume $L=1$ and $\mu \in [0,1]$.
      
      \subsubsection{Well-conditioned problems ($\mu \geqslant 2/ n$)}
      In this situation, it suffices to show the result for $\mu = 2/n$ since, (a) if a function is $\mu$-strongly convex, then it is $\mu'$-strongly convex for all smaller $\mu'$, and (b) the final inequality does not involve $\mu$. All constraints (when properly multiplied where appropriate by $B_3 B_4$) can be written as rational functions of $x = 1/n$:
      \begin{align*}
& B_3 & \text{(Constraint 3)}\\
& B_4  & \text{(Constraint 4)}\\
 & na_1 + a_2 + n \mu d^2 & \text{(Constraint 5)}\\
 & \big[  \frac{L}{2}  ( 2h -  \gamma)   + c \big] (na_1 + a_2 + n \mu d^2)
-   n
\big( dL + b
\big)^2  & \text{(Constraint 6)}\\
& -C_2 B_3 B_4  & \text{(Constraint 7)}\\
& -(C_1+C_3) B_3 B_4  & \text{(Constraint 8)}\\
& -  B_4 B_3 ( C_0+ (C_1+C_3)   \mu + C_2  \mu^2)  & \text{(Constraint 9)}
      \end{align*}
We only need to check the positivity of these polynomials in $x$. We do this using the function above, via the script below. Note that we also have $B_3 > 0$ and $B_4 > 0$ for $n>1$.

\includegraphics[width=.8\textwidth]{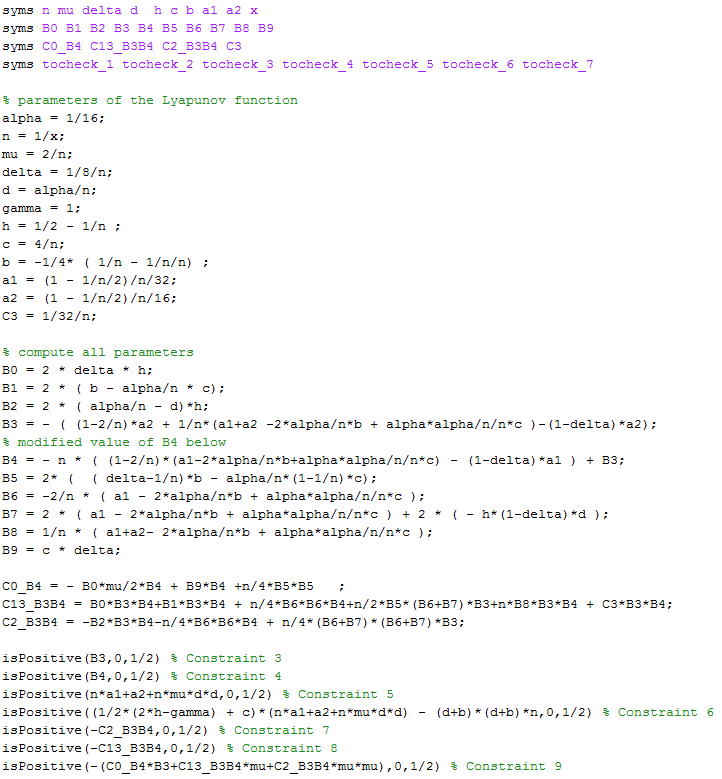}    
      
      \subsubsection{Ill-conditioned problems ($\mu \leqslant 2 / n$)}
      
      We consider the variables $x =1/n \in [0,1/2]$ and $y = n \mu / 2 \in [0,1]$, so that $\mu = 2y/n$. We may express all quantities using $x$ and $y$. The difficulty here is that we have two variables.
We first check the dependency in terms of $y$ of the expressions $B_3$, $B_4$, and $B_6+B_7$. These are univariate polynomials with an affine dependence on $y$, so their positivity can be checked by checking positivity for $y=1$ or $y=0$. Again making use of symbolic computation, we can deduce that
\begin{align*}
B_3 \text{is non-negative and decreasing in $y$,} & & \text{(Constraint 3)}\\
B_4 \text{is non-negative and decreasing in $y$,} & & \text{(Constraint 4)}\\
B_6+B_7 \text{is non-negative and increasing in $y$,}
\end{align*}
We include below the additional commands to verify these properties:

\includegraphics[width=.8\textwidth]{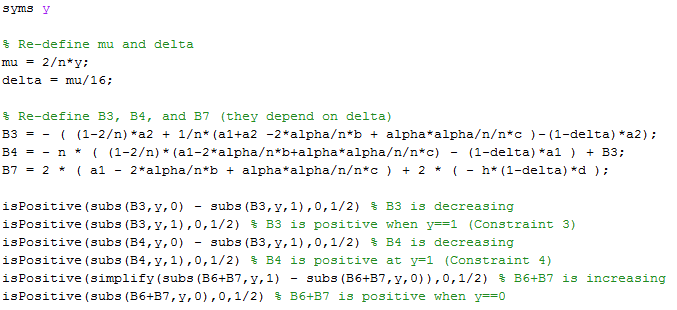}    \\
      Given the monotonicity of the bounds in $B_3$ and $B_4$, we only need to check our results for the smaller values of $B_3$ and $B_4$, i.e.,  for $y=1$. Similarly, because of the monotonicity in the term $(B_6+B_7)^2$, we may replace 
       $(B_6 + B_7)^2$ by $(B_6 + B_7)$ times its upper-bound (i.e., its value at $y=1$). Also note that $B_5$ is divisble by $y$, and we have $B_3 > 0$ and $B_4 > 0$ for $n > 1$. Using this, we can show that Constraints 5-9 are satisfied by checking the positivity of the following polynomials:
      \begin{align*}
 & na_1 + a_2 + n \mu d^2 & \text{(Constraint 5)}\\
 & \big[  \frac{L}{2}  ( 2h -  \gamma)   + c \big] (na_1 + a_2 + n \mu d^2)
-   n
\big( dL + b
\big)^2  & \text{(Constraint 6)}\\
& -C_2 B_3 B_4  & \text{(Constraint 7)}\\
& -(C_1+C_3) B_3 B_4  & \text{(Constraint 8)}\\
&  -B_4B_3 ( C_0 /\mu+ (C_1+C_3)    + C_2 \mu)  & \text{(Constraint 9)}
      \end{align*}
Constraints 5-7 are affine in $y$ and we thus only need to check positivity for $y=0$ and $y=1$. 

\includegraphics[width=.8\textwidth]{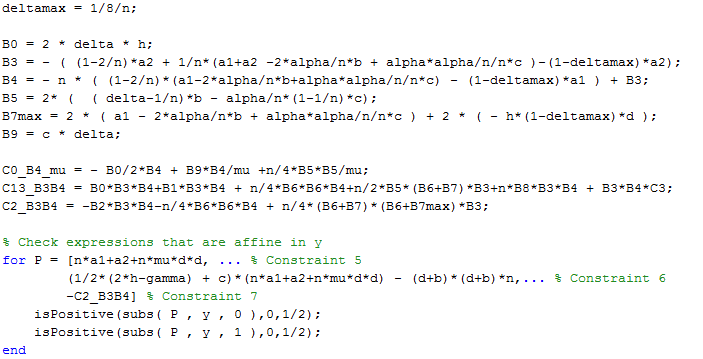}   \\
The other two expressions (Constraints 8 and 9) are more complicated as there are second-order polynomials in $y$. If $n \geqslant 5$, they  have positive second derivatives, negative derivatives at $y=0$ and $y=1$, and positive values at $y=1$. They are thus positive, which is verified by the code below.

\includegraphics[width=.75\textwidth]{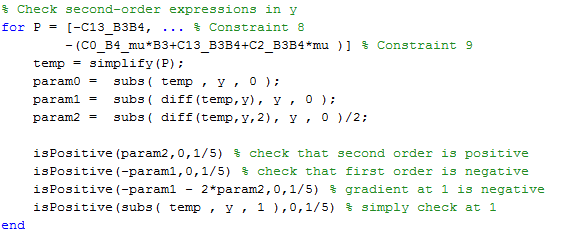}   \\
 For the remaining scenarios where $n \in \{2,3,4\}$, we have a fixed $x$ so we can check the positivity of the polynomials in $y$. A Matlab script that does the steps above in addition to checking these not-particularly-interesting cases is available from the first author's web page.

\subsection{Convergence Rate}

We have that
\[
g(x^k) - g(x^\ast) \leqslant \mathcal{L}(\theta^k),
\]
and that
\[
\mathbb{E}(\mathcal{L}(\theta^k) | \mathcal{F}_{k-1}) - (1-\delta)\mathcal{L}(\theta^{k-1}) \leqslant -\frac{1}{32n}(x^{k-1} - x^\ast)^\top g'(x^{k-1}) \leqslant 0.
\]
In the strongly-convex case ($\delta > 0$), combining these gives us the linear convergence rate
\[
\mathbb{E} (g(x^k)) - g(x^\ast) \leqslant (1-\delta)^k\mathcal{L}(\theta^0).
\]
In the convex case ($\mu=0$ and $\delta = 0$), we have by convexity that
\[
-\frac{1}{32n}(x^{k-1} - x^\ast)^\top g'(x^{k-1}) \leqslant \frac{1}{32n}[g(x^\ast) - g(x^{k-1})],
\]
We then have
\[
\frac{1}{32n}[g(x^{k-1}) - g(x^\ast)] \leqslant L(\theta^{k-1}) - \mathbb{E}(L(\theta^k)|\mathcal{F}_{k-1}).
\]
Summing this up to iteration $k$ yields
\[
\frac{1}{32n}\sum_{i=1}^k [\mathbb{E}(g(x^{i-1})) - g(x^\ast)] \leqslant \sum_{i=1}^k \mathbb{E}[L(\theta^{i-1}) - L(\theta^i)] = L(\theta^0) - \mathbb{E}[L(\theta^k)] \leqslant L(\theta^0).
\]
Finally, by Jensen's inequality note that
\[
\mathbb{E}\left[g\left(\frac{1}{k}\sum_{i=0}^{k-1}x^i\right)\right] \leq \frac{1}{k}\sum_{i=0}^{k-1}\mathbb{E}[g(x^i)],
\]
so with $\bar{x}^k=\frac{1}{k}\sum_{i=0}^{k-1}x^i$ we have
\[
\mathbb{E}[g(\bar{x}^k)] - g(x^\ast) \leqslant \frac{32n}{k}L(\theta^0).
\]

\subsection{Intial values of Lyapunov function}

To obtain the results of Theorem~\ref{thm}, all that remains is computing the initial value of the Lyapunov function for the two initializations.

\subsubsection{Initialization with the zero vector}  
  
  If we initialize with $y^0 = 0$ we have
  \begin{align*}
\mathcal{L}(\theta^0) = 2h(g( x^0) - g(x^\ast)) & + 
f'(x^\ast)^\top (a_1ee^\top + a_2I)f'(x^\ast)\\
& + 2bf'(x^\ast)^\top e (x^0 - x^\ast) + 
c(x^0 - x^*)^\top(x^0 - x^\ast).
\end{align*}
Plugging in the values of our parameters,
\begin{eqnarray*}
h = \frac{1}{2} - \frac{1}{n},
\qquad
a_1 & = \frac{1}{32nL}\left(1 - \frac{1}{2n}\right),
\qquad
a_2 & = \frac{1}{16nL}\left(1 - \frac{1}{2n}\right),\\
b & = -\frac{1}{4n}\left(1 - \frac{1}{n}\right),
\qquad
c & = \frac{4L}{n},
\end{eqnarray*}
and using $\sigma^2 = \frac{1}{n} \sum_i \|f_i'(x^\ast)\|^2$ 
we obtain (noting that $e^\top f'(x^\ast) = 0$) 
\begin{align*}
\mathcal{L}(\theta^0)  & = \left(1 - \frac{2}{n}\right)(g( x^0) - g(x^\ast)) + 
\left(1 - \frac{1}{2n}\right)f'(x^\ast)^\top (\frac{1}{32nL}ee^\top + \frac{1}{16nL}I)f'(x^\ast) 
\\ 
& \hspace{3cm} - \frac{1}{2n}\left(1 - \frac{1}{n}\right)f'(x^\ast)^\top e (x^0 - x^\ast) + 
\frac{4L}{n}(x^0 - x^\ast)^\top(x^0 - x^\ast) \\
& \leqslant g( x^0) - g(x^\ast) + 
\frac{\sigma^2}{16L}+ 
\frac{4L}{n}\norm{x^0 - x^\ast}^2
\end{align*}

\subsubsection{Initialization with average gradient}

If we initialize with $y_i^0 = y_i^0 = f_i'(x^0) - (1/n)\sum_i f_i'(x^0)$ we have, by noting that we still have $(y^0)^\top e = 0$,
  \begin{align*}
\mathcal{L}(\theta^0)
& = \left(1 - \frac{2}{n}\right)(g( x^0) - g(x^\ast)) + \frac{1}{16nL}\left(1 - \frac{2}{n}\right)\norm{y^0 - f'(x^\ast)}^2 + \frac{4L}{n}\norm{x^0 - x^*}^2.\\
  \end{align*}
  By observing that $y^0 = f'(x^0) - eg'(x^0)$ and using~\citethm{Equations 2.17}{nesterov2004introductory}, we can bound the norm in the second term as follows:
  \begin{align*}
  \norm{y^0 - f'(x^\ast)}^2 &= \norm{(f'(x^0) - f'(x^\ast)) - e(g'(x^0) - g'(x^\ast)}^2\\
  & \leqslant  2\norm{f'(x^0) - f'(x^\ast)}^2 + 2\norm{e(g'(x^0) - g'(x^\ast))}^2
    \\
    & = 2\sum_{i=1}^n\norm{f_i'(x^0) - f_i'(x^\ast)}^2 + 2n\norm{g'(x^0) - g'(x^\ast)}^2
    \\
   & \leqslant 4L\sum_{i=1}^n [f_i(x^0) - f_i(x^\ast) - f_i'(x^\ast)^\top(x_0-x^\ast)] + 4nL(g(x^0) - g(x^\ast))
   \\
   & = 8nL(g(x^0) - g(x^\ast)).
   \end{align*} 
    Using this in the Lyapunov function we obtain
    \begin{align*}
    L(\theta^0) \leqslant \frac{3}{2}(g( x^0) - g(x^\ast)) +\frac{4L}{n}\norm{x^0 - x^*}^2.
    \end{align*}
}



\iftoggle{springer}
{

}
{
\bibliography{bib}   
\bibliographystyle{plainnat}
}

\end{document}